\documentclass[11pt]{article}
\usepackage{latexsym,bm}
\usepackage{color}
\usepackage{amssymb, graphicx, amsmath}  

\newcommand\qed{\qquad $\square$}

\def \[{\begin{equation}}
\def \]{\end{equation}}

\newtheorem{theorem}{Theorem}[section]

\newtheorem{definition}{Definition}[section]

\newtheorem{lemma}{Lemma}[section]

\newtheorem{corollary}{Corollary }[section]

\begin{document}

\begin{center}
{\bf Dirichlet problem of nonlinear second order partial differential equations resolved with any derivatives}

\medskip

  {  Jianfeng Wang }

 Department of Mathematics, Hohai University, Nanjing, 210098,

  P.R. China, email: wjf19702014@163.com.

\end{center}
\bigskip
{\narrower \noindent  {\bf Abstract.} In this paper we will discuss the Dirichlet problem of nonlinear second order partial differential equations resolved with any derivatives. First, we transform it into generalized integral equations. Next, we discuss the existence of the classical solution by Leray-Schauder degree and Sobolev space\ $H^{-m_{1}}(\Omega_{1})$. 
    \par }

\vskip 1.0 true cm
 \section{Introduction}
In this paper, we will discuss the Dirichlet problem of nonlinear second order partial differential equations resolved with any derivatives as follows,
\[v_{j}=f_{j}(v_{m+1},\ v_{m+2},\cdots,\ v_{10m},\ x,\ y,\ z),\ 1\leq j\leq m,\] where\ $v_{1},\ v_{2},\cdots,\ v_{10m}$\ is a permutation of the components of\ $u,\ \partial u,\ \partial^{2}u$,\ $\partial u=(u_{x},\ u_{y},\ u_{z})^{T}$,\ $\partial^{2}u=(u_{xx},\ u_{xy},\ u_{xz},\ u_{yy},\ u_{yz},\ u_{zz})^{T}$,\ $u=(u_{1}(x,\ y,\ z),\cdots,u_{m}(x,\ y,\ z))^{T}\in C^{2}(\overline{\Omega})$,\ $f_{j},\ 1\leq j\leq m,$\ are continuous functions,\ $m$\ is a positive whole number,\ $(x,\ y,\ z)^{T}\in \overline{\Omega}\subset R^{3},\ \overline{\Omega}=\Omega\cup\partial\Omega$,\ $\Omega$\ is a bounded domain, and\ $\partial\Omega\in C^{1,\ \alpha},\ 0<\alpha\leq 1$, moreover,\ $\overline{\Omega}$\ is convex, the boundary conditions should be known as follows, Dirichlet problem,\[u|_{\partial\Omega}\in C^{2}(\partial\Omega).\]
We have assumed\ $u\in C^{2}(\Omega)$. However, that is not enough to ensure\[f_{j}(v_{m+1},\ v_{m+2},\cdots,\ v_{10m},\ x,\ y,\ z)\in L^{2}(\Omega),\ 1\leq j\leq m,\]which we will use.\\
In Section 3, we will transform Eq(1.1) into the equivalent generalized integral equations as follows,\[ Z_{1}=w_{1}(x,\ y,\ z)+w_{2}(x,\ y,\ z).\ast(\psi(Z_{1})),\]
where\ $Z_{1}=(v_{m+1},\ v_{m+2},\cdots,\ v_{10m})^{T},\ .\ast$\ means the matrix convolution as follows,
\[ w_{2}.\ast(\psi(Z_{1}))=\int_{R^{3}}w_{2}(x-x_{1},\ y-y_{1},\ z-z_{1})(\psi(Z_{1}(x_{1},\ y_{1},\ z_{1})))dx_{1}dy_{1}dz_{1}. \]Only a few of these can indeed be transformed into the integral equations, this is why we call them the generalized integral equations. Here generalized comes from\ $w_{1}$\ being related to the Dirac function and\ $w_{2}$\ not being integrable on\ $\overline{\Omega}$.\ The term equivalent is defined in the following.
\begin{definition}\label{definition}The equations\ $f_{1}(x)=0$\ and\ $f_{2}(y)=0$\ are equivalent, if and only if there exist continuous mappings\ $T_{1},\ T_{2}$, such that\ $\forall x,\ y$\ are respectively the solutions of\\ $f_{1}(x)=0,\ f_{2}(y)=0$,\ we have\ $f_{1}(T_{2}(y))=0,\ f_{2}(T_{1}(x))=0,$\ and\ $T_{2}(T_{1}(x))=x$\ ,\ $T_{1}(T_{2}(y))=y$.\end{definition}
In Section 4, we will discuss the existence of the classical solution of Eqs(1.1). We will use the theory on\ $H^{-m_{1}}(\Omega)$, which is defined on page 130 in [10], and a primary theory on the Leray-Schauder degree. You will see that the classical solution of Eqs(1.1) exists in many cases.\\
\section{Boundary of\ $C^{1,\ \alpha}$}\setcounter{equation}{0}
In this section, we will explain why we chose\ $\partial\Omega\in C^{1,\ \alpha},\ 0<\alpha\leq 1,$\ instead of\ $\partial\Omega\in C^{2,\ \alpha},\ 0<\alpha\leq 1,$\ or the Sobolev's imbedding surface, the Lyapunov's surface, the Hopf's surface for the strong maximum principle.\\First of all, we know\ $C^{1,\ \alpha}\supset C^{2,\ \alpha},\ 0<\alpha\leq 1$.\ Simplicity and generality are our eternal pursuit. Secondly,\ $\partial\Omega\in C^{1,\ \alpha},\ 0<\alpha\leq 1$\ is just the Sobolev's imbedding surface from the following.
\begin{theorem} \label{Theorem2-1} If\ $\Omega$\ is bounded,\ $\partial\Omega\in C^{1,\ \alpha},\ 0<\alpha\leq 1,$\ then domain\ $\Omega$\ satisfies a uniform exterior cone condition, that is, there exists a fixed finite right circular cone\ $K$,\ such that each\ $P=(x,\ y,\ z)^{T}\in \partial\Omega$\ is the vertex of a cone\ $K(P)$,\ $\overline{K(P)}\cap\overline{\Omega}=P$, and\ $K(P)$\ is congruent to\ $K$. \end{theorem}
{\it Proof of theorem 2.1}. By using the equivalent definition of\ $C^{1,\ \alpha},\ 0<\alpha\leq 1,$\ on page 94 in [2], if\ $\partial\Omega\in C^{1,\ \alpha}$,\ $0<\alpha\leq 1,$\ then each\ $P_{0}=(x_{0},\ y_{0},\ z_{0})^{T}\in \partial\Omega$,\ there exists a neighborhood\ $U(P_{0},\ \delta_{0}(P_{0})),\ \delta_{0}(P_{0})>0$, where
 $$U(P_{0},\ \delta_{0}(P_{0}))=\{P:\ |P-P_{0}|<\delta_{0}(P_{0})\},\ \overline{U}(P_{0},\ \delta_{0}(P_{0}))=\{P:\ |P-P_{0}|\leq\delta_{0}(P_{0})\},$$
$P=(x,\ y,\ z)^{T}\in R^{3},\ \mid P-P_{0}\mid=\sqrt{(x-x_{0})^{2}+(y-y_{0})^{2}+(z-z_{0})^{2}}$, and\ $U(P_{0},\ \delta_{0}(P_{0}))\cap\partial\Omega$\ is a graph of a\ $C^{1,\ \alpha},\ 0<\alpha\leq 1,$\ function of two of the coordinates\ $x,\ y,\ z$.\\ Without loss of the generality, we assume such a function is\ $f_{0}(x,\ y)\in C^{1,\ \alpha},\ 0<\alpha\leq 1$, and we have the following,\ $z-f_{0}(x,\ y)=0$,\ if\ $P=(x,\ y,\ z)^{T}\in U(P_{0},\ \delta_{0}(P_{0}))\cap\partial\Omega$,\ $z-f_{0}(x,\ y)>0$,\ if\ $P=(x,\ y,\ z)^{T}\in U(P_{0},\ \delta_{0}(P_{0}))\setminus\overline{\Omega}$, and\ $z-f_{0}(x,\ y)<0$,\ if\ $P=(x,\ y,\ z)^{T}\in U(P_{0},\ \delta_{0}(P_{0}))\cap\Omega$.\\ We will get the same results if\ $U(P_{0},\ \delta_{0}(P_{0}))\cap\partial\Omega$\ is a graph of other functions.\\ We introduce the following lemma.
\begin{lemma} \label{lemma1} If\ $\Omega$\ is bounded,\ $\partial\Omega\in C^{1,\ \alpha},\ 0<\alpha\leq 1,$\ then for each\ $P_{0}=(x_{0},\ y_{0},\ z_{0})^{T}\in \partial\Omega$,\ there exists a\ $C(P_{0})>0$, related to\ $f_{0}(x,\ y)$, such that\[\mid {\bf r}_{P_{0}P}\cdot n_{p_{0}}\mid\leq C(P_{0})\mid P-P_{0}\mid^{1+\alpha},\ \forall P=(x,\ y,\ z)^{T}\in U(P_{0},\ \delta_{0}(P_{0}))\cap\partial\Omega,\]where\ ${\bf r}_{P_{0}P}=\overrightarrow{P_{0}P},\ n_{p_{0}}$\ is exterior normal vector to\ $\partial\Omega$\ at point\ $P_{0}$.\end{lemma}
{\it Proof of lemma 2.1}. From\ $n_{p_{0}}=\cfrac{1}{\sqrt{f_{0x}^{2}(x_{0},\ y_{0})+f_{0y}^{2}(x_{0},\ y_{0})+1}}\ (-f_{0x}(x_{0},\ y_{0}),\ -f_{0y}(x_{0},\ y_{0}),\ 1)^{T}$, where\ $f_{0x},\ f_{0y}$\ are partial derivatives of\ $f_{0}$, we get\[{\bf r}_{P_{0}P}\cdot n_{p_{0}}=\cfrac{z-z_{0}-(x-x_{0})f_{0x}(x_{0},\ y_{0})-(y-y_{0})f_{0y}(x_{0},\ y_{0})}{\sqrt{f_{0x}^{2}(x_{0},\ y_{0})+f_{0y}^{2}(x_{0},\ y_{0})+1}}.\]From\ $z-z_{0}=f_{0}(x,\ y)-f_{0}(x_{0},\ y_{0})$, and\ $f_{0}(x,\ y)\in C^{1,\ \alpha}$, we have\[z-z_{0}=(x-x_{0})f_{0x}(tx_{0}+(1-t)x,\ ty_{0}+(1-t)y)+(y-y_{0})f_{0y}(tx_{0}+(1-t)x,\ ty_{0}+(1-t)y),\] where\ $0\leq t\leq 1$. From\ $f_{0}(x,\ y)\in C^{1,\ \alpha}$, we know there exists\ $C_{\alpha}(P_{0})>0$, related to\ $f_{0}(x,\ y)$, such that\ $\forall P_{1}=(x_{1},\ y_{1},\ z_{1})^{T},\ P_{2}=(x_{2},\ y_{2},\ z_{2})^{T}\in \overline{U}(P_{0},\ \delta_{0}(P_{0}))\cap\partial\Omega$,\begin{eqnarray*}&&\mid f_{0x}(x_{1},\ y_{1})-f_{0x}(x_{2},\ y_{2})\mid \leq C_{\alpha}(P_{0})\mid (P_{1}-P_{2})_{1}\mid^{\alpha},\\ &&\mid f_{0y}(x_{1},\ y_{1})-f_{0y}(x_{2},\ y_{2})\mid \leq C_{\alpha}(P_{0})\mid (P_{1}-P_{2})_{1}\mid^{\alpha},\end{eqnarray*} where\
$(P_{1}-P_{2})_{1}=(x_{1}-x_{2},\ y_{1}-y_{2},\ 0)^{T}$.\\
Hence we obtain the following,\begin{eqnarray}&&\mid f_{0x}(tx_{0}+(1-t)x,\ ty_{0}+(1-t)y)-f_{0x}(x_{0},\ y_{0})\mid \leq C_{\alpha}(P_{0})\mid (P-P_{0})_{1}\mid^{\alpha},\\ &&\mid f_{0y}(tx_{0}+(1-t)x,\ ty_{0}+(1-t)y)-f_{0y}(x_{0},\ y_{0})\mid \leq C_{\alpha}(P_{0})\mid (P-P_{0})_{1}\mid^{\alpha},\end{eqnarray} where\
$(P-P_{0})_{1}=(x-x_{0},\ y-y_{0},\ 0)^{T}$.\\ From\ $\sqrt{f_{0x}^{2}(x_{0},\ y_{0})+f_{0y}^{2}(x_{0},\ y_{0})+1}\geq 1$, we obtain,
\[\mid{\bf r}_{P_{0}P}\cdot n_{p_{0}}\mid\leq 2C_{\alpha}(P_{0})\mid (P-P_{0})_{1}\mid^{1+\alpha}\leq C(P_{0})\mid P-P_{0}\mid^{1+\alpha},\] where\ $C(P_{0})=2C_{\alpha}(P_{0})$.
\qed\\
Now we can make an exterior finite right circular cone\ $K(P_{0})$\ at point\ $P_{0}$. We choose a\ $\delta_{1}(P_{0})<\delta_{0}(P_{0})$\ that is small enough such that\ $C(P_{0})(2\delta_{1}(P_{0}))^{\alpha}<1$.\ We let\ $P_{0}$\ be the vertex of a cone\ $K(P_{0})$\ and\ $n_{p_{0}}$\ be the symmetry axis. The polar angle is\ $\theta(P_{0})=\arccos(C(P_{0})(2\delta_{1}(P_{0}))^{\alpha})\in(0,\ \pi/2)$,\ and the length of generatrix is\ $\delta_{1}(P_{0})/3$. \\We will prove that\ $\overline{K(P_{0})}\cap\overline{\Omega}=P_{0}$. \\ From\ $\mid P-P_{0}\mid\leq \delta_{1}(P_{0})/3,\ \forall P\in \overline{K(P_{0})}$, we can obtain\ $\overline{K(P_{0})}\subset U(P_{0},\ \delta_{1}(P_{0}))$.\\ And from\[{\bf r}_{P_{0}P}\cdot n_{p_{0}}\geq\mid P-P_{0}\mid C(P_{0})(2\delta_{1}(P_{0}))^{\alpha}> C(P_{0})\mid P-P_{0}\mid^{1+\alpha},\ \forall P\in \overline{K(P_{0})},\ P\neq P_{0}, \]we can obtain\ $\overline{K(P_{0})}\cap \partial\Omega=P_{0}$.\\ Finally, if there exists\ $P\in \overline{K(P_{0})}\cap\Omega$,\ then we can obtain\ $z<f_{0}(x,\ y)$, where\ $P=(x,\ y,\ z)^{T}$. From (2.2), we can get\ ${\bf r}_{P_{0}P}\cdot n_{p_{0}}< C(P_{0})\mid P-P_{0}\mid^{1+\alpha}$. This contradicts (2.7). Hence\ $\overline{K(P_{0})}\cap\overline{\Omega}=P_{0}$.\\
Then, we can make a uniform exterior finite right circular cone\ $K$\ for each\ $P\in U(P_{0},\ \delta_{1}(P_{0})/3)\cap\partial\Omega$.\ $K$\ is congruent to\ $K(P_{0})$, if we let\ $P$\ be the vertex of a cone\ $K(P)$,\ $n_{p}$\ the symmetry axis, a polar angle of\ $\theta(P)=\theta(P_{0})=\arccos(C(P_{0})(2\delta_{1}(P_{0}))^{\alpha})$,\ and the length of generatrix\ $\delta_{1}(P_{0})/3$. From lemma 2.1, we obtain the following,\ $\forall P_{1}=(x_{1},\ y_{1},\ z_{1})^{T},\ P_{2}=(x_{2},\ y_{2},\ z_{2})^{T}\in \overline{U}(P_{0},\ \delta_{0}(P_{0}))\cap\partial\Omega$,
\[\mid{\bf r}_{P_{2}P_{1}}\cdot n_{p_{2}}\mid\leq C(P_{0})\mid P_{1}-P_{2}\mid^{1+\alpha}.\]
If there exists\ $P^{\prime}\neq P\in K(P)\cap \overline{\Omega}$, then it will contradict the following,
\[{\bf r}_{PP^{\prime}}\cdot n_{p}\geq\mid P^{\prime}-P\mid C(P_{0})(2\delta_{1}(P_{0}))^{\alpha}> C(P_{0})\mid P^{\prime}-P\mid^{1+\alpha},\ \forall P^{\prime}\in \overline{K(P)},\ P^{\prime}\neq P. \]
From arbitrary\ $P_{0}$, we obtain that\[\bigcup_{P_{0}\in \partial\Omega}[U(P_{0},\ \delta_{1}(P_{0})/3)\cap\partial\Omega]\] is an open cover for\ $\partial\Omega$. From the Heine-Borel theorem, we see that there exists a finite sub-cover for\ $\partial\Omega$ as follows,\ $\exists N>0,\ \exists P_{k}\in\partial\Omega,\ \exists \delta_{1}(P_{k})>0,\ 1\leq k\leq N$, such that
\[ \bigcup_{k=1}^{N}[U(P_{k},\ \delta_{1}(P_{k})/3)\cap\partial\Omega]\supset\partial\Omega,\] and the definitions of\ $\delta_{0}(P_{k}),\ C_{\alpha}(P_{k}),\ C(P_{k}),\ \delta_{1}(P_{k}),\ \theta(P_{k})$\ are the same as\ $\delta_{0}(P_{0}),\ C_{\alpha}(P_{0}),\ C(P_{0})$,\ $\delta_{1}(P_{0}),\ \theta(P_{0})$,\ $1\leq k\leq N$.\\ So we can make a uniform exterior finite right circular cone\ $K$\ for each\ $P\in \partial\Omega$.\ All the\ $K(P)$\ are congruent, if we let\ $P$\ be the vertex of a cone\ $K(P)$,\ $n_{p}$\ be the symmetry axis, have a polar angle of\ $\theta^{\ast}$,\ and length of generatrix of\ $\delta_{1}^{\ast}$, where
 $$ \theta^{\ast}=\min_{1\leq k\leq N}\theta(P_{k}),\ \delta_{1}^{\ast}=\min_{1\leq k\leq N}\delta_{1}(P_{k})/3.$$ Hence the statement stands. \qed
\\If we choose\ $-n_{p_{0}}$\ as the symmetry axis of the cone, then we can transform the uniform exterior cone into a uniform interior cone. Hence, domain\ $\Omega$\ satisfies a uniform exterior and interior cone condition, if\ $\Omega$\ is bounded,\ $\partial\Omega\in C^{1,\ \alpha},\ 0<\alpha\leq 1.$\ It looks like\ $C^{1,\ \alpha},\ 0<\alpha\leq 1,$\ is better than\ $C^{0,\ 1}$\ for the Sobolev's imbedding. However, we do not discuss the weak solution of Eqs(1.1) and (1.2) directly here. We will be discussing the classical solution directly. We want to obtain the equivalent equations that the classical solution should satisfy.\\ Thirdly, the following two theorems demonstrate that\ $\partial\Omega\in C^{1,\ \alpha},\ 0<\alpha\leq 1,$\ may be taken as the Lyapunov's surface. We will use them in section 3. We require two lemmas as follows.\begin{lemma} \label{lemma2.2} If\ $\Omega$\ is bounded,\ $\partial\Omega\in C^{1,\ \alpha},\ 0<\alpha\leq 1$,\ $f(M,\ P),\ M\neq P$\ is continuous, and\ $\forall M_{0}\in \partial\Omega,\ \forall \epsilon>0,\ \exists \delta>0$,\ such that\ $\forall M\in U(M_{0},\ \delta)$, we have
$$|\int_{U(M_{0},\ \delta)\cap\partial\Omega}f(M,\ P)dS_{P}|\leq \epsilon, $$ then \[\omega(M)=\int_{\partial\Omega}f(M,\ P)dS_{P},\ \forall M=(x,\ y,\ z)^{T}\in R^{3},\] will be continuous. In particular, if\ $\forall M_{0}\in \partial\Omega$,\ there exists a neighbourhood\ $U(M_{0})$,\ and\ $\delta_{1}\in(0,\ 1],\ C>0$, such that\[ |f(M,\ P)|\leq \cfrac{C}{r_{MP}^{2-\delta_{1}}}, \ \forall M\in U(M_{0}),\ \forall P\in  U(M_{0})\cap\partial\Omega ,\]where\ $r_{MP}=|M-P|$, then\ $\omega(M)$\ will also be continuous.\end{lemma}The proof is available on pages 178 to 180 of [6].
\begin{lemma} \label{lemma2.3} If\ $\Omega$\ is bounded,\ $\partial\Omega\in C^{1,\ \alpha},\ 0<\alpha\leq 1,$\ then there exist\ $\delta_{0}>0,\ C_{0}>0$, such that for each\ $P_{0}=(x_{0},\ y_{0},\ z_{0})^{T}\in \partial\Omega$,\ we have the following,\[\mid n_{p}-n_{p_{0}}\mid\leq C_{0}r_{PP_{0}}^{\alpha},\ \forall P=(x,\ y,\ z)^{T}\in U(P_{0},\ \delta_{0})\cap\partial\Omega.\]\end{lemma}
{\it Proof of lemma 2.3}. If\ $\Omega$\ is bounded,\ $\partial\Omega\in C^{1,\ \alpha},\ 0<\alpha\leq 1,$\ then from Theorem 2.1, for each\ $P_{0}=(x_{0},\ y_{0},\ z_{0})^{T}\in \partial\Omega$,\ there exist\ $P_{k}\in \partial\Omega,\ 1\leq k\leq N$, such that\ $P_{0}\in U(P_{k},\ \delta_{1}(P_{k})/3)\cap\partial\Omega$. And\ $U(P_{k},\ \delta_{1}(P_{k}))\cap\partial\Omega$\ is a graph of a\ $C^{1,\ \alpha},\ 0<\alpha\leq 1,$\ function of two of the coordinates\ $x,\ y,\ z$. \\Without loss of the generality, we assume such a function is\ $f_{k}(x,\ y)\in C^{1,\ \alpha},\ 0<\alpha\leq 1$, and we have the following,\ $z-f_{k}(x,\ y)=0$,\ if\ $P=(x,\ y,\ z)^{T}\in U(P_{k},\ \delta_{1}(P_{k}))\cap\partial\Omega$,\ $z-f_{k}(x,\ y)>0$,\ if\ $P=(x,\ y,\ z)^{T}\in U(P_{k},\ \delta_{1}(P_{k}))\setminus\overline{\Omega}$, and\ $z-f_{k}(x,\ y)<0$,\ if\ $P=(x,\ y,\ z)^{T}\in U(P_{k},\ \delta_{1}(P_{k}))\cap\Omega$.\\ We obtain the following results if\ $U(P_{k},\ \delta_{1}(P_{k}))\cap\partial\Omega$\ is a graph of other functions.\\ If we let\ $\delta_{0}=\delta_{1}^{\ast}$, then from\ $\delta_{1}^{\ast}\leq \delta_{1}(P_{k})/3$, we obtain the following,
$$U(P_{0},\ \delta_{0})\cap\partial\Omega\subset U(P_{k},\ \delta_{1}(P_{k}))\cap\partial\Omega. $$
Then\ $\forall P=(x,\ y,\ z)^{T}\in U(P_{0},\ \delta_{0})\cap\partial\Omega$, we obtain the following,
\[n_{p}=\cfrac{1}{\varphi(x,\ y)}\ (-f_{kx}(x,\ y),\ -f_{ky}(x,\ y),\ 1)^{T},\ n_{p_{0}}=\cfrac{1}{\varphi(x_{0},\ y_{0})}\ (-f_{kx}(x_{0},\ y_{0}),\ -f_{ky}(x_{0},\ y_{0}),\ 1)^{T},\] where\ $\varphi(x,\ y)=\sqrt{f_{kx}^{2}(x,\ y)+f_{ky}^{2}(x,\ y)+1}$,\ $f_{kx},\ f_{ky}$\ are partial derivatives of\ $f_{k}$.\\ We can deduce the following,
\begin{eqnarray*}\cfrac{f_{kx}(x,\ y)}{\varphi(x,\ y)}-\cfrac{f_{kx}(x_{0},\ y_{0})}{\varphi(x_{0},\ y_{0})}&=&\cfrac{\varphi(x_{0},\ y_{0})f_{kx}(x,\ y)-\varphi(x,\ y)f_{kx}(x_{0},\ y_{0})}{\varphi(x,\ y)\varphi(x_{0},\ y_{0})}\\&=&\cfrac{(\varphi(x_{0},\ y_{0})-\varphi(x,\ y))f_{kx}(x,\ y)+\varphi(x,\ y)(f_{kx}(x,\ y)-f_{kx}(x_{0},\ y_{0}))}{\varphi(x,\ y)\varphi(x_{0},\ y_{0})},\end{eqnarray*}
\begin{eqnarray*}\cfrac{f_{ky}(x,\ y)}{\varphi(x,\ y)}-\cfrac{f_{ky}(x_{0},\ y_{0})}{\varphi(x_{0},\ y_{0})}&=&\cfrac{\varphi(x_{0},\ y_{0})f_{ky}(x,\ y)-\varphi(x,\ y)f_{ky}(x_{0},\ y_{0})}{\varphi(x,\ y)\varphi(x_{0},\ y_{0})}\\&=&\cfrac{(\varphi(x_{0},\ y_{0})-\varphi(x,\ y))f_{ky}(x,\ y)+\varphi(x,\ y)(f_{ky}(x,\ y)-f_{ky}(x_{0},\ y_{0}))}{\varphi(x,\ y)\varphi(x_{0},\ y_{0})},\end{eqnarray*}
\begin{eqnarray*} \varphi(x_{0},\ y_{0})-\varphi(x,\ y)&=&\cfrac{\varphi^{2}(x_{0},\ y_{0})-\varphi^{2}(x,\ y)}{\varphi(x_{0},\ y_{0})+\varphi(x,\ y)}\\ \varphi^{2}(x_{0},\ y_{0})-\varphi^{2}(x,\ y)&=&f^{2}_{kx}(x_{0},\ y_{0})-f_{kx}^{2}(x,\ y)+f^{2}_{ky}(x_{0},\ y_{0})-f_{ky}^{2}(x,\ y).\end{eqnarray*}
If we assume,\[ M_{k}=\max_{P=(x,\ y,\ z)^{T}\in \overline{U}(P_{k},\ \delta_{1}(P_{k}))\cap\partial\Omega}(|f_{kx}(x,\ y)|,\ |f_{ky}(x,\ y)|),\ 1\leq k\leq N,\]then from\ $\varphi(x,\ y)\geq 1$, we can obtain,\[\mid n_{p}-n_{p_{0}}\mid\leq (4M_{k}^{2}+4M_{k}+1)C_{\alpha}(P_{k})r_{PP_{0}}^{\alpha},\ \forall P=(x,\ y,\ z)^{T}\in U(P_{0},\ \delta_{0})\cap\partial\Omega,\] where\ $C_{\alpha}(P_{k})$\ is defined in the same way as in Theorem 2.1.\\
If we denote,\[ M_{0}=\max_{1\leq k\leq N}M_{k},\ C_{\alpha}=\max_{1\leq k\leq N}C_{\alpha}(P_{k}),\ C_{0}=(4M_{0}^{2}+4M_{0}+1)C_{\alpha},\]
then we know (2.14) holds.\qed\\
\begin{theorem} \label{Theorem2-2} If\ $\Omega$\ is bounded,\ $\partial\Omega\in C^{1,\ \alpha},\ 0<\alpha\leq 1,$\ then we have the following for an ababsolute solid angle that is defined in Lyapunov's potential theory on page 182 of [6],
\[ \max_{M\in R^{3}}\int_{\partial\Omega}\cfrac{|{\bf r}_{MP}\cdot n_{p}|}{r_{MP}^{3}}dS_{P}<+\infty.\]\end{theorem}
{\it Proof of theorem 2.2}. If\ $M=(x,\ y,\ z)^{T}\in\partial\Omega$, then from lemma 2.1, we know there exists a neighborhood\ $U(M,\ \delta(M)),\ \delta(M)>0$,\ and\ $C_{\alpha}>0$, such that\ $|{\bf r}_{MP}\cdot n_{p}|\leq 2C_{\alpha} r_{MP}^{(1+\alpha)}$,\ $\forall P\in U(M,\ \delta(M))\cap \partial\Omega$, where\ $C_{\alpha}$\ is defined in the same way as in lemma 2.3. \\We obtain the following,$$\cfrac{|{\bf r}_{MP}\cdot n_{p}|}{r_{MP}^{3}}\leq \cfrac{2C_{\alpha}}{r_{MP}^{2-\alpha}},\ 0<\alpha\leq 1,\ \forall P\in U(M,\ \delta(M))\cap \partial\Omega.$$
From lemma 2.2, we can see that\[ \int_{\partial\Omega}\cfrac{|{\bf r}_{MP}\cdot n_{p}|}{r_{MP}^{3}}dS_{P},\] is continuous on\ $\partial\Omega$. Hence we get the following,\[ \max_{M\in \partial\Omega}\int_{\partial\Omega}\cfrac{|{\bf r}_{MP}\cdot n_{p}|}{r_{MP}^{3}}dS_{P}<+\infty.\]
Next we assume,
\[\Omega_{1}=\{ M\in R^{3}\setminus\overline{\Omega}:\ dist(M,\ \partial\Omega)\leq \delta_{2}^{\ast}\},\ \Omega_{2}=\{ M\in \Omega:\ dist(M,\ \partial\Omega)\leq \delta_{2}^{\ast}\},\]where$$\delta_{2}^{\ast}=[\delta_{1}^{\ast}\cos(\theta^{\ast})]/2,\ dist(M,\ \partial\Omega)=\min_{P\in\partial\Omega}|M-P|,$$ $\delta_{1}^{\ast},\ \theta^{\ast}$\ are defined in the same way as in Theorem 2.1.
\\If\ $M=(x,\ y,\ z)^{T}\in R^{3}\setminus(\Omega_{1}\cup\Omega_{2})$, then we obtain the following,
$$\int_{\partial\Omega}\cfrac{|{\bf r}_{MP}\cdot n_{p}|}{r_{MP}^{3}}dS_{P}\leq \cfrac{S(\partial\Omega)}{(\delta_{2}^{\ast})^{2}}, $$ where
$$S(\partial\Omega)=\int_{\partial\Omega}dS_{P}.$$Hence we get as follows,\[ \max_{M\in R^{3}\setminus(\Omega_{1}\cup\Omega_{2})}\int_{\partial\Omega}\cfrac{|{\bf r}_{MP}\cdot n_{p}|}{r_{MP}^{3}}dS_{P}<+\infty.\]
\\If\ $M=(x,\ y,\ z)^{T}\in \Omega_{1}$, then there exists\ $P_{0}\in\partial\Omega$, such that\ $|M-P_{0}|=\ dist(M,\ \partial\Omega)$. We discuss a smooth curve on\ $\partial\Omega$\ that passes through\ $P_{0}$. The parameter coordinates of the point\ $P$\ on the curve are\ $(x(\theta),\ y(\theta),\ z(\theta))$.The parameter coordinates of \ $P_{0}$\ are\ $(x(\theta_{0}),\ y(\theta_{0}),\ z(\theta_{0}))$. We assume the tangent vector at\ $P_{0}$\ is as follows,
$$ s_{0}=(x^{\prime}(\theta_{0}),\ y^{\prime}(\theta_{0}),\ z^{\prime}(\theta_{0}))^{T}.$$ We denote\ $f(\theta)$\ as follows,
$$f(\theta)=r_{MP}^{2}=(x-x(\theta))^{2}+(y-y(\theta))^{2}+(z-z(\theta))^{2}.$$ Then\ $f(\theta)$\ attains the minimum at\ $\theta_{0}$. Since\ $f(\theta)$\ is smooth, we see that$$ f^{\prime}(\theta)|_{\theta=\theta_{0}}=-2(M-P_{0})\cdot s_{0}=0.$$ So\ ${\bf r}_{MP_{0}}$\ is perpendicular to\ $s_{0}$. From the arbitrary of tangent vector\ $s_{0}$, we obtain that\ ${\bf r}_{MP_{0}}$\ is parallel to\ $n_{p_{0}}$. \\
From Theorem 2.1, we know that there exists a uniform exterior finite right circular cone\ $K$\ for\ $P_{0}$,\ $P_{0}$\ is the vertex of a cone\ $K(P_{0})$,\ $n_{p_{0}}$\ is the symmetry axis, the polar angle is\ $\theta^{\ast}$,\ and the length of generatrix is\ $\delta_{1}^{\ast}$. From\ $|M-P_{0}|\leq \delta_{2}^{\ast}=[\delta_{1}^{\ast}\cos(\theta^{\ast})]/2$,\ and\ $M\in\Omega_{1},\ {\bf r}_{MP_{0}}$\ is parallel to\ $n_{p_{0}}$, and we can obtain\ $M\in K(P_{0})$, and\ $M$\ is on the symmetry axis\ $n_{p_{0}}$.\\
Again from Theorem 2.1, we know there exists\ $P_{k}\in \partial\Omega,\ 1\leq k\leq N$, such that\ $P_{0}\in U(P_{k},\ \delta_{1}(P_{k})/3)\cap\partial\Omega$. And\ $U(P_{k},\ \delta_{1}(P_{k}))\cap\partial\Omega$\ is a graph of a\ $C^{1,\ \alpha},\ 0<\alpha\leq 1,$\ function of two of the coordinates\ $x,\ y,\ z$. \\From\ $\delta_{1}^{\ast}\leq \delta_{1}(P_{k})/3$, we get\ $U(P_{0},\ \delta_{1}^{\ast})\cap\partial\Omega\subset U(P_{k},\ \delta_{1}(P_{k}))\cap\partial\Omega.$\\
If we denote\ $\theta_{1}=({\bf r}_{P_{0}M},\ {\bf r}_{P_{0}P}),\ \forall P\in U(P_{0},\ \delta_{1}^{\ast})\cap\partial\Omega$, where\ $({\bf r}_{P_{0}M},\ {\bf r}_{P_{0}P})$\ is the angle between\ ${\bf r}_{P_{0}M}$\ and\ ${\bf r}_{P_{0}P}$, then we obtain\ $\theta_{1}> \theta^{\ast}$. We denote\ $\theta_{2}=({\bf r}_{P_{0}M},\ {\bf r}_{PM})$.\\ If\ $\theta_{1}<\pi/2$, then we obtain
$$r_{P_{0}P}=\cfrac{ r_{MP}\sin\theta_{2}}{\sin\theta_{1}}\leq \cfrac{1}{\sin\theta^{\ast}} r_{MP}.$$ If\ $\theta_{1}\geq\pi/2$, then we obtain
$$r_{P_{0}P}\leq  r_{MP}\leq\cfrac{1}{\sin\theta^{\ast}} r_{MP},\ \mbox{where}\ \theta^{\ast}\in(0,\ \cfrac{\pi}{2}).$$ Hence we obtain
\[r_{P_{0}P}\leq\cfrac{1}{\sin\theta^{\ast}} r_{MP},\ \forall P\in U(P_{0},\ \delta_{1}^{\ast})\cap\partial\Omega.\]
From lemma 2.3, we can obtain\ $\forall P\in U(P_{0},\ \delta_{1}^{\ast})\cap\partial\Omega,$
\begin{eqnarray*}|\cfrac{{\bf r}_{MP}\cdot n_{p}}{r_{MP}^{3}}-\cfrac{{\bf r}_{MP}\cdot n_{p_{0}}}{r_{MP}^{3}}|&=&|\cfrac{{\bf r}_{MP}\cdot (n_{p}-n_{p_{0}})}{r_{MP}^{3}}|\\&\leq& C_{0}\cfrac{r_{P_{0}P}^{\alpha}}{r_{MP}^{2}}\leq \cfrac{C_{0}}{(\sin\theta^{\ast})^{\alpha}}\cfrac{1}{r_{MP}^{2-\alpha}}.\end{eqnarray*}
And from lemma 2.1, we can get\ $\forall P\in U(P_{0},\ \delta_{1}^{\ast})\cap\partial\Omega,$\begin{eqnarray*}|\cfrac{{\bf r}_{MP}\cdot n_{p_{0}}}{r_{MP}^{3}}-\cfrac{{\bf r}_{MP_{0}}\cdot n_{p_{0}}}{r_{MP}^{3}}|&=&|\cfrac{{\bf r}_{P_{0}P}\cdot n_{p_{0}}}{r_{MP}^{3}}|\\&\leq& 2C_{\alpha}\cfrac{r_{P_{0}P}^{1+\alpha}}{r_{MP}^{3}}\leq \cfrac{2C_{\alpha}}{(\sin\theta^{\ast})^{1+\alpha}}\cfrac{1}{r_{MP}^{2-\alpha}}.\end{eqnarray*}
Hence we can see the following,\begin{eqnarray*}\int_{U(P_{0},\ \delta_{1}^{\ast})\cap\partial\Omega}\cfrac{|{\bf r}_{MP}\cdot n_{p}|}{r_{MP}^{3}}dS_{P}&\leq&\int_{U(P_{0},\ \delta_{1}^{\ast})\cap\partial\Omega}\cfrac{|{\bf r}_{MP_{0}}\cdot n_{p_{0}}|}{r_{MP}^{3}}dS_{P}+\\&&[\cfrac{C_{0}}{(\sin\theta^{\ast})^{\alpha}}+\cfrac{2C_{\alpha}}{(\sin\theta^{\ast})^{1+\alpha}}]\int_{U(P_{0},\ \delta_{1}^{\ast})\cap\partial\Omega}\cfrac{1}{r_{MP}^{2-\alpha}}dS_{P}.\end{eqnarray*}From lemma 2.2, we know$$\int_{\partial\Omega}\cfrac{1}{r_{MP}^{2-\alpha}}dS_{P}$$ is continuous on\ $M$. So we can assume,\[ C_{1}=\max_{M\in \Omega_{1}\cup\Omega_{2}\cup\partial\Omega}\int_{\partial\Omega}\cfrac{1}{r_{MP}^{2-\alpha}}dS_{P}.\]
Hence we can obtain\[\int_{U(P_{0},\ \delta_{1}^{\ast})\cap\partial\Omega}\cfrac{|{\bf r}_{MP}\cdot n_{p}|}{r_{MP}^{3}}dS_{P}\leq\int_{U(P_{0},\ \delta_{1}^{\ast})\cap\partial\Omega}\cfrac{|{\bf r}_{MP_{0}}\cdot n_{p_{0}}|}{r_{MP}^{3}}dS_{P}+[\cfrac{C_{0}}{(\sin\theta^{\ast})^{\alpha}}+\cfrac{2C_{\alpha}}{(\sin\theta^{\ast})^{1+\alpha}}]C_{1}.\]
From\ $M$\ is on\ $n_{p_{0}}$, we can get that\ $|{\bf r}_{MP_{0}}\cdot n_{p_{0}}|=r_{MP_{0}}$. From the cosine law, we can obtain
\begin{eqnarray*}r_{MP}^{2}&=&r_{MP_{0}}^{2}+r_{P_{0}P}^{2}-2r_{MP_{0}}r_{P_{0}P}\cos\theta_{1}\\
&\geq&r_{MP_{0}}^{2}+r_{P_{0}P}^{2}-2r_{MP_{0}}r_{P_{0}P}\cos\theta^{\ast}\\&\geq&(1-\cos\theta^{\ast})(r_{MP_{0}}^{2}+r_{P_{0}P}^{2}).\end{eqnarray*}
We assume\ $U(P_{j},\ \delta_{1}(P_{j}))\cap\partial\Omega$\ is a graph of a\ $C^{1,\ \alpha},\ 0<\alpha\leq 1,$\ function\ $f_{j}(x_{j1},\ x_{j2})$, where\ $x_{j1},\ x_{j2}$\ are two of the coordinates\ $x,\ y,\ z,\ 1\leq j\leq N$. We assume
$$ C_{2}=\max_{1\leq j\leq N}C_{2}(P_{j}),\ C_{2}(P_{j})=\max_{P=(x,\ y,\ z)^{T}\in \overline{U}(P_{j},\ \delta_{1}(P_{j}))}\sqrt{f^{2}_{j1}(P(x_{j1}),\ P(x_{j2}))+f^{2}_{j2}(P(x_{j1}),\ P(x_{j2}))+1},$$where\ $f_{j1},\ f_{j2}$\ are partial derivatives of\ $f_{j},\ P(x_{jl})$\ is the value of coordinate\ $x_{jl}$\ at point\ $P$,\ $l=1,\ 2,\ 1\leq j\leq N$.\\Now we can obtain the following,
\begin{eqnarray*} \int_{U(P_{0},\ \delta_{1}^{\ast})\cap\partial\Omega}\cfrac{r_{MP_{0}}}{r_{MP}^{3}}dS_{P}&\leq &\cfrac{1}{(1-\cos\theta^{\ast})^{3/2}} \int_{U(P_{0},\ \delta_{1}^{\ast})\cap\partial\Omega}\cfrac{r_{MP_{0}}}{(r_{MP_{0}}^{2}+r_{P_{0}P}^{2})^{3/2}}dS_{P} \\(a=r_{MP_{0}})
&\leq& \cfrac{2\pi C_{2}}{(1-\cos\theta^{\ast})^{3/2}}\int_{0}^{\delta_{1}^{\ast}}\cfrac{ardr}{(r^{2}+a^{2})^{3/2}}.\end{eqnarray*}
And we can work out the following,\ $\forall a>0$,
$$ \int_{0}^{\delta_{1}^{\ast}}\cfrac{ardr}{(r^{2}+a^{2})^{3/2}}= \cfrac{-a}{\sqrt{r^{2}+a^{2}}}|_{0}^{\delta_{1}^{\ast}}
= 1-\cfrac{a}{\sqrt{(\delta_{1}^{\ast})^{2}+a^{2}}}\leq1.$$
From (2.26), we can obtain,\[\int_{U(P_{0},\ \delta_{1}^{\ast})\cap\partial\Omega}\cfrac{|{\bf r}_{MP}\cdot n_{p}|}{r_{MP}^{3}}dS_{P}\leq\cfrac{2\pi C_{2}}{(1-\cos\theta^{\ast})^{3/2}}+[\cfrac{C_{0}}{(\sin\theta^{\ast})^{\alpha}}+\cfrac{2C_{\alpha}}{(\sin\theta^{\ast})^{1+\alpha}}]C_{1}.\]
And from\ $M\in\Omega_{1}$, we can get\ $r_{MP_{0}}\leq \delta_{2}^{\ast}$. Hence we can obtain$$r_{MP}\geq r_{P_{0}P}-r_{MP_{0}}\geq \delta_{1}^{\ast}-\delta_{2}^{\ast},\ \forall P\in\partial\Omega\setminus U(P_{0},\ \delta_{1}^{\ast}).$$ So we get as follows,
 \[\int_{\partial\Omega\setminus U(P_{0},\ \delta_{1}^{\ast})}\cfrac{|{\bf r}_{MP}\cdot n_{p}|}{r_{MP}^{3}}dS_{P}\leq  \cfrac{S(\partial\Omega)}{(\delta_{1}^{\ast}-\delta_{2}^{\ast})^{2}}, \]
where\ $$S(\partial\Omega)=\int_{\partial\Omega}dS_{P}.$$Hence we get the following,\ $\forall M\in\Omega_{1}$,
\[\int_{\partial\Omega}\cfrac{|{\bf r}_{MP}\cdot n_{p}|}{r_{MP}^{3}}dS_{P}\leq  \cfrac{2\pi C_{2}}{(1-\cos\theta^{\ast})^{3/2}}+[\cfrac{C_{0}}{(\sin\theta^{\ast})^{\alpha}}
+\cfrac{2C_{\alpha}}{(\sin\theta^{\ast})^{1+\alpha}}]C_{1}+\cfrac{S(\partial\Omega)}{(\delta_{1}^{\ast}-\delta_{2}^{\ast})^{2}}. \]In the same way, we can obtain following,\ $\forall M\in\Omega_{2}$,
\[\int_{\partial\Omega}\cfrac{|{\bf r}_{MP}\cdot n_{p}|}{r_{MP}^{3}}dS_{P}\leq  \cfrac{2\pi C_{2}}{(1-\cos\theta^{\ast})^{3/2}}+[\cfrac{C_{0}}{(\sin\theta^{\ast})^{\alpha}}
+\cfrac{2C_{\alpha}}{(\sin\theta^{\ast})^{1+\alpha}}]C_{1}+\cfrac{S(\partial\Omega)}{(\delta_{1}^{\ast}-\delta_{2}^{\ast})^{2}}, \]
which proves the statement.\qed \\
\begin{corollary} If\ $\Omega$\ is bounded,\ $\partial\Omega\in C^{1,\ \alpha},\ 0<\alpha\leq 1,$\ then double layer potential
\[u(M)=\int_{\partial\Omega}v(P)\cfrac{\partial\cfrac{1}{r_{PM}}}{\partial n_{p}}dS_{P},\] is continuous on\ $R^{3}\setminus \partial\Omega$, where\ $v(P)\in C(\partial\Omega)$, moreover\ $\forall P_{0}\in \partial\Omega$, we have\begin{eqnarray}&& \lim_{M\rightarrow P_{0}+}u(M)=u(P_{0})-2\pi v(P_{0}),\\&&  \lim_{M\rightarrow P_{0}-}u(M)=u(P_{0})+2\pi v(P_{0}),\end{eqnarray}where\ $M\rightarrow P_{0}+$\ means\ $M$\ is near to\ $P_{0}$\ from the interior of\ $\Omega$\ and\ $M\rightarrow P_{0}-$\ means\ $M$\ is near to\ $P_{0}$\ from the exterior of\ $\Omega$. \end{corollary}
{\it Proof of corollary 2.1}. From lemma 2.2 and the previous Theorem, we may get that\ $\forall P_{0}\in \partial\Omega$,
$$u_{0}(M)=\int_{\partial\Omega}(v(P)-v(P_{0}))\cfrac{\partial\cfrac{1}{r_{PM}}}{\partial n_{p}}dS_{P}$$ is continuous at\ $P_{0}$. From the potential theory, we know$$\int_{\partial\Omega}\cfrac{\partial\cfrac{1}{r_{PM}}}{\partial n_{p}}dS_{P}
=\begin{cases}-4\pi,\ M\in \Omega,\\-2\pi,\ M\in \partial\Omega,\\0,\ M\in R^{3}\setminus\overline{\Omega}.\end{cases}
$$Hence, the statement holds.\qed \\
\begin{theorem} \label{Theorem2-3} If\ $\Omega$\ is bounded,\ $\partial\Omega\in C^{1,\ \alpha},\ 0<\alpha\leq 1,$\ then simple layer potential
\[u(M)=\int_{\partial\Omega}\cfrac{v(P)}{r_{PM}}dS_{P},\] where\ $v(P)\in C(\partial\Omega)$, satisfies the following,\ $\forall P_{0}\in \partial\Omega$,
\begin{eqnarray}&& \cfrac{\partial u(P_{0})}{\partial n_{p_{0}}^{+}}=\int_{\partial\Omega}v(P)\cfrac{\partial\cfrac{1}{r_{PP_{0}}}}{\partial n_{p_{0}}}dS_{P}-2\pi v(P_{0}),\\&&  \cfrac{\partial u(P_{0})}{\partial n_{p_{0}}^{-}}=\int_{\partial\Omega}v(P)\cfrac{\partial\cfrac{1}{r_{PP_{0}}}}{\partial n_{p_{0}}}dS_{P}+2\pi v(P_{0}),\end{eqnarray} where\[ \cfrac{\partial u(P_{0})}{\partial n_{p_{0}}^{+}}=\lim_{M\rightarrow n_{p_{0}}^{+}}\cfrac{u(M)-u(P_{0})}{r_{MP_{0}}},\ \cfrac{\partial u(P_{0})}{\partial n_{p_{0}}^{-}}=\lim_{M\rightarrow n_{p_{0}}^{-}}\cfrac{u(P_{0})-u(M)}{r_{P_{0}M}},\]
 where\ $M\rightarrow n_{p_{0}}^{+}$\ means\ $M$\ is near to\ $P_{0}$\ along\ $n_{p_{0}}$\ from the exterior of\ $\Omega$\ and\ $M\rightarrow n_{p_{0}}^{-}$\ means\ $M$\ is near to\ $P_{0}$\ along\ $n_{p_{0}}$\ from the interior of\ $\Omega$. \end{theorem}
{\it Proof of theorem 2.3}. We refer to the proof on pages 190 to page 193 of [6]. If\ $M\in n_{p_{0}}\setminus \partial\Omega$, then a directional deriviative\ $\partial u(M)/\partial n_{p_{0}}$\ exists, and we can work it out through the integral as follows.
\[\cfrac{\partial u(M)}{\partial n_{p_{0}}}=\int_{\partial\Omega}v(P)\cfrac{\partial\cfrac{1}{r_{PM}}}{\partial n_{p_{0}}}dS_{P}
=-\int_{\partial\Omega}v(P)\cfrac{\cos({\bf r}_{PM},\ n_{p_{0}})}{r_{PM}^{2}}dS_{P},\] where\ $({\bf r}_{PM},\ n_{p_{0}})$\ is the angle between\ ${\bf r}_{PM}$\ and\ $n_{p_{0}}$. Together with double potential\[u_{1}(M)=\int_{\partial\Omega}v(P)\cfrac{\partial\cfrac{1}{r_{PM}}}{\partial n_{p}}dS_{P}=\int_{\partial\Omega}v(P)\cfrac{\cos({\bf r}_{PM},\ n_{p})}{r_{PM}^{2}}dS_{P},\] we have
\[\cfrac{\partial u(M)}{\partial n_{p_{0}}}+u_{1}(M)=\int_{\partial\Omega}v(P)\cfrac{\cos({\bf r}_{PM},\ n_{p})-\cos({\bf r}_{PM},\ n_{p_{0}})}{r_{PM}^{2}}dS_{P}.\]
We want to prove the right side of (2.40) is continuous when\ $M$\ is near to\ $P_{0}$\ along\ $n_{P_{0}}$. We only need to prove\ $\forall \epsilon>0,\ \exists\delta>0$, such that\ $\forall M\in U(P_{0},\ \delta)\cap n_{p_{0}}$,\[\mid\int_{(\partial\Omega)_{\delta}}v(P)\cfrac{\cos({\bf r}_{PM},\ n_{p})-\cos({\bf r}_{PM},\ n_{p_{0}})}{r_{PM}^{2}}dS_{P}\mid\leq\epsilon, \]
where\ $(\partial\Omega)_{\delta}=U(P_{0},\ \delta)\cap\partial\Omega$.\\ If we assume\[\max_{P\in \partial\Omega}\mid v(P)\mid=C_{3},\]
then we have\begin{eqnarray}\mid v(P)\cfrac{\cos({\bf r}_{PM},\ n_{p})-\cos({\bf r}_{PM},\ n_{p_{0}})}{r_{PM}^{2}}\mid &\leq& C_{3}\cfrac{\mid \cos({\bf r}_{PM},\ n_{p})-\cos({\bf r}_{PM},\ n_{p_{0}})\mid}{r_{PM}^{2}}\\&\leq& 2C_{3}\cfrac{\mid \sin\cfrac{({\bf r}_{PM},\ n_{p})-({\bf r}_{PM},\ n_{p_{0}})}{2}\mid}{r_{PM}^{2}}\\
&\leq&2C_{3}\cfrac{\mid \sin\cfrac{(n_{p},\ n_{p_{0}})}{2}\mid}{r_{PM}^{2}}=C_{3}\cfrac{\mid n_{p}-n_{p_{0}}\mid}{r_{PM}^{2}},\end{eqnarray}
where $$\mid \cfrac{({\bf r}_{PM},\ n_{p})-({\bf r}_{PM},\ n_{p_{0}})}{2}\mid\\
\leq\mid \cfrac{(n_{p},\ n_{p_{0}})}{2}\mid, $$ is obtained from the sum of two angles of the trihedral being no less than the third one, and\ $(n_{p},\ n_{p_{0}})/2$\ being in\ $[0,\ \pi/2]$.\\
From lemma 2.3 and (2.24), we have\ $\forall P=(x,\ y,\ z)^{T}\in U(P_{0},\ \delta_{1}^{\ast})\cap\partial\Omega,\ \forall M\in n_{p_{0}},\ |M-P_{0}|\leq \delta_{2}^{\ast}$, \[\mid v(P)\cfrac{\cos({\bf r}_{PM},\ n_{p})-\cos({\bf r}_{PM},\ n_{p_{0}})}{r_{PM}^{2}}\mid \leq \cfrac{C_{0}C_{3}}{(\sin\theta^{\ast})^{\alpha}}\cfrac{1}{r_{PM}^{2-\alpha}}.\]
From lemma 2.2, we know (2.41) is true and the right side of (2.40) is continuous when\ $M$\ is near to\ $P_{0}$\ along\ $n_{P_{0}}$. \\
Since the continuity, we have the following,\begin{eqnarray}&&\lim_{M^{\prime}\rightarrow n_{p_{0}}^{+}}(\cfrac{\partial u(M^{\prime})}{\partial n_{p_{0}}}+u_{1}(M^{\prime}))=\lim_{M^{\prime\prime}\rightarrow n_{p_{0}}^{-}}(\cfrac{\partial u(M^{\prime\prime})}{\partial n_{p_{0}}}+u_{1}(M^{\prime\prime}))\\&&=\int_{\partial\Omega}v(P)\cfrac{\cos({\bf r}_{PP_{0}},\ n_{p})-\cos({\bf r}_{PP_{0}},\ n_{p_{0}})}{r_{PP_{0}}^{2}}dS_{P}.\end{eqnarray}
From corollary 2.1, we have\begin{eqnarray}&&\lim_{M^{\prime}\rightarrow n_{p_{0}}^{+}}u_{1}(M^{\prime})=\lim_{M^{\prime}\rightarrow P_{0}-}u_{1}(M^{\prime})
=\int_{\partial\Omega}v(P)\cfrac{\partial\cfrac{1}{r_{PP_{0}}}}{\partial n_{p}}dS_{P}+2\pi v(P_{0})\\&&=\int_{\partial\Omega}v(P)\cfrac{\cos({\bf r}_{PP_{0}},\ n_{p})}{r_{PP_{0}}^{2}}dS_{P}+2\pi v(P_{0}),\end{eqnarray}\begin{eqnarray}&&\lim_{M^{\prime\prime}\rightarrow n_{p_{0}}^{+}}u_{1}(M^{\prime\prime})=\lim_{M^{\prime\prime}\rightarrow P_{0}-}u_{1}(M^{\prime\prime})
=\int_{\partial\Omega}v(P)\cfrac{\partial\cfrac{1}{r_{PP_{0}}}}{\partial n_{p}}dS_{P}-2\pi v(P_{0})\\&&=\int_{\partial\Omega}v(P)\cfrac{\cos({\bf r}_{PP_{0}},\ n_{p})}{r_{PP_{0}}^{2}}dS_{P}-2\pi v(P_{0}).\end{eqnarray}
And from (2.47) and (2.48) we can push out the limits\[\cfrac{\partial u(P_{0})}{\partial n_{p_{0}}^{+}}=\lim_{M^{\prime}\rightarrow n_{p_{0}}^{+}}\cfrac{\partial u(M^{\prime})}{\partial n_{p_{0}}},\ \cfrac{\partial u(P_{0})}{\partial n_{p_{0}}^{-}}=\lim_{M^{\prime\prime}\rightarrow n_{p_{0}}^{-}}\cfrac{\partial u(M^{\prime\prime})}{\partial n_{p_{0}}},\]and integrals\[\int_{\partial\Omega}v(P)\cfrac{-\cos({\bf r}_{PP_{0}},\ n_{p_{0}})}{r_{PP_{0}}^{2}}dS_{P}\]all exist, moreover\begin{eqnarray}&&\cfrac{\partial u(P_{0})}{\partial n_{p_{0}}^{+}}=\lim_{M^{\prime}\rightarrow n_{p_{0}}^{+}}\cfrac{\partial u(M^{\prime})}{\partial n_{p_{0}}}=\int_{\partial\Omega}v(P)\cfrac{-\cos({\bf r}_{PP_{0}},\ n_{p_{0}})}{r_{PP_{0}}^{2}}dS_{P}-2\pi v(P_{0})\\&&=\int_{\partial\Omega}v(P)\cfrac{\partial\cfrac{1}{r_{PP_{0}}}}{\partial n_{p_{0}}}dS_{P}-2\pi v(P_{0}),\\&&\cfrac{\partial u(P_{0})}{\partial n_{p_{0}}^{-}}=\lim_{M^{\prime\prime}\rightarrow n_{p_{0}}^{-}}\cfrac{\partial u(M^{\prime\prime})}{\partial n_{p_{0}}}=\int_{\partial\Omega}v(P)\cfrac{-\cos({\bf r}_{PP_{0}},\ n_{p_{0}})}{r_{PP_{0}}^{2}}dS_{P}+2\pi v(P_{0})\\&&=\int_{\partial\Omega}v(P)\cfrac{\partial\cfrac{1}{r_{PP_{0}}}}{\partial n_{p_{0}}}dS_{P}+2\pi v(P_{0}).\end{eqnarray}
That's the end of proof.\qed \\These last two theorems are the main results in Lyapunov's potential theory on pages 173 to 193 of [6]. We see that\ $\partial\Omega\in C^{1,\ \alpha},\ 0<\alpha\leq 1,$ can play the role of Lyapunov's surface.\\
Finally, we see whether\ $\partial\Omega\in C^{1,\ \alpha},\ 0<\alpha\leq 1$, is the Hopf's surface for the strong maximum principle. From the example on page 35 of [2], we know it is not the Hopf's surface if domain\ $\Omega$\ only satisfies an interior cone condition. However, we can obtain a stronger condition as follows.\begin{theorem} \label{Theorem2-4} If\ $\Omega$\ is bounded,\ $\partial\Omega\in C^{1,\ \alpha},\ 0<\alpha\leq 1,$\ then domain\ $\Omega$\ satisfies a uniform interior oblate spheroid condition, that is,\ $\exists \delta>0,\ \forall P_{0}\in \partial\Omega$, there exists a finite right oblate spheroid as follows,\[K_{\delta}(P_{0})=\{P: {\bf r}_{P_{0}P}\cdot(-n_{p_{0}})\geq 4C_{\alpha}r_{P_{0}P}^{1+\alpha},\ r_{P_{0}P}\leq \delta\},\]
where\ $C_{\alpha}$\ is defined in the same way as in lemma 2.3,\ $K_{\delta}(P_{0})\cap(R^{3}\setminus\Omega)=P_{0}$. \end{theorem}
{\it Proof of theorem 2.4}. If we select\ $\delta\leq\delta_{1}^{\ast}$, then we will get\ $K_{\delta}(P_{0})\cap(R^{3}\setminus\Omega)=P_{0}$.\\ If there exists\ $P\in K_{\delta}(P_{0})\cap(R^{3}\setminus\Omega)$, and\ $P\neq P_{0}$, then from lemma 2.1, we obtain\[ {\bf r}_{P_{0}P}\cdot(-n_{p_{0}})\leq 2C_{\alpha}r_{P_{0}P}^{1+\alpha}.\] This contradicts (2.59).\qed \\ Now we can see that domain\ $\Omega$\ satisfies not only a uniform exterior and interior cone condition, but also a uniform exterior and interior oblate spheroid condition, if\ $\Omega$\ is bounded,\ $\partial\Omega\in C^{1,\ \alpha},\ 0<\alpha\leq 1$. Therefore, we arrive at the following.\begin{theorem} \label{Theorem2-5} If\ $\Omega$\ is bounded,\ $\partial\Omega\in C^{1,\ \alpha},\ 0<\alpha\leq 1$, supposing\ $u\in C^{2}(\Omega)$\ and\ $\triangle u=0$\ in\ $\Omega$, letting\ $P_{0}\in \partial\Omega$\ be such that\ $u$\ is continuous at\ $P_{0}$,\ $u(P_{0})>u(P)$\ for all\ $P\in \Omega$, then the exterior normal derivative of\ $u$\ at\ $P_{0}$, if it exists, satisfies the strict inequality\[ \cfrac{\partial u}{\partial n}(P_{0})=\lim_{P\rightarrow n_{p_{0}}^{-}}\cfrac{u(P_{0})-u(P)}{r_{P_{0}P}}>0.\]\end{theorem}
{\it Proof of theorem 2.5}. From the previous theorem, we know\ $\forall \delta\in (0,\ \delta_{1}^{\ast}]$, there exists an oblate spheroid\ $K_{\delta}(P_{0})\subset \Omega$. We introduce an auxiliary function\ $v$\ by defining\[ v(P)=e^{-(4C_{\alpha}r_{P_{0}P}^{1+\alpha})^{\gamma}}-e^{-({\bf r}_{P_{0}P}\cdot(-n_{p_{0}}))^{\gamma}}+e^{-(4C_{\alpha}r_{P_{0}P}^{1+\alpha})}-e^{-({\bf r}_{P_{0}P}\cdot(-n_{p_{0}}))},\]
where\ $\gamma\in(1,\ 1+\alpha)$. Direct calculation gives
\begin{eqnarray*}\triangle v&=&e^{-(4C_{\alpha}r_{P_{0}P}^{1+\alpha})^{\gamma}}[(4C_{\alpha})^{2\gamma}(\gamma+\gamma\alpha)^{2}r_{P_{0}P}^{2\gamma+2\gamma\alpha-2}
-(4C_{\alpha})^{\gamma}(\gamma+\gamma\alpha)(\gamma+\gamma\alpha-1)r_{P_{0}P}^{\gamma+\gamma\alpha-2}]\\ &&-e^{-({\bf r}_{P_{0}P}\cdot(-n_{p_{0}}))^{\gamma}}[\gamma^{2}({\bf r}_{P_{0}P}\cdot(-n_{p_{0}}))^{2\gamma-2}-\gamma(\gamma-1)({\bf r}_{P_{0}P}\cdot(-n_{p_{0}}))^{\gamma-2}]\\&&+e^{-(4C_{\alpha}r_{P_{0}P}^{1+\alpha})}[(4C_{\alpha})^{2}(1+\alpha)^{2}r_{P_{0}P}^{2\alpha}
-(4C_{\alpha})(1+\alpha)\alpha r_{P_{0}P}^{\alpha-1}]-e^{-({\bf r}_{P_{0}P}\cdot(-n_{p_{0}}))}.\end{eqnarray*}
From$$ \lim_{r_{P_{0}P}\rightarrow 0}({\bf r}_{P_{0}P}\cdot(-n_{p_{0}}))^{2-\gamma}r_{P_{0}P}^{\alpha-1}=\lim_{r_{P_{0}P}\rightarrow 0}({\bf r}_{P_{0}P}\cdot(-n_{p_{0}}))^{2-\gamma}r_{P_{0}P}^{\gamma+\gamma\alpha-2}=0,$$we obtain$$\lim_{r_{P_{0}P}\rightarrow 0}\triangle v=+\infty.$$
Hence, there exists\ $ \delta_{1}\in (0,\ \delta_{1}^{\ast}]$, such that\ $\triangle v>0$, throughout oblate spheroid\ $K_{\delta_{1}}(P_{0})$.\\
If\ ${\bf r}_{P_{0}P}\cdot(-n_{p_{0}})=4C_{\alpha}r_{P_{0}P}^{1+\alpha}$, then\ $v(P)=0$. Therefore, there exists a constant\ $\epsilon>0$, for which\ $u-u(P_{0})+\epsilon v\leq 0$, on\ $\partial K_{\delta_{1}}(P_{0})$. The weak maximum principle now implies that\ $u-u(P_{0})+\epsilon v\leq 0$, throughout oblate spheroid\ $K_{\delta_{1}}(P_{0})$. Taking the exterior normal derivative of at\ $P_{0}$, we obtain, as required,
\[\cfrac{\partial u}{\partial n}(P_{0})\geq -\epsilon\cfrac{\partial v}{\partial n}(P_{0})=\epsilon>0.\] That's the end of the proof.
\qed \\
\section{Equivalence} \setcounter{equation}{0}
We transform Eq.(1.1) into the following,\[v_{j}-\sum_{k=1}^{9m}c_{jk}v_{m+k}-s_{j}=0,\ 1\leq j\leq m,\]where\ $c_{jk},\ 1\leq j\leq m,\ 1\leq k\leq 9m,$\ are all real constants to be determined,\[s_{j}=f_{j}(v_{m+1},\ v_{m+2},\cdots,\ v_{10m},\ x,\ y,\ z)-\sum_{k=1}^{9m}c_{jk}v_{m+k},\ 1\leq j\leq m.\]
Let's introduce\ $Z=(z_{j}(x,\ y,\ z))_{10m\times 1}=(V_{2}^{T},\ V_{3}^{T},\ \cdots,\ V_{10}^{T},\ S^{T})^{T}$, where\[V_{n}=(v_{m(n-1)+1},\ v_{m(n-1)+2},\cdots,\ v_{mn})^{T},\ 1\leq n\leq10,\ S=(s_{1},\ s_{2},\cdots,\ s_{m})^{T}.\] Eq.(1.1) is equivalent to\[V_{1}=\beta^{T}Z,\] where\ $ \beta^{T}=(C_{1},\ C_{2},\ \cdots,\ C_{9},\ E),\ C_{l}=(c_{j,(m(l-1)+k)})_{m\times m},\ 1\leq l\leq 9$,\ $E$\ is\ $m$\ order identity matrix. \\We should discuss\ $Z$\ as follows,\begin{eqnarray}&&V_{1}=\beta^{T}Z,\ V_{j}=E_{j-1}^{T}Z,\ 2\leq j\leq 10,\ S=E_{10}^{T}Z,\end{eqnarray} where\ $E_{n}=(e_{m(n-1)+1},\ e_{m(n-1)+2},\cdots,\ e_{mn}),\ 1\leq n\leq 10,\ e_{k}$\ is the\ $k$th unit coordinate vector,\ $1\leq k\leq 10m$.\\ Because\ $v_{1},\ v_{2},\cdots,\ v_{10m}$\ is a permutation of the components of\ $u,\ \partial u,\ \partial^{2}u$,\ $Z$\ should satisfies the following,
\begin{eqnarray*}&&u=A_{0}Z,\ u_{x}=A_{01}Z,\ u_{y}=A_{02}Z,\ u_{z}=A_{03}Z,\ u_{xx}=A_{11}Z,\\&& u_{xy}=A_{12}Z,\ u_{xz}=A_{13}Z,\ u_{yy}=A_{22}Z,\ u_{yz}=A_{23}Z,\ u_{zz}=A_{33}Z,\end{eqnarray*} where the rows of\ $A_{0},\ A_{01},\ A_{02},\ A_{03},\ A_{11},\ A_{12},\ A_{13},\ A_{22},\ A_{23},\ A_{33}$\ is a permutation of the rows of\ $\beta^{T},\ E_{j}^{T},\ 1\leq j\leq 9$, which is correspondent to\ $v_{j},\ 1\leq j\leq 10m$.\\Hence,\ $Z$\ should satisfies the following system,\begin{eqnarray}&& \cfrac{\partial A_{0}Z}{\partial x}=A_{01}Z,\\&& \cfrac{\partial A_{0}Z}{\partial y}=A_{02}Z,\\&& \cfrac{\partial A_{0}Z}{\partial z}=A_{03}Z,\\&&\cfrac{\partial^{2} A_{0}Z}{\partial x^{2}}=A_{11}Z,\\&& \cfrac{\partial^{2} A_{0}Z}{\partial x\partial y}=A_{12}Z,\\&& \cfrac{\partial^{2} A_{0}Z}{\partial x\partial z}=A_{13}Z,\\&& \cfrac{\partial^{2} A_{0}Z}{\partial y^{2}}=A_{22}Z,\\&& \cfrac{\partial^{2} A_{0}Z}{\partial y\partial z}=A_{23}Z,\\&& \cfrac{\partial^{2} A_{0}Z}{\partial z^{2}}=A_{33}Z,\end{eqnarray}\begin{eqnarray} f(E_{1}^{T}Z,\ E_{2}^{T}Z,\ \cdots,\ E_{9}^{T}Z,\ x,\ y,\ z)-\sum_{l=1}^{9}C_{l}E_{l}^{T}Z=E_{10}^{T}Z,\ &&\end{eqnarray}
 where\ $Z=(z_{j}(x,\ y,\ z))_{10m\times 1}\in C(\overline{\Omega}),\ A_{0}Z\in C^{2}(\overline{\Omega}),\ f=(f_{1},\ f_{2},\cdots,\ f_{m})^{T}$.\\
Moreover we have\begin{theorem} \label{Theorem3-1} Eqs.(1.1) is equivalent to the system from Eq.(3.6) to Eq.(3.15) respect to\ $Z$. \end{theorem}
 {\it Proof of theorem 3.1}. If\ $u$\ satisfies Eqs.(1.1), then letting\ $Z=T_{1}(u)=(V_{2},\ V_{3},\ \cdots,\ V_{10},\ S)^{T}$, we obtain\begin{eqnarray*}&&V_{1}=\beta^{T}Z,\ V_{j}=E_{j-1}^{T}Z,\ 2\leq j\leq 10,\ S=E_{10}^{T}Z.\ \end{eqnarray*}Therefore we get\begin{eqnarray*}&&u=A_{0}Z,\ u_{x}=A_{01}Z,\ u_{y}=A_{02}Z,\ u_{z}=A_{03}Z,\ u_{xx}=A_{11}Z,\\&& u_{xy}=A_{12}Z,\ u_{xz}=A_{13}Z,\ u_{yy}=A_{22}Z,\ u_{yz}=A_{23}Z,\ u_{zz}=A_{33}Z,\end{eqnarray*}where the rows of\ $A_{0},\ A_{01},\ A_{02},\ A_{03},\ A_{11},\ A_{12},\ A_{13},\ A_{22},\ A_{23},\ A_{33}$\ is a permutation of the rows of\ $\beta^{T},\ E_{j}^{T},\ 1\leq j\leq 9$, which is correspondent to\ $v_{j},\ 1\leq j\leq 10m$. Hence\ $Z$\ satisfies Eq.(3.6) to Eq.(3.15).\\  If\ $Z$\ satisfies Eq.(3.6) to Eq.(3.15), then letting\ $u=T_{2}(Z)=A_{0}Z$, we obtain
\begin{eqnarray*}&&u=A_{0}Z,\ u_{x}=A_{01}Z,\ u_{y}=A_{02}Z,\ u_{z}=A_{03}Z,\ u_{xx}=A_{11}Z,\\&& u_{xy}=A_{12}Z,\ u_{xz}=A_{13}Z,\ u_{yy}=A_{22}Z,\ u_{yz}=A_{23}Z,\ u_{zz}=A_{33}Z.\end{eqnarray*}Because the rows of\ $A_{0},\ A_{01},\ A_{02},\ A_{03},\ A_{11},\ A_{12},\ A_{13},\ A_{22},\ A_{23},\ A_{33}$\ is a permutation of the rows of\ $\beta^{T},\ E_{j}^{T},\ 1\leq j\leq 9$, which is correspondent to\ $v_{j},\ 1\leq j\leq 10m$, we get\begin{eqnarray*}&&V_{1}=\beta^{T}Z,\ V_{j}=E_{j-1}^{T}Z,\ 2\leq j\leq 10,\ S=E_{10}^{T}Z.\ \end{eqnarray*}It follows that\ $u$\ satisfies Eqs.(1.1).\\Obviously\ $T_{1},\ T_{2}$\ are continuous. Moreover\ $T_{1}(T_{2}(Z))=Z,\ T_{2}(T_{1}(u))=u.$ From definition 1.1, we know the statement stands.\qed\\
 We notice that Eqs(3.6) to (3.14) are good, because they are the first order linear partial differential equations with constant coefficients. Eq.(3.15) could be considered as complicated, but if we assume
$$ Z=\left(
       \begin{array}{c}
         Z_{1} \\
         Z_{2} \\
       \end{array}
     \right),\ \mbox{where}\ Z_{1}\ \mbox{is the first\ $9m$\ componenets of}\ Z,$$ then we will obtain\ $Z_{2}=\psi(Z_{1})$\ from Eq.(3.15).  Next we want to get\ $Z_{1}=T_{0}(Z_{2})$\ from Eqs.(3.6) to (3.14). This is not difficult. From our experience, we guess that\ $T_{0}$\ should be related to the integral equations. We will obtain\ $T_{0}$\ by the Fourier transform on\ $\overline{\Omega}\times [0,\ T]$. At last we will transform Eqs.(3.6) to (3.15) into the equivalent generalized integral equations\ $Z_{1}=T_{0}(\psi(Z_{1}))$.\\ We apply the Fourier transform on\ $\overline{\Omega}\times [0,\ T]$\ on both sides from Eqs.(3.6) to (3.14) as follows, \begin{eqnarray*}\int_{\Omega}e^{-ix\xi_{1}-iy\xi_{2}-iz\xi_{3}}\cfrac{\partial A_{0}Z}{\partial t}dxdydz=\int_{\Omega}e^{-ix\xi_{1}-iy\xi_{2}-iz\xi_{3}}A_{00}Zdxdydz,&&\\ \int_{\Omega}e^{-ix\xi_{1}-iy\xi_{2}-iz\xi_{3}} \cfrac{\partial A_{0}Z}{\partial x}dxdydz=\int_{\Omega}e^{-ix\xi_{1}-iy\xi_{2}-iz\xi_{3}}A_{01}Zdxdydz,&&\\ \int_{\Omega}e^{-ix\xi_{1}-iy\xi_{2}-iz\xi_{3}}\cfrac{\partial A_{0}Z}{\partial y}dxdydz=\int_{\Omega}e^{-ix\xi_{1}-iy\xi_{2}-iz\xi_{3}}A_{02}Zdxdydz,&&\\ \int_{\Omega}e^{-ix\xi_{1}-iy\xi_{2}-iz\xi_{3}} \cfrac{\partial A_{0}Z}{\partial z}dxdydz=\int_{\Omega}e^{-ix\xi_{1}-iy\xi_{2}-iz\xi_{3}}A_{03}Zdxdydz,&&\\ \int_{\Omega}e^{-ix\xi_{1}-iy\xi_{2}-iz\xi_{3}}\cfrac{\partial^{2} A_{0}Z}{\partial x^{2}}dxdydz=\int_{\Omega}e^{-ix\xi_{1}-iy\xi_{2}-iz\xi_{3}}A_{11}Zdxdydz,&&\\ \int_{\Omega}e^{-ix\xi_{1}-iy\xi_{2}-iz\xi_{3}} \cfrac{\partial^{2} A_{0}Z}{\partial x\partial y}dxdydz=\int_{\Omega}e^{-ix\xi_{1}-iy\xi_{2}-iz\xi_{3}}A_{12}Zdxdydz,&&\\ \int_{\Omega}e^{-ix\xi_{1}-iy\xi_{2}-iz\xi_{3}}\cfrac{\partial^{2} A_{0}Z}{\partial x\partial z}dxdydz=\int_{\Omega}e^{-ix\xi_{1}-iy\xi_{2}-iz\xi_{3}}A_{13}Zdxdydz,&&\\ \int_{\Omega}e^{-ix\xi_{1}-iy\xi_{2}-iz\xi_{3}} \cfrac{\partial^{2} A_{0}Z}{\partial y^{2}}dxdydz=\int_{\Omega}e^{-ix\xi_{1}-iy\xi_{2}-iz\xi_{3}}A_{22}Zdxdydz.&&\\ \int_{\Omega}e^{-ix\xi_{1}-iy\xi_{2}-iz\xi_{3}}\cfrac{\partial^{2} A_{0}Z}{\partial y\partial z}dxdydz=\int_{\Omega}e^{-ix\xi_{1}-iy\xi_{2}-iz\xi_{3}}A_{23}Zdxdydz,&&\\ \int_{\Omega}e^{-ix\xi_{1}-iy\xi_{2}-iz\xi_{3}} \cfrac{\partial^{2} A_{0}Z}{\partial z^{2}}dxdydz=\int_{\Omega}e^{-ix\xi_{1}-iy\xi_{2}-iz\xi_{3}}A_{33}Zdxdydz.&&\end{eqnarray*}
The reason why we apply the Fourier transform on\ $\overline{\Omega}$\ instead of\ $R^{4}$\ is that Eqs(1.1) only holds on\ $\overline{\Omega}$. It is possible that they will not stand outside of\ $\overline{\Omega}$. We can't apply the Fourier transform on\ $R^{4}$\ to both sides of Eq(3.6) to (3.14). \\In order to denote this easily, we define the Fourier transform on\ $\overline{\Omega}$\ as
follows.\begin{definition}\label{definition} $\forall\ f(x,\ y,\ z)\in L^{2}(\overline{\Omega}),\ \forall (\xi_{1},\ \xi_{2},\ \xi_{3})^{T}\in R^{3}$,
\begin{eqnarray*}FI(f(x,\ y,\ z))&=&\int_{\overline{\Omega}}f(x,\ y,\ z)e^{-ix\xi_{1}-iy\xi_{2}-iz\xi_{3}}dxdydz\\&=&F(f(x,\ y,\ z)I_{\overline{\Omega}}(x,\ y,\ z)),\end{eqnarray*}
where\ $F$\ means the Fourier transform and\ $I_{\overline{\Omega}}(x,\ y,\ z)$\
is the characteristic function. In the following, we write\ $I_{\overline{\Omega}}(x,\ y,\ z)$\ into\ $I_{\overline{\Omega}}$.\end{definition}
By using divergence theorem, we obtain\begin{eqnarray*}
FI(\cfrac{\partial A_{0}Z}{\partial x})&=&\int_{\Omega}e^{-ix\xi_{1}-iy\xi_{2}-iz\xi_{3}}(\cfrac{\partial A_{0}Z}{\partial x})dxdydz
\\&=&\int_{0}^{T}\int_{\partial\Omega}(A_{0}Z)n_{1}e^{-ix\xi_{1}-iy\xi_{2}-iz\xi_{3}}dxdydz
+\\&&i\xi_{1}\int_{\Omega}e^{-ix\xi_{1}-iy\xi_{2}-iz\xi_{3}}(A_{0}Z)dxdydz
\\&=&g_{1}+ i\xi_{1} FI(A_{0}Z),\end{eqnarray*} \begin{eqnarray*}
FI(\cfrac{\partial A_{0}Z}{\partial y})&=&\int_{\Omega}e^{-ix\xi_{1}-iy\xi_{2}-iz\xi_{3}}(\cfrac{\partial A_{0}Z}{\partial y})dxdydz
\\&=&\int_{0}^{T}\int_{\partial\Omega}(A_{0}Z)n_{2}e^{-ix\xi_{1}-iy\xi_{2}-iz\xi_{3}}dxdydz
+\\&&i\xi_{2}\int_{\Omega}e^{-ix\xi_{1}-iy\xi_{2}-iz\xi_{3}}(A_{0}Z)dxdydz
\\&=&g_{2}+ i\xi_{2} FI(A_{0}Z),\end{eqnarray*} \begin{eqnarray*}
FI(\cfrac{\partial A_{0}Z}{\partial z})&=&\int_{\Omega}e^{-ix\xi_{1}-iy\xi_{2}-iz\xi_{3}}(\cfrac{\partial A_{0}Z}{\partial z})dxdydz
\\&=&\int_{0}^{T}\int_{\partial\Omega}(A_{0}Z)n_{3}e^{-ix\xi_{1}-iy\xi_{2}-iz\xi_{3}}dxdydz
+\\&&i\xi_{3}\int_{\Omega}e^{-ix\xi_{1}-iy\xi_{2}-iz\xi_{3}}(A_{0}Z)dxdydz
\\&=&g_{3}+ i\xi_{3} FI(A_{0}Z),\end{eqnarray*}\begin{eqnarray*}
FI(\cfrac{\partial^{2} A_{0}Z}{\partial x^{2}})&=&\int_{\overline{\Omega}}e^{-ix\xi_{1}-iy\xi_{2}-iz\xi_{3}}(\cfrac{\partial^{2} A_{0}Z}{\partial x^{2}})]dxdydz
\\&=&\int_{\partial\Omega}(\cfrac{\partial A_{0}Z}{\partial x})n_{1}e^{-ix\xi_{1}-iy\xi_{2}-iz\xi_{3}}dS
+\\&&i\xi_{1}\int_{\overline{\Omega}}e^{-ix\xi_{1}-iy\xi_{2}-iz\xi_{3}}(\cfrac{\partial A_{0}Z}{\partial x})dxdydz
\\&=&g_{11}+ i\xi_{1} FI(\cfrac{\partial A_{0}Z}{\partial x})=g_{11}+ i\xi_{1}(g_{1}+ i\xi_{1} FI(A_{0}Z)),\end{eqnarray*}
\begin{eqnarray*}
FI(\cfrac{\partial^{2} A_{0}Z}{\partial x\partial y})&=&\int_{\overline{\Omega}}e^{-ix\xi_{1}-iy\xi_{2}-iz\xi_{3}}(\cfrac{\partial^{2} A_{0}Z}{\partial x\partial y})dxdydz
\\&=&\int_{\partial\Omega}(\cfrac{\partial A_{0}Z}{\partial y})n_{1}e^{-ix\xi_{1}-iy\xi_{2}-iz\xi_{3}}dS
+\\&&i\xi_{1}\int_{\overline{\Omega}}e^{-ix\xi_{1}-iy\xi_{2}-iz\xi_{3}}(\cfrac{\partial A_{0}Z}{\partial y})dxdydz
\\&=&g_{21}+ i\xi_{1} FI(\cfrac{\partial A_{0}Z}{\partial y})=g_{21}+ i\xi_{1}(g_{2}+ i\xi_{2} FI(A_{0}Z)),\end{eqnarray*}
\begin{eqnarray*}
FI(\cfrac{\partial^{2} A_{0}Z}{\partial x\partial z})&=&\int_{\overline{\Omega}}e^{-ix\xi_{1}-iy\xi_{2}-iz\xi_{3}}(\cfrac{\partial^{2} A_{0}Z}{\partial x\partial z})dxdydz
\\&=&\int_{\partial\Omega}(\cfrac{\partial A_{0}Z}{\partial z})n_{1}e^{-ix\xi_{1}-iy\xi_{2}-iz\xi_{3}}dS
+\\&&i\xi_{1}\int_{\overline{\Omega}}e^{-ix\xi_{1}-iy\xi_{2}-iz\xi_{3}}(\cfrac{\partial A_{0}Z}{\partial z})dxdydz
\\&=&g_{31}+ i\xi_{1} FI(\cfrac{\partial A_{0}Z}{\partial z})=g_{31}+ i\xi_{1}(g_{3}+ i\xi_{3} FI(A_{0}Z)),\end{eqnarray*}
\begin{eqnarray*}
FI(\cfrac{\partial^{2} A_{0}Z}{\partial y^{2}})&=&\int_{\overline{\Omega}}e^{-ix\xi_{1}-iy\xi_{2}-iz\xi_{3}}(\cfrac{\partial^{2} A_{0}Z}{\partial y^{2}})dxdydz
\\&=&\int_{\partial\Omega}(\cfrac{\partial A_{0}Z}{\partial y})n_{2}e^{-ix\xi_{1}-iy\xi_{2}-iz\xi_{3}}dS
+\\&&i\xi_{2}\int_{\overline{\Omega}}e^{-ix\xi_{1}-iy\xi_{2}-iz\xi_{3}}(\cfrac{\partial A_{0}Z}{\partial y})dxdydz
\\&=&g_{22}+ i\xi_{2} FI(\cfrac{\partial A_{0}Z}{\partial y})=g_{22}+ i\xi_{2}(g_{2}+ i\xi_{2} FI(A_{0}Z)),\end{eqnarray*}
\begin{eqnarray*}
FI(\cfrac{\partial^{2} A_{0}Z}{\partial y\partial z})&=&\int_{\overline{\Omega}}e^{-ix\xi_{1}-iy\xi_{2}-iz\xi_{3}}(\cfrac{\partial^{2} A_{0}Z}{\partial y\partial z})dxdydz
\\&=&\int_{\partial\Omega}(\cfrac{\partial A_{0}Z}{\partial z})n_{2}e^{-ix\xi_{1}-iy\xi_{2}-iz\xi_{3}}dS
+\\&&i\xi_{2}\int_{\overline{\Omega}}e^{-ix\xi_{1}-iy\xi_{2}-iz\xi_{3}}(\cfrac{\partial A_{0}Z}{\partial z})dxdydz
\\&=&g_{32}+ i\xi_{2} FI(\cfrac{\partial A_{0}Z}{\partial z})=g_{32}+ i\xi_{2}(g_{3}+ i\xi_{3} FI(A_{0}Z)),\end{eqnarray*}
\begin{eqnarray*}
FI(\cfrac{\partial^{2} A_{0}Z}{\partial z^{2}})&=&\int_{\overline{\Omega}}e^{-ix\xi_{1}-iy\xi_{2}-iz\xi_{3}}(\cfrac{\partial^{2} A_{0}Z}{\partial z^{2}})dxdydz
\\&=&\int_{\partial\Omega}(\cfrac{\partial A_{0}Z}{\partial z})n_{3}e^{-ix\xi_{1}-iy\xi_{2}-iz\xi_{3}}dS
+\\&&i\xi_{3}\int_{\overline{\Omega}}e^{-ix\xi_{1}-iy\xi_{2}-iz\xi_{3}}(\cfrac{\partial A_{0}Z}{\partial z})dxdydz
\\&=&g_{33}+ i\xi_{3} FI(\cfrac{\partial A_{0}Z}{\partial z})=g_{33}+ i\xi_{3}(g_{3}+ i\xi_{3} FI(A_{0}Z)),\end{eqnarray*}
where\begin{eqnarray*}
g_{1}&=&\int_{\partial \Omega}A_{1}n_{1}e^{-ix\xi_{1}-iy\xi_{2}-iz\xi_{3}}dS,\\
g_{2}&=&\int_{\partial \Omega}A_{1}n_{2}e^{-ix\xi_{1}-iy\xi_{2}-iz\xi_{3}}dS,\\
g_{3}&=&\int_{\partial \Omega}A_{1}n_{3}e^{-ix\xi_{1}-iy\xi_{2}-iz\xi_{3}}dS\\
g_{11}&=&\int_{\partial\Omega}A_{2}n_{1}e^{-ix\xi_{1}-iy\xi_{2}-iz\xi_{3}}dS,\\
g_{21}&=&\int_{\partial\Omega}A_{3}n_{1}e^{-ix\xi_{1}-iy\xi_{2}-iz\xi_{3}}dS,\\
g_{31}&=&\int_{\partial\Omega}A_{4}n_{1}e^{-ix\xi_{1}-iy\xi_{2}-iz\xi_{3}}dS,\\
g_{22}&=&\int_{\partial\Omega}A_{3}n_{2}e^{-ix\xi_{1}-iy\xi_{2}-iz\xi_{3}}dS,\end{eqnarray*}
 \begin{eqnarray*}
g_{32}&=&\int_{\partial\Omega}A_{4}n_{2}e^{-ix\xi_{1}-iy\xi_{2}-iz\xi_{3}}dS,\\
g_{33}&=&\int_{\partial\Omega}A_{4}n_{3}e^{-ix\xi_{1}-iy\xi_{2}-iz\xi_{3}}dS,\\A_{1}&=&u|_{\partial \Omega},\ A_{2}=u_{x}|_{\partial \Omega},\ A_{3}=u_{y}|_{\partial \Omega},\ A_{4}=u_{z}|_{\partial \Omega}.
\end{eqnarray*} Now we have transformed Eqs(3.6) to (3.14) into the following,\ $
BFI(Z)=\beta_{1}$, where\begin{eqnarray*}B&=&\left(
                                         \begin{array}{c}
                                         i\xi_{1}A_{0}-A_{01} \\
                                         i\xi_{2}A_{0}-A_{02} \\
                                         i\xi_{3}A_{0}-A_{03} \\
                                         (i\xi_{1})^{2}A_{0}-A_{11} \\
                                         i\xi_{1}i\xi_{2}A_{0}-A_{12} \\
                                         i\xi_{1}i\xi_{3}A_{0}-A_{13} \\
                                         (i\xi_{2})^{2}A_{0}-A_{22} \\
                                         i\xi_{2}i\xi_{3}A_{0}-A_{23} \\
                                         (i\xi_{3})^{2}A_{0}-A_{33} \\
                                         \end{array}
                                       \right)_{9m\times10m}
=(B_{1},\ -B_{2}),\\ \beta_{1}&=&(-g_{1}^{T},\ -g_{2}^{T},\ -g_{3}^{T},\ -g_{11}^{T}-i\xi_{1}g_{1}^{T},\ -g_{21}^{T}-i\xi_{1}g_{2}^{T},\ -g_{31}^{T}-i\xi_{1}g_{3}^{T},\\&& -g_{22}^{T}-i\xi_{2}g_{2}^{T},\ -g_{32}^{T}-i\xi_{2}g_{3}^{T},\ -g_{33}^{T}-i\xi_{3}g_{3}^{T})^{T},\end{eqnarray*} $B_{1}$\ is the first\ $9m$\ columns of\ $B$,\ $-B_{2}$\ is the last\ $m$\ columns of\ $B$.\\If we assume
$$ Z=\left(
       \begin{array}{c}
         Z_{1} \\
         Z_{2} \\
       \end{array}
     \right),\ \mbox{where}\ Z_{1}\ \mbox{is the first\ $9m$\ componenets of}\ Z,$$ then we obtain\ $ B_{1}FI(Z_{1})=\beta_{1}+B_{2}FI(Z_{2}).$\\
                             Now we determine the parameters. We choose\ $C_{k},\ 1\leq k\leq 9$, such that\ $det(B_{1})\neq 0$. Moreover, there exists\ $a_{1}$\ is a polynomial of\ $i\xi_{1},\ i\xi_{2},\ i\xi_{3},\ F^{-1}(a_{1}^{-1})$\ is in\ $S^{\prime}$\ and locally integrable, such that\ $a_{1}B_{1}^{-1}$\ is also a polynomial of\ $i\xi_{1},\ i\xi_{2},\ i\xi_{3}$, where\ $S^{\prime}$\ is the dual space of the Schwartz space.\\ We may select\ $a_{1}$\ is just\ $det(B_{1})$\ or the factor of it.\\In next section, we need two conditions more
                             \begin{eqnarray}&&\forall\ \varphi\in S,\ \varphi.\ast F^{-1}(a_{1}^{-1})\in C(R^{3}),\\&& (1+|\xi|^{2})^{-a-3}B_{1}^{-1},\ (1+|\xi|^{2})^{-a-3}B_{1}^{-1}B_{2}\in L^{1}(R^{3}),
                             \end{eqnarray}where\ $S$\ is the Schwartz space,\ $\xi=(\xi_{1},\ \xi_{2},\ \xi_{3})^{T}$,\ $a=\max\{\partial(a_{1}B^{-1}_{1}),\ \partial(a_{1}B^{-1}_{1}B_{2})\}$, here\ $\partial(\cdot)$\ means the highest degree.\\
Let's introduce three examples.\\
                             (1)\ $u$\ is resolved, Eqs(1.1) is just as follows,\[u-\sum_{k=1}^{9}C_{k}V_{1+k}-S=0,\]
                             where\ $(V_{2},\ V_{3},\cdots,\ V_{10})=(\partial u,\ \partial^{2}u)$.\\
                             We can obtain\ $B_{1}=\alpha_{0}\beta_{0}-E_{9m}$, where\[\alpha_{0}=(i\xi_{j}E,\ 1\leq j\leq 3,\ i\xi_{j}i\xi_{k}E,\ 1\leq j\leq k\leq 3)^{T},\] $\beta_{0}=(C_{1},\ C_{2},\cdots,\ C_{9})$,\ $E$\ is the\ $m$th order identity matrix,\ $E_{9m}$\ is the\ $9m$th order identity matrix.\\ From the results in linear Algebra,\ $\forall\ A\in R^{n_{1}\times n_{2}},\ B\in R^{n_{2}\times n_{1}},\ n_{1}\geq n_{2}$, we have
                              \[ det(AB-E_{n_{1}})=(-1)^{n_{1}-n_{2}}det(BA-E_{n_{2}}),\ (AB-E_{n_{1}})^{-1}=A(BA-E_{n_{2}})^{-1}B-E_{n_{1}}.\]
                              Hence, we can obtain\[ det(B_{1})=det(\beta_{0}\alpha_{0}-E),\ B_{1}^{-1}=\alpha_{0}(\beta_{0}\alpha_{0}-E)^{-1}\beta_{0}-E_{9m},\]
                              where\[\beta_{0}\alpha_{0}=\sum_{j=1}^{3}i\xi_{j}C_{j}+\sum_{1\leq j\leq k\leq3}i\xi_{j}i\xi_{k}C_{3+k+(6-j)(j-1)/2}.\]
                            Letting\ $C_{1}=-E$, the others are\ $0$, then\ $det(B_{1})\neq 0$, moreover, (3.16) and (3.17) hold. We denote it as\ $u_{x}+u$.\\
                             (2)\ $u_{x}$\ is resolved, Eqs(1.1) is just as follows,\[u_{x}-\sum_{k=1}^{9}C_{k}V_{1+k}-S=0,\]
                             where\ $(V_{2},\ V_{3},\cdots,\ V_{10})=(u,\ \partial u\setminus u_{x},\ \partial^{2}u)$.\\
                             We can obtain\begin{eqnarray*}B_{1}&=&\left(
                                         \begin{array}{c}
                                         i\xi_{1}E_{1}^{T}-(C_{1},\ C_{2},\cdots,\ C_{9}) \\
                                         i\xi_{2}E_{1}^{T}-E_{2}^{T} \\
                                         i\xi_{3}E_{1}^{T}-E_{3}^{T} \\
                                         (i\xi_{1})^{2}E_{1}^{T}-E_{4}^{T}\\
                                         i\xi_{1}i\xi_{2}E_{1}^{T}-E_{5}^{T} \\
                                         i\xi_{1}i\xi_{3}E_{1}^{T}-E_{6}^{T} \\
                                         (i\xi_{2})^{2}E_{1}^{T}-E_{7}^{T}\\
                                         i\xi_{2}i\xi_{3}E_{1}^{T}-E_{8}^{T} \\
                                         (i\xi_{3})^{2}E_{1}^{T}-E_{9}^{T} \\
                                         \end{array}
                                       \right)_{9m\times9m}.\end{eqnarray*} We can work out\ $det(B_{1})=det(i\xi_{1}E-C_{1}-(C_{2},\cdots,\ C_{9})\alpha_{1})$, where
                             \[\alpha_{1}=(i\xi_{j}E,\ 2\leq j\leq 3,\ i\xi_{j}i\xi_{k}E,\ 1\leq j\leq k\leq 3)^{T}.\]
                             Moreover,\begin{eqnarray*}B_{1}^{-1}&=&-\left(
                                         \begin{array}{c}
                                         R_{1} \\
                                         i\xi_{2}R_{1}-E_{2}^{T} \\
                                         i\xi_{3}R_{1}-E_{3}^{T} \\
                                         (i\xi_{1})^{2}R_{1}-E_{4}^{T}\\
                                         i\xi_{1}i\xi_{2}R_{1}-E_{5}^{T} \\
                                         i\xi_{1}i\xi_{3}R_{1}-E_{6}^{T} \\
                                         (i\xi_{2})^{2}R_{1}-E_{7}^{T}\\
                                         i\xi_{2}i\xi_{3}R_{1}-E_{8}^{T} \\
                                         (i\xi_{3})^{2}R_{1}-E_{9}^{T} \\
                                         \end{array}
                                       \right)_{9m\times9m},\end{eqnarray*} where\ $R_{1}=(-B_{1,\ 1}^{-1},\ B_{1,\ 1}^{-1}C_{2},\ B_{1,\ 1}^{-1}C_{3},\cdots,\ \ B_{1,\ 1}^{-1}C_{9}),\ B_{1,\ 1}=i\xi_{1}E-C_{1}-(C_{2},\cdots,\ C_{9})\alpha_{1}.$\\
                             Letting\ $C_{1}=-E$, the others are\ $0$, then\ $det(B_{1})\neq 0$, moreover, (3.16) and (3.17) hold. We also denote it as\ $u_{x}+u$.\\
                             (3)\ $u_{xx}$\ is resolved, Eqs(1.1) is just as follows,\[u_{xx}-\sum_{k=1}^{9}C_{k}V_{1+k}-s=0,\]
                             where\ $(V_{2},\ V_{3},\cdots,\ V_{10})=(u,\ \partial u,\ \partial^{2}u\setminus u_{xx})$.\\
                             We can obtain\begin{eqnarray*}B_{1}&=&\left(
                                         \begin{array}{c}
                                         (i\xi_{1})^{2}E_{1}^{T}-(C_{1},\ C_{2},\cdots,\ C_{9}) \\
                                         i\xi_{1}E_{1}^{T}-E_{2}^{T} \\
                                         i\xi_{2}E_{1}^{T}-E_{3}^{T} \\
                                         i\xi_{3}E_{1}^{T}-E_{4}^{T} \\
                                         i\xi_{1}i\xi_{2}E_{1}^{T}-E_{5}^{T} \\
                                         i\xi_{1}i\xi_{3}E_{1}^{T}-E_{6}^{T} \\
                                         (i\xi_{2})^{2}E_{1}^{T}-E_{7}^{T}\\
                                         i\xi_{2}i\xi_{3}E_{1}^{T}-E_{8}^{T} \\
                                         (i\xi_{3})^{2}E_{1}^{T}-E_{9}^{T} \\
                                         \end{array}
                                       \right)_{9m\times9m}.\end{eqnarray*} We can work out\ $det(B_{1})=det((i\xi_{1})^{2}E-C_{1}-(C_{2},\cdots,\ C_{9})\alpha_{2})$, where
                             \[\alpha_{2}=(i\xi_{j}E,\ 1\leq j\leq 3,\ i\xi_{j}i\xi_{k}E,\ 1\leq j\leq k\leq 3,\ (j,\ k)\neq(1,\ 1))^{T}.\]
                             Moreover\begin{eqnarray*}B_{1}^{-1}&=&-\left(
                                         \begin{array}{c}
                                         R_{1} \\
                                         i\xi_{1}R_{1}-E_{2}^{T} \\
                                         i\xi_{2}R_{1}-E_{3}^{T} \\
                                         i\xi_{3}R_{1}-E_{4}^{T} \\
                                         i\xi_{1}i\xi_{2}R_{1}-E_{5}^{T} \\
                                         i\xi_{1}i\xi_{3}R_{1}-E_{6}^{T} \\
                                         (i\xi_{2})^{2}R_{1}-E_{7}^{T}\\
                                         i\xi_{2}i\xi_{3}R_{1}-E_{8}^{T} \\
                                         (i\xi_{3})^{2}R_{1}-E_{9}^{T} \\
                                         \end{array}
                                       \right)_{9m\times9m},\end{eqnarray*} where\ $R_{1}=(-B_{1,\ 2}^{-1},\ B_{1,\ 2}^{-1}C_{2},\ B_{1,\ 2}^{-1}C_{3},\cdots,\ \ B_{1,\ 2}^{-1}C_{9}),\ B_{1,\ 2}=(i\xi_{1})^{2}E-C_{1}-(C_{2},\cdots,\ C_{9})\alpha_{2}$.\\
                             Letting\ $C_{7}=C_{9}=-E$, the others are\ $0$, then\ $det(B_{1})\neq 0$, moreover, (3.16) and (3.17) hold. We denote it as\ $\triangle u$.\\
In fact,\ $B_{1}$\ is in the form of\ $\alpha_{0}A_{0}(B_{1})-A(B_{1})$, where\ $A_{0}(B_{1})$\ is the first\ $9m$\ columns of\ $A_{0},\ A(B_{1})=\alpha_{0}A_{0}(B_{1})-B_{1}$. There is at least\ $E_{8m}$\ in\ $A(B_{1})$. \\With the parameters that\ $A(B_{1})$\ is invertible, we have
                             \begin{eqnarray}det(B_{1})&=&(-1)^{m}det(A(B_{1}))det(A^{\prime}_{0}(B_{1})\alpha_{0}-E),\\ B_{1}^{-1}&=&(A(B_{1}))^{-1}(\alpha_{0}(A^{\prime}_{0}(B_{1})\alpha_{0}-E)^{-1}A^{\prime}_{0}(B_{1})-E_{9m}),\end{eqnarray}
                             where\ $A^{\prime}_{0}(B_{1})=A_{0}(B_{1})(A(B_{1}))^{-1}$.\\
                             Assuming that\ $A^{\prime}_{0}(B_{1})=(C^{\prime}_{1},\ C^{\prime}_{2},\cdots,\ C^{\prime}_{9})$, where\ $C^{\prime}_{k}\in R^{m\times m},\ 1\leq k\leq 9$, then we obtain \[A^{\prime}_{0}(B_{1})\alpha_{0}=\sum_{j=1}^{3}i\xi_{j}C^{\prime}_{j}+\sum_{1\leq j\leq k\leq3}i\xi_{j}i\xi_{k}C^{\prime}_{3+k+(6-j)(j-1)/2}.\]
                             From\ $A_{0}(B_{1})=(C^{\prime}_{1},\ C^{\prime}_{2},\cdots,\ C^{\prime}_{9})A(B_{1})$, we have
                             \[(E,\ -C^{\prime}_{1},\ -C^{\prime}_{2},\cdots,\ -C^{\prime}_{9})\left(\begin{array}{c}
                             A_{0}(B_{1})\\
                              A(B_{1})\\\end{array}
                                       \right)_{10m\times 9m}=0.\]Because the rows of\ $A_{0},\ A_{01},\ A_{02},\ A_{03},\ A_{11},\ A_{12},\ A_{13},\ A_{22},\ A_{23},\ A_{33}$\ is a permutation of the rows of\ $\beta^{T},\ E_{j}^{T},\ 1\leq j\leq 9$, there exists a permutation matrix\ $P_{10m\times 10m}$, such that$$ P\left(\begin{array}{c}
                             A_{0}(B_{1})\\
                              A(B_{1})\\\end{array}
                                       \right)=\left(
          \begin{array}{cccccccccc}
            C_{1},&C_{2}, & C_{3}, &C_{4}, & C_{5}, & C_{6}, &C_{7}, &C_{8}, & C_{9}\\
            E, & 0, & 0, & 0, & 0,& 0, & 0, & 0, & 0\\
           0, & E, & 0, & 0, & 0,& 0, & 0, & 0, & 0 \\
            0, & 0, & E, & 0, & 0,& 0, & 0, & 0, & 0 \\
           0, & 0, & 0, & E, & 0,& 0, & 0, & 0, & 0 \\
            0, & 0, & 0, & 0, & E,& 0, & 0, & 0, & 0 \\
            0, & 0, & 0, & 0, & 0,& E, & 0, & 0, & 0 \\
            0, & 0, & 0, & 0, & 0,& 0, &E, & 0, & 0 \\
            0, & 0, & 0, & 0, & 0,&0, &0, & E, & 0 \\
            0,& 0, & 0, & 0, & 0, &  0, & 0, & 0, & E \\
          \end{array}
        \right).$$Assuming that\[ (E,\ -C^{\prime}_{1},\ -C^{\prime}_{2},\cdots,\ -C^{\prime}_{9})P^{-1}=(C^{\prime\prime}_{1},\ C^{\prime\prime}_{2},\cdots,\ C^{\prime\prime}_{10}),\]we can obtain\[ C^{\prime\prime}_{1}C_{k}+C^{\prime\prime}_{k+1}=0,\ 1\leq k\leq 9.\] We can work out\ $C_{k},\ 1\leq k\leq 9$, if\ $C^{\prime\prime}_{1}$\ is invertible, where\ $C^{\prime\prime}_{1}=(E,\ -C^{\prime}_{1},\ -C^{\prime}_{2},\cdots,\ -C^{\prime}_{9})P_{1},\ P_{1}$\ is the first\ $m$ columns of\ $P^{-1}$.\ From (3.31) and\ $P^{-1}=P^{T}$, we know that\ $C^{\prime\prime}_{k}$\ is composed of\ $m$\ different columns of\ $(E,\ -C^{\prime}_{1},\ -C^{\prime}_{2},\cdots,\ -C^{\prime}_{9})$,\ $1\leq k\leq 10$.\\Assuming that $u_{1},\ u_{2},\cdots,\ u_{r},\ 0\leq r\leq m$, are in the left hand sides of the Eqs.(1.1), moreover,\ $u_{r+1},\ u_{r+2},\cdots,\ u_{m}$, are in the beginning of\ $Z_{1}$, then there exists another permutation matrix\ $Q_{9m\times 9m}$, such that
        \[ QA(B_{1})=\left(\begin{array}{cc}
         C_{1,\ m-r},&C^{\ast}_{2}\\
         0_{(8m+r)\times (m-r)},& E_{8m+r}\end{array}
        \right),\]where\ $C_{1,\ m-r}$\ is the block of\ $C_{1}$\ which is composed of the elements at the intersections of the last\ $m-r$\ rows and the first\ $m-r$\ columns, moreover\ $(C_{1,\ m-r},\ C^{\ast}_{2})$\ is the last\ $m-r$\ rows of the\ $(C_{1},\ C_{2},\cdots,\ C_{9})$. We can see that\ $A(B_{1})$\ is invertible, if\ $C_{1,\ m-r}$\ is invertible.\\
        From (3.32), we know that\[C_{1,\ m-r}=-[(C^{\prime\prime}_{1})^{-1}(m-r)][C^{\prime\prime}_{2}(m-r)],\]where\ $(C^{\prime\prime}_{1})^{-1}(m-r)$\ is composed of the last\ $m-r$\ rows of the\ $(C^{\prime\prime}_{1})^{-1}$,\ $C^{\prime\prime}_{2}(m-r)$\ is composed of the first\ $m-r$\ columns of the\ $C^{\prime\prime}_{2}$.\\
         It's obviously that the set of\ $(C^{\prime}_{1},\ C^{\prime}_{2},\cdots,\ C^{\prime}_{9})$\ which satisfies (3.16) and (3.17) has the inner point.\\ From the perturbation theory, we can select\ $(C^{\prime}_{1},\ C^{\prime}_{2},\cdots,\ C^{\prime}_{9})$\ which satisfies (3.16) and (3.17), moreover,\ $$C^{\prime\prime}_{1},\ [(C^{\prime\prime}_{1})^{-1}(m-r)][C^{\prime\prime}_{2}(m-r)],$$ are invertible.\\
        Hence, we can determine the parameters\ $C_{k},\ 1\leq k\leq 9$, satisfying\ $det(B_{1})\neq 0$, moreover, (3.16) and (3.17) hold.\\
Assuming that\ $C_{0}=\{\xi|det(B_{1})=0\}$, then the measure of\ $C_{0}$\ is\ $0$, and we obtain
 \[FI(Z_{1})(1-I_{C_{0}}(\xi))=B_{1}^{-1}\beta_{1}(1-I_{C_{0}}(\xi))+B_{1}^{-1}B_{2}FI(Z_{2})(1-I_{C_{0}}(\xi)).\]
We want to do the inverse Fourier transform to both sides of (3.35). We need the following lemmas.
\begin{lemma} \label{lemma1}(Plancherel) If\ $f(x,\ y,\ z)\in L^{2}(R^{3})$, then\ $F(f(x,\ y,\ z))$\ exists, moreover\\(1)$\parallel F(f(x,\ y,\ z))\parallel_{L^{2}}=\parallel f(x,\ y,\ z)\parallel_{L^{2}}$,\\(2)$F^{-1}[F(f(x,\ y,\ z))]=f(x,\ y,\ z)$.\end{lemma}
We can see the proof on page 194 of [7].
\begin{lemma} \label{lemma1} If\ $f(x,\ y,\ z)\in L^{2}(R^{3})$,\ $C_{0}\subset R^{3}$, and the measure of\ $C_{0}$\ is\ $0$, then \[F^{-1}([F(f(x,\ y,\ z))](1-I_{C_{0}}(\xi)))=f(x,\ y,\ z).\]\end{lemma}
{\it Proof of lemma 3.2}. From the lemma 3.1, we know\ $F(f(x,\ y,\ z))\in L^{2}(R^{3})$. Therefore,
$$\int_{C_{0}}F(f(x,\ y,\ z))e^{ix\xi_{1}+iy\xi_{2}+iz\xi_{3}}d\xi_{1}d\xi_{2}d\xi_{3}=0.$$We obtain
$$ F^{-1}([F(f(x,\ y,\ z))](1-I_{C_{0}}(\xi)))=F^{-1}[F(f(x,\ y,\ z))]=f(x,\ y,\ z),$$which proves the statement.\qed\\
From these two lemmas, we obtain\[F^{-1}[FI(Z_{1})(1-I_{C_{0}}(\xi))]=Z_{1}I_{\overline{\Omega}}.\] This means that\ $B_{1}^{-1}\beta_{1}(1-I_{C_{0}}(\xi))+B_{1}^{-1}B_{2}FI(Z_{2})(1-I_{C_{0}}(\xi))$\ is the Fourier transform of continuous functions on\ $\overline{\Omega}$. We have the Paley-Wiener-Schwartz theorem from [1] as follows,
\begin{lemma} \label{lemma1}(Paley-Wiener-Schwartz) Let\ $K$\ be a convex compact set of\ $R^{n}$\ with support function\ $H(\xi)=\sup_{x\in K}<x,\ \xi>,\ \forall \xi\in R^{n}$. If\ $u$\ is a distribution with support contained in\ $K$,\ then there exists\ $C>0,\ N$\ is a positive whole number, such that\[|F(u)(\zeta)|\leq C(1+|\zeta|)^{N}e^{H(Im \zeta)}, \forall\ \zeta\in C^{n}.\]
Conversely, every entire analytic function in\ $C^{n}$\ satisfying an estimate of the form (3.38) is the Fourier-Laplace transform of a distribution with support contained in\ $K$.\end{lemma}We may read the proof on pages 181 to 182 in [1]. H$\ddot{o}$rmander theorem from [1] is as follows,
\begin{lemma} \label{lemma1}(H$\ddot{o}$rmander) If\ $F(u)(\zeta)$\ is an entire analytic function in\ $C^{n}$\ satisfying an estimate of the form (3.38),\ $p(\zeta)$\ is a polynomial,\ $F(u)(\zeta)/p(\zeta)$\ is an entire function, then\ $F(u)(\zeta)/p(\zeta)$\ satisfies an estimate of the form (3.38), too.\end{lemma}We can see the proof on page 183 of [1]. \\
From (3.37) we know\ $B_{1}^{-1}\beta_{1}(1-I_{C}(\xi))+B_{1}^{-1}B_{2}FI(Z_{2})(1-I_{C}(\xi))$\ is the Fourier transform of continuous functions on\ $\overline{\Omega}$. Because\ $a_{1}B_{1}^{-1}$\ is entire, from lemma 2.2, Paley-Wiener-Schwartz theorem and H$\ddot{o}$rmander theorem in [1], we know that \[B_{1}^{-1}\beta_{1}+B_{1}^{-1}B_{2}FI(Z_{2})\] should be entire, too. \\
Under this condition, we can obtain$$Z_{1}I_{\overline{\Omega}}=F^{-1}(B_{1}^{-1}\beta_{1}+B_{1}^{-1}B_{2}FI(Z_{2})).$$We can also work out the following, \begin{eqnarray*}B_{1}^{-1}\beta_{1}&=&(a_{1}^{-1})(a_{1}B_{1}^{-1}\beta_{1}),\\
B_{1}^{-1}B_{2}&=&(a_{1}^{-1})(a_{1}B_{1}^{-1}B_{2}).\end{eqnarray*}
We find that\ $F^{-1}(a_{1}^{-1})\in S^{\prime},\ a_{1}B_{1}^{-1}\beta_{1}\ \mbox{and}\ a_{1}B_{1}^{-1}B_{2}$\ satisfy the Paley-Wiener-Schwartz form (3.38). Hence\ $F^{-1}(a_{1}B_{1}^{-1}\beta_{1})\in \varepsilon^{\prime},\ F^{-1}(a_{1}B_{1}^{-1}B_{2})\in \varepsilon^{\prime}$, and their compact support is contained in\ $\overline{\Omega}$, where\ $S^{\prime}$\ is the dual space of the Schwartz space, and\ $\varepsilon^{\prime}$\ is the dual space of\ $C^{\infty}(R^{3})$. We need a lemma as follows.
\begin{lemma} \label{lemma1} If\ $v_{1}\in S^{\prime},\ v_{2}\in \varepsilon^{\prime}$, it follows that\ $v_{1}\ast v_{2}\in S^{\prime}$\ and that
$$F(v_{1}\ast v_{2})=F(v_{1})F(v_{2}).$$\end{lemma} We can see the proof on page 166 in [1]. \\From this lemma we know\begin{eqnarray*}F^{-1}(B_{1}^{-1}\beta_{1})&=&F^{-1}(a_{1}^{-1}).\ast F^{-1}(a_{1}B_{1}^{-1}\beta_{1}),\\ F^{-1}(B_{1}^{-1}B_{2})&=&F^{-1}(a_{1}^{-1}).\ast F^{-1}(a_{1}B_{1}^{-1}B_{2}),\end{eqnarray*} all exist and are in\ $S^{\prime}$, where\ $.\ast$\ is the matrix convolution.\\
If we assume\begin{eqnarray*}w_{1}(x,\ y,\ z)&=&F^{-1}(B_{1}^{-1}\beta_{1}),\\ w_{2}(x,\ y,\ z)&=&F^{-1}(B_{1}^{-1}B_{2}),\end{eqnarray*}
then we obtain\[Z_{1}I_{\overline{\Omega}}=w_{1}(x,\ y,\ z)+w_{2}(x,\ y,\ z).\ast Z_{2}I_{\overline{\Omega}},\]
where\begin{eqnarray*}Z_{1}I_{\overline{\Omega}}&=&(E_{1}^{T}ZI_{\overline{\Omega}},\ E_{2}^{T}ZI_{\overline{\Omega}},\ E_{3}^{T}ZI_{\overline{\Omega}},\ E_{4}^{T}ZI_{\overline{\Omega}},\ E_{5}^{T}ZI_{\overline{\Omega}},\\&& E_{6}^{T}ZI_{\overline{\Omega}},\ E_{7}^{T}ZI_{\overline{\Omega}},\ E_{8}^{T}ZI_{\overline{\Omega}},\ E_{9}^{T}ZI_{\overline{\Omega}})^{T},\\Z_{2}I_{\overline{\Omega}}&=&E_{10}^{T}ZI_{\overline{\Omega}}=f(E_{1}^{T}ZI_{\overline{\Omega}},\ E_{2}^{T}ZI_{\overline{\Omega}},\ E_{3}^{T}ZI_{\overline{\Omega}},\ E_{4}^{T}ZI_{\overline{\Omega}},\\&& E_{5}^{T}ZI_{\overline{\Omega}},\ E_{6}^{T}ZI_{\overline{\Omega}},\ E_{7}^{T}ZI_{\overline{\Omega}},\ E_{8}^{T}ZI_{\overline{\Omega}},\ E_{9}^{T}ZI_{\overline{\Omega}},\\&&  x,\ y,\ z)-\sum_{l=1}^{9}C_{l}E_{l}^{T}Z.\end{eqnarray*}It is obvious\ $\exists\ \psi$, such that\ $Z_{2}I_{\overline{\Omega}}=\psi(Z_{1}I_{\overline{\Omega}})$. Therefore, we attain
\[ Z_{1}I_{\overline{\Omega}}=w_{1}(x,\ y,\ z)+w_{2}(x,\ y,\ z).\ast(\psi(Z_{1}I_{\overline{\Omega}})).\]
Now we arrive at the second equivalent result as follows,\begin{theorem} \label{Theorem3-2}\ $ w_1,\ w_2,\ \psi,$\ as we described, then Eqs(1.1) is equivalent to Eqs(3.41). \end{theorem}
 {\it Proof of theorem 3.2}. If\ $u$\ satisfies Eqs(1.1), then from the theorem 3.1,\ $Z=(V_{2},\ V_{3},\ \cdots,\ V_{10},\ S)^{T}$\ satisfies Eq.(3.6) to Eq.(3.15). Hence we can obtain the following by Fourier transform on\ $\overline{\Omega}$,
 $$B FI(Z)=\beta_{1},\ B_{1}FI(Z_{1})=\beta_{1}+B_{2}FI(Z_{2}),\ FI(Z_{1})=B_{1}^{-1}\beta_{1}+B_{1}^{-1}B_{2}FI(Z_{2}).$$ After we do the inverse Fourier transform, we obtain\ $Z_{1}=S_{1}(u)=(V_{2},\ V_{3},\ \cdots,\ V_{10})^{T}$\ satisfies Eqs(3.41). \\If\ $Z_{1}$\ satisfies Eqs(3.41), then letting\ $Z_{2}=\psi(Z_{1}),\ Z=(Z_{1},\ Z_{2})^{T}$,\ we obtain the following by the Fourier transform,
 $$FI(Z_{1})=B_{1}^{-1}\beta_{1}+B_{1}^{-1}B_{2}FI(Z_{2}),\ B_{1}FI(Z_{1})=\beta_{1}+B_{2}FI(Z_{2}),\ B FI(Z)=\beta_{1}.$$ After we do the inverse Fourier transform, we obtain\ $Z$\ satisfies Eq.(3.6) to Eq.(3.15) on\ $\overline{\Omega}$. From the theorem 3.1, we know\ $u=S_{2}(Z_{1})=A_{0}(Z_{1},\ \psi(Z_{1}))^{T}$\ is the solution of Eqs(1.1).\\ Obviously\ $S_{1},\ S_{2}$\ are continuous. Moreover\ $S_{1}(S_{2}( Z_{1}))=Z_{1},\ S_{2}(S_{1}(u))=u.$\ From definition 1.1, we know the the statement stands.\qed\\We denote Eqs(3.41) as\ $Z_{1}=T_{0}(Z_{1})$, where\ $T_{0}(Z_{1})=w_{1}+w_{2}.\ast\psi(Z_{1})$. We will discuss the continuous fix-point of it.\\
Next we determinate all the boundary conditions.\\
If the component of\ $u_{xx},\ u_{yy},\ u_{zz}$\ is in the left hand side of Eqs(1.1), then it is continuous function of\ $Z_{1}$. Else if the component of\ $u_{xx},\ u_{yy},\ u_{zz}$\ is in the right hand side of Eqs(1.1), then it is component of\ $Z_{1}$. Hence, there exists continuous function\ $\psi_{1}$, such that
\[ -\triangle u=\psi_{1}(Z_{1}).\]After we do Fourier transform on\ $\overline{\Omega}$, we can get the following,
\[-(i\xi_{1})^{2}-(i\xi_{2})^{2}-(i\xi_{3})^{2})FI(u)=g_{11}+g_{22}+g_{33}+i\xi_{1}g_{1}+i\xi_{2}g_{2}+i\xi_{3}g_{3}+FI(\psi_{1}(Z_{1})).\]
We get that\[(-(i\xi_{1})^{2}-(i\xi_{2})^{2}-(i\xi_{3})^{2})^{-1}(g_{11}+g_{22}+g_{33}+i\xi_{1}g_{1}+i\xi_{2}g_{2}+i\xi_{3}g_{3}+FI(\psi_{1}(Z_{1})))\]should be entire. Under this condition, if we assume\ $a_{2}=-(i\xi_{1})^{2}-(i\xi_{2})^{2}-(i\xi_{3})^{2}$, then we obtain that
\begin{eqnarray*}uI_{\overline{\Omega}}&=&F^{-1}[a_{2}^{-1}(g_{11}+g_{22}+g_{33}+i\xi_{1}g_{1}+i\xi_{2}g_{2}+i\xi_{3}g_{3}+FI(\psi_{1}(Z_{1})(x,\ y,\ z)))]\ \mbox{a.e.}\\F^{-1}[a_{2}^{-1}]&=&h(x,\ y,\ z)=\cfrac{1}{4\pi}\ \cfrac{1}{\sqrt{x^{2}+y^{2}+z^{2}}},
        \end{eqnarray*}where\ $a.e.$\ means almost everywhere,
         \begin{eqnarray*}F^{-1}[a_{2}^{-1}(g_{11}+g_{22}+g_{33})]&=&\int_{\partial\Omega}h(x-x_{1},\ y-y_{1},\ z-z_{1})A_{5}(x_{1},\ y_{1},\ z_{1})dS,\ A_{5}=\cfrac{\partial u}{\partial n}|_{\partial\Omega},\end{eqnarray*}
         $$F^{-1}[a_{2}^{-1}(i\xi_{1}g_{1}+i\xi_{2}g_{2}+i\xi_{3}g_{3})]=\int_{\partial\Omega}\cfrac{\partial h(x-x_{1},\ y-y_{1},\ z-z_{1})}{\partial n(x_{1},\ y_{1},\ z_{1})}A_{1}(x_{1},\ y_{1},\ z_{1})dS,$$
         $$F^{-1}[a_{2}^{-1}FI(\psi_{1}(Z_{1})(x,\ y,\ z))]=\int_{\overline{\Omega}} h(x-x_{1},\ y-y_{1},\ z-z_{1})\psi_{1}(Z_{1})(x_{1},\ y_{1},\ z_{1})dx_{1}dy_{1}dz_{1},$$
         and\begin{eqnarray*} \cfrac{\partial h(x-x_{1},\ y-y_{1},\ z-z_{1})}{\partial n(x_{1},\ y_{1},\ z_{1})}&=&\cfrac{\partial h(x-x_{1},\ y-y_{1},\ z-z_{1})}{\partial x}n_{1}(x_{1},\ y_{1},\ z_{1})+\\&&\cfrac{\partial h(x-x_{1},\ y-y_{1},\ z-z_{1})}{\partial y} n_{2}(x_{1},\ y_{1},\ z_{1})+\\&&\cfrac{\partial h(x-x_{1},\ y-y_{1},\ z-z_{1})}{\partial z}n_{3}(x_{1},\ y_{1},\ z_{1}).\end{eqnarray*}
         We denote it in an easy way as follows,\begin{eqnarray*}F^{-1}[a_{2}^{-1}(g_{11}+g_{22}+g_{33})]&=&h.\ast_{\partial\Omega}A_{5},
         \\F^{-1}[a_{2}^{-1}(i\xi_{1}g_{1}+i\xi_{2}g_{2}+i\xi_{3}g_{3})]&=&\cfrac{\partial h}{\partial n_{p_{1}}}.\ast_{\partial\Omega}A_{1},\\F^{-1}[a_{2}^{-1}FI(\psi_{1}(Z_{1})(x,\ y,\ z))]&=&h.\ast \psi_{1}(Z_{1}),\end{eqnarray*}where\ $P_{1}=(x_{1},\ y_{1},\ z_{1})^{T}$, then we get the following,\[uI_{\overline{\Omega}}=h.\ast_{\partial\Omega}A_{5}+\cfrac{\partial h}{\partial n_{p_{1}}}.\ast_{\partial\Omega}A_{1}+h.\ast \psi_{1}(Z_{1}).\ \mbox{a.e.}\]Maybe you will say that all items in right hand side of (3.45) are continuous on\ $\overline{\Omega}$. We don't need\ $a.e.$. You are right except\ $\partial h/\partial n_{p_{1}}.\ast_{\partial\Omega}A_{1}$. We say a little more on it.\begin{eqnarray*}\cfrac{\partial h(M-P_{1})}{\partial n_{p_{1}}} &=&\nabla_{M}h(M-P_{1})\cdot n_{p_{1}}\\
         &=&-\nabla_{P_{1}}h(M-P_{1})\cdot n_{p_{1}},\end{eqnarray*}
         where\ $M=(x,\ y,\ z)^{T}$, and\begin{eqnarray*}\nabla_{M}h(M-P_{1})&=&(\cfrac{\partial h(M-P_{1})}{\partial x},\ \cfrac{\partial h(M-P_{1})}{\partial y},\ \cfrac{\partial h(M-P_{1})}{\partial z})^{T},\\ \nabla_{P_{1}}h(M-P_{1})&=&(\cfrac{\partial h(M-P_{1})}{\partial x_{1}},\ \cfrac{\partial h(M-P_{1})}{\partial y_{1}},\ \cfrac{\partial h(M-P_{1})}{\partial z_{1}})^{T}.\end{eqnarray*}We see\ $\nabla_{P_{1}}h(M-P_{1})\cdot n_{p_{1}}$\ in some books. That will increase a sign of minus.\\
         From Corollary 2.1 and Theorem 2.3 in last section, we have two results as follows,
         \begin{lemma} \label{lamma1} If\ $\Omega$\ is bounded, $\partial\Omega\in C^{1,\ \beta},\ 0<\beta\leq 1,\ A_{1}(P_{1})\in C(\partial\Omega)$, then\ $\partial h/\partial n_{p_{1}}.\ast_{\partial\Omega}A_{1}$\ is continuous on\ $R^{3}\setminus \partial\Omega$\ and\ $\partial\Omega$, but\ $\forall P_{0}\in \partial\Omega$, we have\begin{eqnarray}&& \lim_{M\rightarrow P_{0}+}(\cfrac{\partial h}{\partial n_{p_{1}}}.\ast_{\partial\Omega}A_{1})(M)=(\cfrac{\partial h}{\partial n_{p_{1}}}.\ast_{\partial\Omega}A_{1})(P_{0})+\cfrac{1}{2} A_{1}(P_{0}),\\&&  \lim_{M\rightarrow P_{0}-}(\cfrac{\partial h}{\partial n_{p_{1}}}.\ast_{\partial\Omega}A_{1})(M)=(\cfrac{\partial h}{\partial n_{p_{1}}}.\ast_{\partial\Omega}A_{1})(P_{0})-\cfrac{1}{2} A_{1}(P_{0}),\end{eqnarray}where$$(\cfrac{\partial h}{\partial n_{p_{1}}}.\ast_{\partial\Omega}A_{1})(M)=\int_{\partial\Omega}[\nabla_{M}h(M-P_{1})\cdot n_{p_{1}}]A_{1}(P_{1}) dS_{P_{1}}, $$
         $M\rightarrow P_{0}+$\ means\ $M$\ is near to\ $P_{0}$\ from interior of\ $\Omega$\ and\ $M\rightarrow P_{0}-$\ means\ $M$\ is near to\ $P_{0}$\ from exterior of\ $\Omega$. \end{lemma}
\begin{lemma}\label{lamma1} If\ $\Omega$\ is bounded, $\partial\Omega\in C^{1,\ \beta},\ 0<\beta\leq 1,\ A_{1}(P_{1})\in C(\partial\Omega)$, then\ $\partial h/\partial n_{p_{0}}.\ast_{\partial\Omega}A_{1}$\ is continuous on\ $R^{3}\setminus \partial\Omega$\ and\ $\partial\Omega$, but\ $\forall P_{0}\in \partial\Omega$, we have\begin{eqnarray}&& \lim_{M\rightarrow n_{p_{0}}^{+}}(\cfrac{\partial h}{\partial n_{p_{0}}}.\ast_{\partial\Omega}A_{1})(M)=(\cfrac{\partial h}{\partial n_{p_{0}}}.\ast_{\partial\Omega}A_{1})(P_{0})-\cfrac{1}{2} A_{1}(P_{0}),\\&&  \lim_{M\rightarrow n_{p_{0}}^{-}}(\cfrac{\partial h}{\partial n_{p_{0}}}.\ast_{\partial\Omega}A_{1})(M)=(\cfrac{\partial h}{\partial n_{p_{0}}}.\ast_{\partial\Omega}A_{1})(P_{0})+\cfrac{1}{2} A_{1}(P_{0}),\end{eqnarray}
 where$$(\cfrac{\partial h}{\partial n_{p_{0}}}.\ast_{\partial\Omega}A_{1})(M)=\int_{\partial\Omega}[\nabla_{M}h(M-P_{1})\cdot n_{p_{0}}]A_{1}(P_{1}) dS_{P_{1}}, $$
 $M\rightarrow n_{p_{0}}^{+}$\ means\ $M$\ is near to\ $P_{0}$\ along\ $n_{P_{0}}$\ from exterior of\ $\Omega$\ and\ $M\rightarrow n_{p_{0}}^{-}$\ means\ $M$\ is near to\ $P_{0}$\ along\ $n_{P_{0}}$\ from interior of\ $\Omega$.\end{lemma}
 And from page 54 to page 56 in [2], we may get another result as the following,
\begin{lemma}\label{lamma1} If\ $\Omega$\ is bounded, $\partial\Omega\in C^{1,\ \beta},\ 0<\beta\leq 1,\ f$\ is continuous, then\ $w=h.\ast f\in C^{1}(\overline{\Omega})$.
\end{lemma}
Now we can see that \[uI_{\Omega}=[h.\ast_{\partial\Omega}A_{5}+\cfrac{\partial h}{\partial n_{p_{1}}}.\ast_{\partial\Omega}A_{1}+h.\ast \psi_{1}(Z_{1})]I_{\Omega}\in C^{1}(\Omega).\]
We find there is only\ $A_{1}$\ is known in (3.50). We should determinate\ $A_{5}$. Following the Lyapunov's potential theory from page 173 to page 201 in [6], we discuss it as follows.\\
If\ $A_{1}$\ is known, then from (3.50) we can get the following,\ $\forall M\in \Omega$, \[ \cfrac{\partial u(M)I_{\Omega}}{\partial n_{p_{0}}}=(\cfrac{\partial g_{1}}{\partial n_{p_{0}}}+\cfrac{\partial h(M)}{\partial n_{p_{0}}}.\ast_{\partial\Omega}A_{5}+\cfrac{\partial h(M)}{\partial n_{p_{0}}}.\ast \psi_{1}(Z_{1}))I_{\Omega},\] where$$ g_{01}=\cfrac{\partial h}{\partial n_{p_{1}}}.\ast_{\partial\Omega}A_{1}.$$
Because\ $u(M)\in C^{2}(\overline{\Omega})$, and from lemma 3.6, we can get as follows,
$$ \lim_{M\rightarrow n_{p_{0}}^{-}} \cfrac{\partial u(M)I_{\Omega}}{\partial n_{p_{0}}}=A_{5},\ \lim_{M\rightarrow n_{p_{0}}^{-}}\cfrac{\partial h(M)}{\partial n_{p_{0}}}.\ast_{\partial\Omega}A_{5}=\cfrac{\partial h(P_{0})}{\partial n_{p_{0}}}.\ast_{\partial\Omega}A_{5}+\cfrac{1}{2}A_{5}. $$ Hence, there exist\ $g_{02}\in C^{1}(\partial\Omega)$, such that
\[  \lim_{M\rightarrow n_{p_{0}}^{-}} \cfrac{\partial g_{01}(M)}{\partial n_{p_{0}}}=g_{02}(P_{0}),\ \forall P_{0}\in \partial\Omega.\] We can get the second type linear Fredholm integral equations which\ $A_{5}$\ should satisfy,
\[ \cfrac{1}{2}A_{5}=\cfrac{\partial h(P_{0})}{\partial n_{p_{0}}}.\ast_{\partial\Omega}A_{5}+g_{2}+\cfrac{\partial h(P_{0})}{\partial n_{p_{0}}}.\ast \psi_{1}(Z_{1}).\]
We will prove that there is only\ $0$\ for the homogeneous equations as follows, \[ \cfrac{1}{2}A_{5}=\cfrac{\partial h(P_{0})}{\partial n_{p_{0}}}.\ast_{\partial\Omega}A_{5}.\] If\ $A_{5}$\ satisfies (3.54), then we may let\[ W(M)=h(M).\ast_{\partial\Omega}A_{5}=\int_{\partial\Omega}h(M-P_{1})A_{5}(P_{1})dS_{P_{1}},\]
and we can obtain that\ $W$\ is continuous on\ $R^{3}$, moreover
\[-\triangle W=0,\ \mbox{on}\ (R^{3}\setminus\partial\Omega).\] From\ $A_{5}$\ satisfies (3.54), we can get
\[ \cfrac{\partial W}{\partial  n_{p_{0}}^{+}}|_{\partial \Omega}=\lim_{M\rightarrow n_{p_{0}}^{+}}\cfrac{\partial  W(M)}{\partial n_{p_{0}}} =\cfrac{\partial h(P_{0})}{\partial n_{p_{0}}}.\ast_{\partial\Omega}A_{5}-\cfrac{1}{2}A_{5}=0.\]
From\ $W(M)\rightarrow 0$,\ if\ $M\rightarrow +\infty$, we know\ $W\equiv0$, on\ $R^{3}$. This means that\[ \cfrac{\partial W}{\partial  n_{p_{0}}^{+}}|_{\partial \Omega}-\cfrac{\partial W}{\partial  n_{p_{0}}^{-}}|_{\partial \Omega}=-A_{5}=0.\]
Hence the solution of (3.54) is only\ $0$. From Fredholm integral equation theory, we know there is only one solution for (3.53). Moreover, there exists an entire analytic function\ $\Gamma_{1}$\ which is only related to\ $\partial h(P_{0})/\partial n_{p_{0}}$, such that
\[ \cfrac{A_{5}(P_{0})}{2}=g_{02}(P_{0})+\cfrac{\partial h(P_{0})}{\partial n_{p_{0}}}.\ast \psi_{1}(Z_{1})+\int_{\partial\Omega}\Gamma_{1}(P_{0},\ P_{1})[g_{02}(P_{1})+\cfrac{\partial h(P_{1})}{\partial n_{p_{1}}}.\ast \psi_{1}(Z_{1})]dS_{P_{1}}.\]
Maybe you will point out\ $\partial h(P_{0})/\partial n_{p_{0}}(P_{0}-P_{1})$\ is not continuous, if\ $P_{0}=P_{1}$. Yes, that's true. It's lucky that $$|P_{0}-P_{1}|^{(\epsilon_{0}+5/2)}\cfrac{\partial h(P_{0})}{\partial n_{p_{0}}}(P_{0}-P_{1})$$
is continuous, if\ $\epsilon_{0}\in(0,\ 0.5]$. Hence,\ $\partial h(P_{0})/\partial n_{p_{0}}(P_{0}-P_{1})$\ is weak singular kernel. From Fredholm theorem, we know (3.59) still stands.\\From\ $ A_{5}\in C^{1}(\partial\Omega)$, we know\ $g_{02}$\ should be in\ $C^{1}(\partial\Omega)$, either.\\
At the end of this section, we introduce an easy way to get\ $A_{2},\ A_{3},\ A_{4}$.\\
For any\ $P_{0}\in \partial\Omega$, we discuss a smooth curve on\ $\partial\Omega$\ that passes through\ $P_{0}$. The parameter coordinates of the point\ $P$\ on the curve are\ $(x(\theta),\ y(\theta),\ z(\theta))$. Moreover, the parameter coordinates of \ $P_{0}$\ are\ $(x(\theta_{0}),\ y(\theta_{0}),\ z(\theta_{0}))$. We assume the tangent vector at\ $P_{0}$\ is as follows,
$$ s_{0}=(x^{\prime}(\theta_{0}),\ y^{\prime}(\theta_{0}),\ z^{\prime}(\theta_{0}))^{T}.$$
If we select a\ $P(x(\theta),\ y(\theta),\ z(\theta))$\ on the curve that is near to\ $P_{0}$, then we can get\[ u(P)-u(P_{0})=u_{x}(P_{0})(x(\theta)-x(\theta_{0}))+u_{y}(P_{0})(y(\theta)-y(\theta_{0}))+u_{z}(P_{0})(z(\theta)-z(\theta_{0}))+o(\rho),\] where\ $\rho=|P-P_{0}|=\sqrt{(x(\theta)-x(\theta_{0}))^{2}+(y(\theta)-y(\theta_{0}))^{2}+(z(\theta)-z(\theta_{0}))^{2}}$.\\
From$$u(P)|_{P\in\partial\Omega}=A_{1}((P),$$ we get,\[A_{1}(P)-A_{1}(P_{0})=A_{1x}(P_{0})(x(\theta)-x(\theta_{0}))+A_{1y}(P_{0})(y(\theta)-y(\theta_{0}))+A_{1z}(P_{0})(z(\theta)-z(\theta_{0}))+o(\rho),\]where
$$A_{1x}=\cfrac{\partial A_{1}}{\partial x},\ A_{1y}=\cfrac{\partial A_{1}}{\partial y},\ A_{1z}=\cfrac{\partial A_{1}}{\partial z}.$$
By the subtraction of (3.60) and (3.61), we can obtain \begin{eqnarray}&&(u_{x}(P_{0})-A_{1x}(P_{0}))(x(\theta)-x(\theta_{0}))+(u_{y}(P_{0})-A_{1y}(P_{0}))(y(\theta)-y(\theta_{0}))\nonumber\\&&+(u_{z}(P_{0})-A_{1z}(P_{0}))(z(\theta)-z(\theta_{0}))+o(\rho)=0.\end{eqnarray}
If we divide\ $\rho$\ on the left hand side of (3.62), and let\ $\theta\rightarrow \theta_{0}$, then from the following,
\begin{eqnarray*}\lim_{\theta\rightarrow \theta_{0}}\cfrac{x(\theta)-x(\theta_{0})}{\rho}&=&\lim_{\theta\rightarrow \theta_{0}}\cfrac{(x(\theta)-x(\theta_{0}))/(\theta-\theta_{0})}{\rho/(\theta-\theta_{0})}=\cfrac{x^{\prime}(\theta_{0})}{\parallel s_{0}\parallel},\\
\lim_{\theta\rightarrow \theta_{0}}\cfrac{y(\theta)-y(\theta_{0})}{\rho}&=&\lim_{\theta\rightarrow \theta_{0}}\cfrac{(y(\theta)-y(\theta_{0}))/(\theta-\theta_{0})}{\rho/(\theta-\theta_{0})}=\cfrac{y^{\prime}(\theta_{0})}{\parallel s_{0}\parallel},\\
\lim_{\theta\rightarrow \theta_{0}}\cfrac{z(\theta)-z(\theta_{0})}{\rho}&=&\lim_{\theta\rightarrow \theta_{0}}\cfrac{(z(\theta)-z(\theta_{0}))/(\theta-\theta_{0})}{\rho/(\theta-\theta_{0})}=\cfrac{z^{\prime}(\theta_{0})}{\parallel s_{0}\parallel},\end{eqnarray*}
where\ $\parallel s_{0}\parallel=\sqrt{(x^{\prime}(\theta_{0}))^{2}+(y^{\prime}(\theta_{0}))^{2}+(z^{\prime}(\theta_{0}))^{2}}$, we get,\[(u_{x}(P_{0})-A_{1x}(P_{0}))x^{\prime}(\theta_{0})+(u_{y}(P_{0})-A_{1y}(P_{0}))y^{\prime}(\theta_{0})+(u_{z}(P_{0})-A_{1z}(P_{0}))z^{\prime}(\theta_{0})=0.\]
(3.63) will still stand even if\ $\parallel s_{0}\parallel=0$. From the arbitrary nature of\ $s_{0}$, we know that$$((u_{x}(P_{0})-A_{1x}(P_{0}))_{k},\ (u_{y}(P_{0})-A_{1y}(P_{0}))_{k},\ (u_{z}(P_{0})-A_{1z}(P_{0}))_{k})^{T},\ 1\leq k\leq m,$$should be parallel to the exterior normal vector\ $n_{p_{0}}=(n_{1}(P_{0}),\ n_{2}(P_{0}),\ n_{3}(P_{0}))^{T}$, where\ $(\ast)_{k}$\ is the\ $k$th component of the vector,\ $1\leq k\leq m$. \\ Hence there exists\ $\lambda(P_{0})$, such that\begin{eqnarray*}u_{x}(P_{0})-A_{1x}(P_{0})&=&\lambda(P_{0})n_{1}(P_{0}),\\u_{y}(P_{0})-A_{1y}(P_{0})&=&\lambda(P_{0})n_{2}(P_{0}),\\u_{z}(P_{0})-A_{1z}(P_{0})&=&\lambda(P_{0})n_{3}(P_{0}).\end{eqnarray*}
From\ $n_{1}^{2}(P_{0})+n_{2}^{2}(P_{0})+n_{3}^{2}(P_{0})=1$, we can get,
\begin{eqnarray*}\lambda(P_{0})&=&(u_{x}(P_{0})-A_{1x}(P_{0}))n_{1}(P_{0})+(u_{y}(P_{0})-A_{1y}(P_{0}))n_{2}(P_{0})+\\&&(u_{z}(P_{0})-A_{1z}(P_{0}))n_{3}(P_{0})\\&=&A_{5}(P_{0})-\cfrac{\partial A_{1}}{\partial n_{p_{0}}}(P_{0}).\end{eqnarray*}
Hence, we can work out\ $A_{2},\ A_{3},\ A_{4}$\ as follows, \begin{eqnarray}A_{2}&=&u_{x}|_{\partial\Omega}=\cfrac{\partial A_{1}}{\partial x}+(A_{5}-\cfrac{\partial A_{1}}{\partial n})n_{1},\\A_{3}&=&u_{y}|_{\partial\Omega}=\cfrac{\partial A_{1}}{\partial y}+(A_{5}-\cfrac{\partial A_{1}}{\partial n})n_{2},\\A_{4}&=&u_{z}|_{\partial\Omega}=\cfrac{\partial A_{1}}{\partial z}+(A_{5}-\cfrac{\partial A_{1}}{\partial n})n_{3},\end{eqnarray}
where\ $n=(n_{1},\ n_{2},\ n_{3})^{T}$\ is the exterior normal vector to\ $\partial\Omega$.\\
It looks that the classical solution of Eqs(3.41) would be locally exist and unique. So were the Eqs(1.1). But\ $F^{-1}[g_{1}]$\ will cause some troubles.\\
          We can work out$$F^{-1}[g_{1}]=\int_{\partial\Omega}\delta(x-x_{1},\ y-y_{1},\ z-z_{1})A_{1}(x_{1},\ y_{1},\ z_{1})n_{1}(x_{1},\ y_{1},\ z_{1})dS,$$where\ $\delta$\ is Dirac function. We don't know what the end will be.\\ Because\ $\partial\Omega\in C^{1,\ \beta}$, we can get the following from Theorem 2.1, $$\partial\Omega=\bigcup_{k=1}^{N}\partial\Omega_{k},$$ where\ $\partial\Omega_{k}$\ is the graph of a\ $C^{1,\ \beta}$\ function of\ $2$\ of the coordinates\ $x,\ y,\ z$.\\ Without loss of the generality, we assume\ $\partial\Omega_{k}$\ is the graph of a\ $C^{1,\ \beta}$\ function\ $z=f_{k}(x,\ y),\ (x,\ y)\in D_{k}\subset R^{2}$, then we may solve inverse Fourier transform of\ $g_{1}$\ on\ $\partial\Omega_{k}$\ as follows,\begin{eqnarray*}&&F^{-1}[\int_{\partial\Omega_{k}}A_{1}(x_{1},\ y_{1},\ z_{1})n_{1}(x_{1},\ y_{1},\ z_{1})e^{-i\xi_{1}x_{1}-i\xi_{2}y_{1}-i\xi_{3}z_{1}}dS]\\
          &=&F^{-1}[\int_{D_{k}}A_{1}(x_{1},\ y_{1},\ f_{k}(x_{1},\ y_{1}))n_{1}(x_{1},\ y_{1},\ f_{k}(x_{1},\ y_{1}))e^{-i\xi_{1}x_{1}-i\xi_{2}y_{1}-i\xi_{3}f_{k}(x_{1},\ y_{1})}\\&&\sqrt{1+f^{2}_{k1}(x_{1},\ y_{1})+f^{2}_{k2}(x_{1},\ y_{1})}dx_{1}dy_{1}]\\&=&A_{1}(x,\ y,\ f_{k}(x,\ y))n_{1}(x,\ y,\ f_{k}(x,\ y))\sqrt{1+f^{2}_{k1}(x,\ y)+f^{2}_{k2}(x,\ y)}I_{D_{k}}\\&&\cfrac{1}{2\pi}\int_{-\infty}^{+\infty}e^{-i\xi_{3}f_{k}(x,\ y)}e^{i\xi_{3}z}d\xi_{3}\\&=&A_{1}(x,\ y,\ f_{k}(x,\ y))n_{1}(x,\ y,\ f_{k}(x,\ y))\sqrt{1+f^{2}_{k1}(x,\ y)+f^{2}_{k2}(x,\ y)}I_{D_{k}}\delta[z-f_{k}(x,\ y)],\end{eqnarray*}
          where$$f_{k1}(x,\ y)=\cfrac{\partial f_{k}(x,\ y)}{\partial x},\ f_{k2}(x,\ y)=\cfrac{\partial f_{k}(x,\ y)}{\partial y}.$$ We are surprised to see that inverse Fourier transform of\ $g_{1}$\ on\ $\partial\Omega_{k}$\ is related with\ $\delta[z-f_{k}(x,\ y)]$. This means that\ $F^{-1}[g_{1}]$\ will be related with\ $\delta(\partial\Omega)$. And it is possible that\ $\partial F^{-1}(a_{1}^{-1})$\ or\ $\partial^{2} F^{-1}(a_{1}^{-1})$\ is not integrable in\ $w_{1},\ w_{2}$.\\
       We regret that we lost local existence and uniqueness for the solution of\ $Z_{1}=T_{0}(Z_{1})$. However, it would be too easy if there were no\ $\delta(\partial\Omega)$. In the next section, we will discuss\ $Z_{1}=T_{0}(Z_{1})$\ by Lerry-Schauder degree and the Sobolev space\ $H^{-m_{1}}(\Omega)$.\\
\section{Existence} \setcounter{equation}{0}
In this section, we will discuss the existence for classical solution of Eqs(4.1). We include remarks if it is necessary. First of all, we introduce our ideas as follows.\\
First, we will construct a norm\ $\|\cdot\|_{-m_{1}}$\ for\ $T_{0}(Z_{1})$.\\ Second, since\ $Z_{1}=T_{0}(Z_{1})$\ are generalized integral equations, we make approximate ordinary integral equations\ $Z_{1}=T_{0\epsilon}(Z_{1}),\ \forall \epsilon>0$,
such that\ $\forall Z_{1}\in C(\overline{\Omega}),\ \|Z_{1}\|_{\infty}\leq M,\ M$\ is given, we have the following,
\[ \lim_{\epsilon\rightarrow 0}\|T_{0\epsilon}(Z_{1})-T_{0}(Z_{1})\|_{-m_{1}}=0,\ \mbox{uniformly,}\]where\[\|Z_{1}\|_{\infty}=\max_{1\leq i\leq 9m}\|Z_{1,\ i}\|_{\infty},\ \|Z_{1,\ i}\|_{\infty}=\max_{X=(x,\ y,\ z)^{T}\in\overline{\Omega}}|Z_{1,\ i}(X)|,\] $Z_{1,\ i},\ 1\leq i\leq 9m$, are components of\ $Z_{1}$.\\ This will help us to prove\ $T_{0}(Z_{1k})$\ is sequentially compact under norm\ $\|\cdot\|_{-m_{1}}$, if\ $\forall Z_{1k}\in C(\overline{\Omega}),\ \|Z_{1k}\|_{\infty}\leq M,\ k\geq 1,\ M$\ is given.\\Finally, we use some primary theorems on Leray-Schauder degree to discuss ordinary integral equations\ $Z_{1}=T_{0\epsilon}(Z_{1}),\ \forall \epsilon>0$. If\ $Z_{1\epsilon}$\ satisfies\ $Z_{1}=T_{0\epsilon}(Z_{1}),\ \forall \epsilon>0$,\ and bounded uniformly, then there exists a sequence\ $Z_{1\epsilon_{k}},\ k\geq1,\ \epsilon_{k}\rightarrow 0$, if\ $k\rightarrow+\infty$, such that\ $Z_{1\epsilon_{k}}$\ is convergent to\ $Z_{1}^{\ast}$. We will get\ $Z_{1}^{\ast}=T_{0}(Z_{1}^{\ast})$, which is what we want.\\Can our imagination come true? Let's introduce our answer. The answer is not unique.
\begin{definition}\label{definition}$\forall Z_{1}=(Z_{1,\ i})_{9m\times 1},\ Z_{1,\ i}\in H^{-m_{1}}(\Omega),\ 1\leq i\leq 9m$, we have
\[\|Z_{1}\|_{-m_{1}}=\max_{1\leq i\leq 9m}\|Z_{1,\ i}\|_{-m_{1}},\ \|Z_{1,\ i}\|_{-m_{1}}=\sup_{\varphi\in C_{0}^{\infty}(\Omega)}\cfrac{|<Z_{1,\ i},\ \varphi>|}{\|\varphi\|_{m_{1}}},\] where\ $m_{1}=6+2a,\ a=\max\{\partial(a_{1}B^{-1}_{1}),\ \partial(a_{1}B^{-1}_{1}B_{2})\}$, here\ $\partial(\cdot)$\ means the highest degree,\ $X=(x,\ y,\ z)^{T}$,\ $<Z_{1,\ i},\ \varphi>$\ is the value of the generalized function\ $Z_{1,\ i}$\ on\ $\varphi$, if\ $Z_{1,\ i}$\ is locally integrable, then \begin{eqnarray}<Z_{1,\ i},\ \varphi>&=&\int_{\Omega}Z_{1,\ i}(X)\varphi(X) dX,\ 1\leq i\leq 9m,\\ \|\varphi\|_{m_{1}}&=&(\int_{\Omega}\sum_{|\alpha|\leq m_{1}}|\partial^{\alpha}\varphi(X)|^{2}dX)^{1/2}.\end{eqnarray}\end{definition}
And we can get that\[\|T_{0}(Z_{1})\|_{-m_{1}}=\max_{1\leq i\leq 9m}\sup_{\varphi\in C_{0}^{\infty}(\Omega)}\cfrac{|<T_{0,\ i}(Z_{1}),\ \varphi>|}{\|\varphi\|_{m_{1}}}<+\infty,\]in the following lemma, where\ $ T_{0,\ i}(Z_{1}),\ 1\leq i\leq 9m$, are components of\ $T_{0}(Z_{1})$.\\Next, we see approximate ordinary integral equations,\ $\forall \epsilon>0$,
\[ T_{0\epsilon}(Z_{1})=w_{1\epsilon}+w_{2\epsilon}.\ast(\psi(Z_{1}I_{\overline{\Omega}})),\] where$$w_{1\epsilon}=F^{-1}(\widetilde{\delta_{\epsilon}}B_{1}^{-1}\beta_{1}),\ w_{2\epsilon}=F^{-1}(\widetilde{\delta_{\epsilon}}B_{1}^{-1}B_{2}),\ \widetilde{\delta_{\epsilon}}=F(\delta_{\epsilon}),$$
        \[\delta_{\epsilon}=\cfrac{1}{(\sqrt{\pi\epsilon})^{3}}e^{-|X|^{2}/\epsilon},\ |X|=\sqrt{x^{2}+y^{2}+z^{2}},\ \lim_{\epsilon\rightarrow 0}\delta_{\epsilon}=\delta(X),\ \lim_{\epsilon\rightarrow 0}\widetilde{\delta_{\epsilon}}=1,  \] $\delta(X)$\ is the Dirac function.\\In Dirichlet problem,\ $A_{1}$\ is known, and from (3.59), we can get\ $A_{5}$\ as follows,
         \[ \cfrac{A_{5}(P_{0})}{2}=g_{02}(P_{0})+\cfrac{\partial h(P_{0})}{\partial n_{p_{0}}}.\ast \psi_{1}(Z_{1})+\int_{\partial\Omega}\Gamma_{1}(P_{0},\ P_{1})[g_{02}(P_{1})+\cfrac{\partial h(P_{1})}{\partial n_{p_{1}}}.\ast \psi_{1}(Z_{1})]dS_{P_{1}}.\]
         We can see that\ $A_{5}$\ is a compact operator on\ $Z_{1}$.\\ From (3.64), (3.65), (3.66), we can also get\ $A_{2},\ A_{3},\ A_{4}$\ as the following,
         \begin{eqnarray}A_{2}&=&u_{x}|_{\partial\Omega}=\cfrac{\partial A_{1}}{\partial x}+(A_{5}-\cfrac{\partial A_{1}}{\partial n})n_{1},\\A_{3}&=&u_{y}|_{\partial\Omega}=\cfrac{\partial A_{1}}{\partial y}+(A_{5}-\cfrac{\partial A_{1}}{\partial n})n_{2},\\A_{4}&=&u_{z}|_{\partial\Omega}=\cfrac{\partial A_{1}}{\partial z}+(A_{5}-\cfrac{\partial A_{1}}{\partial n})n_{3},\end{eqnarray}
where\ $n=(n_{1},\ n_{2},\ n_{3})^{T}$\ is the exterior normal vector to\ $\partial\Omega$.\\We can see that\ $A_{2},\ A_{3},\ A_{4}$\ are all compact operators on\ $Z_{1}$.\\
We know\ $v_{1}\ast v_{2}\in S^{\prime}$\ and\ $F(v_{1}\ast v_{2})=F(v_{1})F(v_{2})$\ will still hold if\ $v_{1}\in S,\ v_{2}\in S^{\prime}$.  There is an example of this on pages 118 to 119 of [10]. \\We can see that\ $\forall \epsilon>0$,\ $\delta_{\epsilon}\in S,\ F^{-1}(a_{1}^{-1})\in S^{\prime},\ F^{-1}(g_{j})\in S^{\prime},\ F^{-1}(g_{jk})\in S^{\prime},\ 1\leq j,\ k\leq 3$. Hence we can get the following,\begin{eqnarray*} F^{-1}[\widetilde{\delta_{\epsilon}}(i\xi)^{\alpha}a_{1}^{-1}]&=&(\partial^{\alpha}\delta_{\epsilon}).\ast F^{-1}(a_{1}^{-1}),\\ F^{-1}[\widetilde{\delta_{\epsilon}}(i\xi)^{\alpha}g_{j}]&=&(\partial^{\alpha}\delta_{\epsilon}).\ast_{\partial\Omega}(A_{3}n_{j}),\\ F^{-1}(\widetilde{\delta_{\epsilon}}g_{jk})&=&\delta_{\epsilon}.\ast_{\partial\Omega}(A_{j+3}n_{k}),\end{eqnarray*}where\ $1\leq j,\ k\leq 3.$\\
We may denote that\ $T_{0\epsilon}(Z_{1})=\delta_{\epsilon}.\ast T_{0}(Z_{1})$. From (3.16) and\ $\partial^{\alpha}\delta_{\epsilon}\in S$, we can get$$ (\partial^{\alpha}\delta_{\epsilon}).\ast F^{-1}(a_{1}^{-1})\in C(R^{3}),\ (\partial^{\alpha}\delta_{\epsilon}.\ast F^{-1}(a_{1}^{-1})).\ast(\psi(Z_{1}I_{\overline{\Omega}}))\in C(\overline{\Omega}).$$ We see that\ $\forall \epsilon>0$, there is no\ $\delta(\partial\Omega)$\ or\ $\partial F^{-1}(a_{1}^{-1})$\ in\ $T_{0\epsilon}$\ again. So\ $Z_{1}I_{\overline{\Omega}}=T_{0\epsilon}(Z_{1}I_{\overline{\Omega}})$\ are ordinary integral equations. Moreover\ $T_{0\epsilon}(Z_{1})$\ is bounded uniformly and equicontinuous if\ $Z_{1}\in C(\overline{\Omega})$\ and bounded uniformly. From the Arzela-Ascoli theorem, we know\ $T_{0\epsilon}$\ is a compact operator on\ $Z_{1}$. And the Leray-Schauder degree can work now. \\For the preliminaries, we have a lemma as follows.\begin{lemma} \label{lemma1}(1)There exists\ $C>0$, such that\[ \|Z_{1}\|_{-m_{1}}\leq C \|Z_{1}\|_{\infty},\ \forall Z_{1}\in C(\overline{\Omega}).\]
        (2) $\forall Z_{1}\in C(\overline{\Omega}),\ \|Z_{1}\|_{\infty}\leq M,\ M>0,\ M$\ is given, we have
\[ \|T_{0}(Z_{1})\|_{-m_{1}}<+\infty,\ \mbox{and}\ \lim_{\epsilon\rightarrow 0}\|T_{0\epsilon}(Z_{1})-T_{0}(Z_{1})\|_{-m_{1}}=0,\ \mbox{uniformly.}\]
(3) $T_{0}(Z_{1k}),\ k\geq1,$\ is sequentially compact under norm\ $\|\cdot\|_{-m_{1}}$, if\ $Z_{1k}\in C(\overline{\Omega}),\ \|Z_{1k}\|_{\infty}\leq M,\ k\geq 1,\ M>0,\ M$\ is given.\\(4)$\forall \epsilon_{0}>0,\ \forall Z_{1}\in C(\overline{\Omega}),\ \|Z_{1}\|_{\infty}\leq M,\ M>0,\ M$\ is given, we have\[ \lim_{\epsilon\rightarrow\epsilon_{0}}\|[T_{0\epsilon}(Z_{1})-T_{0\epsilon_{0}}(Z_{1})]I_{\overline{\Omega}}\|_{\infty}=0.\]\end{lemma}
{\it Proof of lemma 4.1}. (1)From\ $\|Z_{1,\ i}\|_{-m_{1}}\leq \|Z_{1,\ i}\|_{L^{2}}\leq \|Z_{1,\ i}\|_{\infty}\sqrt{m(\overline{\Omega})},\ 1\leq i\leq 9m$, where\ $m(\overline{\Omega})$\ is the Lebesgue measure of\ $\overline{\Omega}$. Hence we get\ $\|Z_{1}\|_{-m_{1}}\leq \|Z_{1}\|_{\infty}\sqrt{m(\overline{\Omega})}$. We may let\ $C=\sqrt{m(\overline{\Omega})}$.\\
There will not exist\ $C>0$, such that\[ \|Z_{1}\|_{\infty}\leq C \|Z_{1}\|_{-m_{1}},\ \forall Z_{1}\in C(\overline{\Omega}).\]
We may select\ $ Z_{1k}\in C(\overline{\Omega}),\ k\geq1,\ \|Z_{1k}\|_{\infty}\equiv1$, but\ $Z_{1k}\rightarrow 0,\ \mbox{a.e.}$, here a.e. means that almost everywhere. Then from the Lebesgue dominated convergence theorem,\ $\|Z_{1k}\|_{-m_{1}}\rightarrow 0$. Hence (4.16) will not stand.\ $\|Z_{1}\|_{-m_{1}}$\ and\ $\|Z_{1}\|_{\infty}$\ are not equivalent.\\
(2)We will prove\ $(F[T_{0}(Z_{1})])F^{-1}(\varphi)\in L^{1}(R^{3})$. We can work out the following,
\begin{eqnarray*}<T_{0}(Z_{1}),\ \varphi>&=&<F[T_{0}(Z_{1})],\ F^{-1}(\varphi)>\\ &=&<(B^{-1}_{1}\beta_{1}+B^{-1}_{1}B_{2}FI(Z_{2})),\ F^{-1}(\varphi)>\\&=&<(B^{-1}_{1}\beta_{1}+B^{-1}_{1}B_{2}FI(Z_{2}))(1+|\xi|^{2})^{-a-3},\\&& (1+|\xi|^{2})^{a+3}F^{-1}(\varphi)>,
\end{eqnarray*}where\ $\xi=(\xi_{1},\ \xi_{2},\ \xi_{3})^{T}\in R^{3},\ a_{1}B^{-1}_{1},\ a_{1}B^{-1}_{1}B_{2}$\ are polynomial matrices,\ $Z_{2}=\psi(Z_{1})$.\\ From\ $\psi(Z_{1})$\ is a continuous, we get that\ $\beta_{1},\ Z_{2}$\ are all bounded. \\ We have the following,\begin{eqnarray*}|(i\xi)^{\alpha}F^{-1}(\varphi)|/\|\varphi\|_{m_{1}}&\leq&\int_{\Omega}|\partial^{\alpha}\varphi|dX/\|\varphi\|_{m_{1}}\\
&\leq&(\int_{\Omega}|\partial^{\alpha}\varphi|^{2}dX)^{1/2}(\int_{\Omega}dX)^{1/2}/\|\varphi\|_{m_{1}}\\&\leq&\sqrt{m(\overline{\Omega})},\end{eqnarray*}
where\ $|\alpha|\leq 6+2a,\ a=\max\{\partial(a_{1}B^{-1}_{1}),\ \partial(a_{1}B^{-1}_{1}B_{2})\},\ m_{1}=6+2a$.\\
From (3.17), we can get\ $(F[T_{0}(Z_{1})])F^{-1}(\varphi)\in L^{1}(R^{3})$, moreover\[\|T_{0}(Z_{1})\|_{-m_{1}}\leq C_{3,\ 0},\]
where\ $C_{3,\ 0}>0$,\ is a constant only related to\ $M$. It is not related with\ $Z_{1}$.\\
We can work out as follows,
\begin{eqnarray*}<T_{0\epsilon}(Z_{1})-T_{0}(Z_{1}),\ \varphi>&=&<F[T_{0\epsilon}(Z_{1})-T_{0}(Z_{1})],\ F^{-1}(\varphi)>\\ &=&<(\widetilde{\delta_{\epsilon}}-1)F[T_{0}(Z_{1})],\ F^{-1}(\varphi)>,\\ &=&<(\widetilde{\delta_{\epsilon}}-1)(B^{-1}_{1}\beta_{1}+B^{-1}_{1}B_{2}FI(Z_{2})),\ F^{-1}(\varphi)>\\&=&<(\widetilde{\delta_{\epsilon}}-1)(B^{-1}_{1}\beta_{1}+B^{-1}_{1}B_{2}FI(Z_{2}))(1+|\xi|^{2})^{-a-3},\\&& (1+|\xi|^{2})^{a+3}F^{-1}(\varphi)>.
\end{eqnarray*}
Hence, we can get\[\|T_{0\epsilon}(Z_{1})-T_{0}(Z_{1})\|_{-m_{1}}\leq
\int_{R^{3}}|\widetilde{\delta_{\epsilon}}-1||(B^{-1}_{1}\beta_{1}+B^{-1}_{1}B_{2}FI(Z_{2}))(1+|\xi|^{2})^{-a-3}|\sqrt{m(\overline{\Omega})}d\xi.\]
We can also work out the following,
\begin{eqnarray*}\widetilde{\delta_{\epsilon}}&=&\int_{R^{3}}\cfrac{1}{(\sqrt{\pi\epsilon})^{3}}e^{-|X|^{2}/\epsilon}e^{-i\xi\cdot X}dX,\ (X=\sqrt{\epsilon}Y,\ dX=(\sqrt{\epsilon})^{3}dY)
\\&=&\int_{R^{3}}\cfrac{1}{(\sqrt{\pi})^{3}}e^{-|Y|^{2}}e^{-i\xi\cdot \sqrt{\epsilon}Y}dY,
\\ &=&\int_{R^{3}}\cfrac{1}{(\sqrt{\pi})^{3}}e^{-|Y+i\sqrt{\epsilon}\xi/2 |^{2}}e^{-\epsilon|\xi|^{2}/4}dY,\\ &=&\int_{R^{3}}\cfrac{1}{(\sqrt{\pi})^{3}}e^{-|Y |^{2}}e^{-\epsilon|\xi|^{2}/4}dY=e^{-\epsilon|\xi|^{2}/4}.\end{eqnarray*}Here we use Cauchy contour integral as follows,
\[\int_{-\infty}^{+\infty} e^{-(x+i\xi_{1})^{2}}dx=\int_{-\infty}^{+\infty} e^{-x^{2}}dx,\ \forall \xi_{1}\in R.\]
We can get the following,
$$ \int_{R^{3}}|\widetilde{\delta_{\epsilon}}-1||(B^{-1}_{1}\beta_{1}+B^{-1}_{1}B_{2}FI(Z_{2}))(1+|\xi|^{2})^{-a-3}|\sqrt{m(\overline{\Omega})}d\xi=I_{3,\ 1}+I_{3,\ 2},$$ where\begin{eqnarray} I_{3,\ 1}&=&\int_{|\xi|>M_{0}}|\widetilde{\delta_{\epsilon}}-1||(B^{-1}_{1}\beta_{1}+B^{-1}_{1}B_{2}FI(Z_{2}))(1+|\xi|^{2})^{-a-3}|\sqrt{m(\overline{\Omega})}d\xi,\\ I_{3,\ 2}&=&\int_{|\xi|\leq M_{0}}|\widetilde{\delta_{\epsilon}}-1||(B^{-1}_{1}\beta_{1}+B^{-1}_{1}B_{2}FI(Z_{2}))(1+|\xi|^{2})^{-a-3}|\sqrt{m(\overline{\Omega})}d\xi.\end{eqnarray}
$\forall \epsilon^{\prime}>0$, there exists\ $M_{0}>0$, which is only related to\ $M$,
such that\ $|I_{3,\ 1}|\leq \epsilon^{\prime}/2$. And for such\ $M_{0}$, there exists\ $\delta_{0}>0$, such that\ $|I_{3,\ 2}|\leq \epsilon^{\prime}/2$, if\ $\epsilon\leq \delta_{0}$.\\
Hence, we obtain\ $\|T_{0\epsilon}(Z_{1})-T_{0}(Z_{1})\|_{-m_{1}}\rightarrow0,$\ uniformly.\\
(3)We will prove\ $T_{0}(Z_{1k}),\ k\geq 1,$\ is totally bounded. From (2),\ $\forall \epsilon_{1}>0,\ \exists \delta_{1}>0$, such that \[\|T_{0\epsilon}(Z_{1k})-T_{0}(Z_{1k})\|_{-m_{1}}\leq \epsilon_{1}/3,\ \forall \epsilon\in (0,\ \delta_{1}),\ \forall k,\ k\geq 1.\]  If we choose\ $\epsilon_{0}\in (0,\ \delta_{1})$, then\ $T_{0\epsilon_{0}}(Z_{1k}),\ k\geq 1$\ is sequentially compact. There exist finite\ $\epsilon_{1}/(3C)$\ net\ $T_{0\epsilon_{0}}(Z_{1k_{1}}),\ T_{0\epsilon_{0}}(Z_{1k_{2}}),\cdots,\ T_{0\epsilon_{0}}(Z_{1k_{s}})$, where\ $C$\ is defined in (1).\\ This means that\ $\forall k,\ \exists l,\ 1\leq l\leq s$, such that\ $\|T_{0\epsilon_{0}}(Z_{1k})-T_{0\epsilon_{0}}(Z_{1k_{l}})\|_{\infty}\leq \epsilon_{1}/(3C)$.\\ From\ $\|T_{0\epsilon_{0}}(Z_{1k})-T_{0}(Z_{1k})\|_{-m_{1}}\leq \epsilon_{1}/3,\ \|T_{0\epsilon_{0}}(Z_{1k})-T_{0\epsilon_{0}}(Z_{1k_{l}})\|_{-m_{1}}\leq \epsilon_{1}/3$,\\ $\|T_{0\epsilon_{0}}(Z_{1k_{l}})-T_{0}(Z_{1k_{l}})\|_{-m_{1}}\leq \epsilon_{1}/3$, we obtain\ $\|T_{0}(Z_{1k})-T_{0}(Z_{1k_{l}})\|_{-m_{1}}\leq \epsilon_{1}$. \\So\ $T_{0}(Z_{1k_{1}}),\ T_{0}(Z_{1k_{2}}),\cdots,\ T_{0}(Z_{1k_{s}})$\ is a finite\ $\epsilon_{1}$\ net for\ $T_{0}(Z_{1k}),\ k\geq 1$. \\This means that\ $T_{0}(Z_{1k}),\ k\geq 1,$\ is totally bounded. From Hausdorff theorem on page 14 in [7],\ $T_{0}(Z_{1k}),\ k\geq1,$\ is sequentially compact under norm\ $\|\cdot\|_{-m_{1}}$.\\
(4)From\ $T_{0\epsilon,\ i}(Z_{1})=\delta_{\epsilon}.\ast T_{0,\ i}(Z_{1}),\ 1\leq i\leq 9m$, we obtain
\begin{eqnarray*}T_{0\epsilon,\ i}(Z_{1})-T_{0\epsilon_{0},\ i}(Z_{1})&=&(\delta_{\epsilon}-\delta_{\epsilon_{0}}).\ast T_{0,\ i}(Z_{1})
\\&=&F^{-1}(F((\delta_{\epsilon}-\delta_{\epsilon_{0}}).\ast T_{0,\ i}(Z_{1})))\\&=&\cfrac{1}{(2\pi)^{3}}\int_{R^{3}}(e^{-\epsilon|\xi|^{2}/4}-e^{-\epsilon_{0}|\xi|^{2}/4})F(T_{0,\ i}(Z_{1}))e^{i\xi\cdot X}d\xi
\\&=& I_{3,\ 3}+I_{3,\ 4},\end{eqnarray*}
where\begin{eqnarray} I_{3,\ 3}=\cfrac{1}{(2\pi)^{3}}\int_{|\xi|>M_{0}}(e^{-\epsilon|\xi|^{2}/4}-e^{-\epsilon_{0}|\xi|^{2}/4})(1+|\xi|^{2})^{a+3}(F(T_{0,\ i}(Z_{1}))(1+|\xi|^{2})^{-a-3})e^{i\xi\cdot X}d\xi,&&\\ I_{3,\ 4}=\cfrac{1}{(2\pi)^{3}}\int_{|\xi|\leq M_{0}}(e^{-\epsilon|\xi|^{2}/4}-e^{-\epsilon_{0}|\xi|^{2}/4})(1+|\xi|^{2})^{a+3}(F(T_{0,\ i}(Z_{1}))(1+|\xi|^{2})^{-a-3})e^{i\xi\cdot X}d\xi.&&\end{eqnarray}From (3.17), we know\ $F(T_{0,\ i}(Z_{1}))(1+|\xi|^{2})^{-a-3}\in L^{1}(R^{3}),\ 1\leq i\leq 9m$.\\ If we let\ $\epsilon\in [\epsilon_{0}/2,\ 3\epsilon_{0}/2]$, then\ $\forall \epsilon^{\prime}>0$, there exists\ $M_{0}>0$, which is related to\ $\epsilon_{0}$,
such that\ $|I_{3,\ 3}|\leq \epsilon^{\prime}/2$. And for such\ $M_{0}$, there exists\ $\delta_{0}\in (0,\ \epsilon_{0}/2)$, such that\ $|I_{3,\ 4}|\leq \epsilon^{\prime}/2$, if\ $|\epsilon-\epsilon_{0}|\leq \delta_{0}$. So the statement holds.
\qed\\
We will use result (4) into Leray-Schauder degree. This is the reason why we choose\ $\delta_{\epsilon}$\ instead of\ $I_{\{|X|\leq \epsilon\}}/ |I_{\{|X|\leq \epsilon\}}|$, where \[|I_{\{|X|\leq \epsilon\}}|=\int_{|X|\leq \epsilon}dX.\] The latter will not satisfy (4). Maybe it works after being polished. We haven't tested it yet.\\ Now we introduce Leray-Schauder degree.
\begin{definition}\label{definition}$\Omega_{0}$\ is bounded open set of real Banach space\ $B,\ T:\ \overline{\Omega_{0}}\rightarrow B$\ is totally continuous,\ $f(x)=x-T(x),\ \forall x\in \overline{\Omega_{0}},\ p\in B\setminus f(\partial\Omega_{0})$, \[ \tau=\inf_{x\in\partial\Omega_{0}}\|f(x)-p\|>0.\] There exists\ $B^{(n)}$\ is subspace of\ $B$\ with finite dimensions,\ $p\in B^{(n)}$, and there exists bounded continuous operator\ $T_{n}:\ \overline{\Omega_{0}}\rightarrow B^{(n)}$, such that\[ \|T(x)-T_{n}(x)\|<\tau,\ \forall x\in \overline{\Omega_{0}}.\] Then Leray-Schauder degree of totally continuous field\ $f$\ is
\[ deg(f,\ \Omega_{0},\ p)=deg(f_{n},\ \Omega_{0,\ n},\ p),\]where\ $f_{n}=I-T_{n},\ \Omega_{0,\ n}=B^{(n)}\cap\Omega_{0}$.\end{definition}
Maybe you are not very familiar with Leray-Schauder degree or even you know nothing about it. That's not the problem. We only apply three primary theorems to\ $Z_{1}=T_{0\epsilon}(Z_{1}),\ \forall \epsilon>0$. We write them together into a lemma as follows.
\begin{lemma} \label{lemma1}(1)(Kronecker) If\ $deg(f,\ \Omega_{0},\ p)\neq 0$, then there exists solution for\ $f(x)=p$\ in\ $\Omega_{0}$.\\
(2)(Rothe) If\ $\Omega_{0}$\ is bounded and open convex set in Banach space\ $B$,\ $T: \overline{\Omega_{0}}\rightarrow B$\ is totally continuous,\ $T(\partial\Omega_{0})\subset \overline{\Omega_{0}},$\ and\ $T(x)\neq x,\ \forall x\in \partial\Omega_{0}$, then\ $deg(I-T,\ \Omega_{0},\ 0)\neq 0$.\\
(3) If\ $f=I-T$\ and\ $f_{1}=I-T_{1}$\ are all totally continuous fields mapping from\ $\overline{\Omega_{0}}$\ to Banach space\ $B$,\ $p$\ is not in\ $f(\partial\Omega_{0})\cup f_{1}(\partial\Omega_{0})$, moreover\[ \|T_{1}(x)-T(x)\|\leq \|x-T(x)-p\|,\ \forall x\in \partial\Omega_{0},\]
then\ $deg(f,\ \Omega_{0},\ p)=deg(f_{1},\ \Omega_{0},\ p)$.
\end{lemma}
(3) is called homotopic. You may read the explanation of definition and all the proofs of three primary theorems from page 135 to page 165 in [8]. We won't repeat them again.\\
We denote as follows, \[\tau(M,\ \epsilon,\ \Omega)=\inf_{\|Z_{1}\|_{\infty}=M}\|Z_{1}-T_{0\epsilon}(Z_{1})\|_{\infty}.\]
If\ $\Omega$\ is fixed, then we denote\ $\tau(M,\ \epsilon,\ \Omega)$\ into\ $\tau(M,\ \epsilon)$.\\If\ $\tau(M,\ \epsilon,\ \Omega)=0$, then there exists a sequence\ $Z_{1k},\ \|Z_{1k}\|_{\infty}=M,\ k\geq1$, such that
\[\lim_{k\rightarrow+\infty}\|Z_{1k}-T_{0\epsilon}(Z_{1k})\|_{\infty}=0.\]
Because\ $T_{0\epsilon}$\ is compact, there exist a sub-sequence\ $Z_{1n_{k}},\ k\geq1,\ \mbox{and}\ Z_{1\epsilon}\in C(\overline{\Omega})$, such that\[\lim_{k\rightarrow+\infty}\|Z_{1\epsilon}-T_{0\epsilon}(Z_{1n_{k}})\|_{\infty}=0.\] From (4.31), we obtain\[\lim_{k\rightarrow+\infty}\|Z_{1n_{k}}-T_{0\epsilon}(Z_{1n_{k}})\|_{\infty}=0.\]
And we can get that\[\lim_{k\rightarrow+\infty}\|Z_{1\epsilon}-Z_{1n_{k}}\|_{\infty}=0.\] Because\ $T_{0\epsilon}$\ is continuous, we can obtain
\[\lim_{k\rightarrow+\infty}\|T_{0\epsilon}(Z_{1\epsilon})-T_{0\epsilon}(Z_{1n_{k}})\|_{\infty}=0.\] (4.32) and (4.35) mean that\ $Z_{1\epsilon}=T_{0\epsilon}(Z_{1\epsilon})$, where\ $Z_{1\epsilon}\in C(\overline{\Omega}),\ \|Z_{1\epsilon}\|_{\infty}=M.$\\
Now we see the solution of\ $Z_{1\epsilon}=T_{0\epsilon}(Z_{1\epsilon})$\ as follows.
\begin{theorem} \label{Theorem4-1}(Local existence) If the following condition stands,
\[ \exists M>0,\ \exists \delta>0,\ \exists \delta^{\prime}>0,\ \forall \epsilon\in (0,\ \delta],\ \forall \Omega, m(\Omega)\in (0,\ \delta^{\prime}],\ \mbox{we have}\ \tau(M,\ \epsilon,\ \Omega)>0,\] then\ $\exists \delta^{\prime\prime}\in (0,\ \delta^{\prime}],\ \forall \Omega, m(\Omega)\in (0,\ \delta^{\prime\prime}],\ \forall \epsilon\in (0,\ \delta]$, there exists\ $Z_{1\epsilon}\in C(\overline{\Omega}),\ \|Z_{1\epsilon}\|_{\infty}<M$\, such that\ $Z_{1\epsilon}=T_{0\epsilon}(Z_{1\epsilon})$.\end{theorem}
{\it Proof of theorem 4.1}. If we denote\ $\Omega_{M}=\{Z_{1}\in C(\overline{\Omega}): \|Z_{1}\|_{\infty}<M\}$, then we will see\ $\forall \Omega, m(\Omega)\in (0,\ \delta^{\prime}]$,\ $deg(Z_{1}-T_{0\epsilon}(Z_{1}),\ \Omega_{M},\ 0)$\ keeping constant,\ $\forall \epsilon\in (0,\ \delta]$.\\If we select\ $ \epsilon_{0}\in (0,\ \delta]$, then\ $\forall \epsilon_{1}\in [\epsilon_{0},\ \delta]$, according to (4) in lemma 4.1, we have\ $\exists \delta(\epsilon_{1})>0$, such that$$ \|[T_{0\epsilon_{2}}(Z_{1})-T_{0\epsilon_{1}}(Z_{1})]I_{\overline{\Omega}}\|_{\infty}\leq \tau(M,\ \epsilon_{1},\ \Omega),\ \forall\ Z_{1}\in \Omega_{M},\ \forall  \epsilon_{2}\in U(\epsilon_{1},\ \delta(\epsilon_{1})).  $$
 And from (3) in lemma 4.2, we obtain$$deg(Z_{1}-T_{0\epsilon_{2}}(Z_{1}),\ \Omega_{M},\ 0)=deg(Z_{1}-T_{0\epsilon_{1}}(Z_{1}),\ \Omega_{M},\ 0),\ \forall  \epsilon_{2}\in U(\epsilon_{1},\ \delta(\epsilon_{1})).  $$
Hence\ $ deg(Z_{1}-T_{0\epsilon_{2}}(Z_{1}),\ \Omega_{M},\ 0)$\
keep constant,\ $\forall \epsilon_{2}\in U(\epsilon_{1},\ \delta(\epsilon_{1}))$.\\ We can see that\ $U(\epsilon_{1},\ \delta(\epsilon_{1})),\ \forall \epsilon_{1}\in [\epsilon_{0},\ \delta]$\ is an open cover for\ $[\epsilon_{0},\ \delta]$. From Heine-Borel theorem, there exists finite sub-cover. This means that\ $ deg(Z_{1}-T_{0\epsilon_{1}}(Z_{1}),\ \Omega_{M},\ 0)$\ keep constant,\ $\forall \epsilon_{1}\in [\epsilon_{0},\ \delta]$. From the arbitrary of\ $ \epsilon_{0}$, we know\ $ deg(Z_{1}-T_{0\epsilon}(Z_{1}),\ \Omega_{M},\ 0)$\ keep constant,\ $\forall \epsilon\in (0,\ \delta]$.\\ Now we choose\ $m(\Omega)$\ is sufficient small, such that\ $T_{0\delta}$\ is contract mapping. \\There exists\ $\delta^{\prime\prime}\in (0,\ \delta^{\prime}],\ \forall \Omega, m(\Omega)\in (0,\ \delta^{\prime\prime}]$,\ $T_{0\delta}(\partial\Omega_{M})\subset\overline{\Omega_{M}}$. From Rothe theorem, (2) in lemma 4.2, we get\ $deg(Z_{1}-T_{0\delta}(Z_{1}),\ \Omega_{M},\ 0)\neq 0$. \\Hence we obtain\ $deg(Z_{1}-T_{0\epsilon}(Z_{1}),\ \Omega_{M},\ 0)\neq 0,\ \forall \epsilon\in (0,\ \delta]$.\\
From Kronecker theorem, (1) in lemma 4.2, we know there exist\ $Z_{1\epsilon}\in \Omega_{M}$, such that\ $Z_{1\epsilon}=T_{0\epsilon}(Z_{1\epsilon}),\ \forall \epsilon\in (0,\ \delta],\ \forall \Omega, m(\Omega)\in (0,\ \delta^{\prime\prime}]$.\qed\\
 If (4.36) is not true, then the following will stand,
\[ \forall M>0,\ \exists \epsilon_{k}>0,\ \exists \Omega_{k},\ m(\Omega_{k})>0,\ \lim_{k\rightarrow+\infty}\epsilon_{k}=\lim_{k\rightarrow+\infty}m(\Omega_{k})=0,\ \mbox{such that}\ \tau(M,\ \epsilon_{k},\ \Omega_{k})\equiv 0,\ \forall k\geq1.\]
This means that\ $\forall M>0$, there exists\ $Z_{1}(M,\ \epsilon_{k},\ \Omega_{k})\in C(\overline{\Omega_{k}}),\ k\geq1$, such that
\[ Z_{1}(M,\ \epsilon_{k},\ \Omega_{k})=T_{0\epsilon_{k}}(Z_{1}(M,\ \epsilon_{k},\ \Omega_{k})),\ \|Z_{1}(M,\ \epsilon_{k},\ \Omega_{k})\|_{\infty}=M,\ \forall k\geq1.\]
That's not easy. And if we take\[M_{2}-M_{1}\geq M,\ Z_{1}(M_{1},\ \epsilon_{1},\ \Omega_{1})=Z_{1}(\epsilon_{1},\ \Omega_{1}),\ Z_{1}(M_{2},\ \epsilon_{2},\ \Omega_{2})=Z_{1}(\epsilon_{2},\ \Omega_{2}),\] then we can get as follows,\ $\forall M>0,\ \forall \delta>0,\ \forall \delta^{\prime}>0,\ \exists \epsilon_{1},\ \epsilon_{2}\in (0,\ \delta],\ \exists \Omega_{1},\ \Omega_{2}, m(\Omega_{1}), m(\Omega_{2})\in (0,\ \delta^{\prime}]$, such that\[\|Z_{1}(\epsilon_{1},\ \Omega_{1})-Z_{1}(\epsilon_{2},\ \Omega_{2})\|_{\infty}\geq M,\] where\ $Z_{1}(\epsilon_{j},\ \Omega_{j})=T_{0\epsilon_{j}}(Z_{1}(\epsilon_{j},\ \Omega_{j})),\ j=1,\ 2$. (4.40) looks something related with blow-up is happening.
\begin{corollary}(Global existence) If the following condition stands,
\[ \exists M>0,\ \exists \delta>0,\ \forall \epsilon\in (0,\ \delta],\ \tau(M,\ \epsilon)>0,\ \mbox{and}
\ \exists \epsilon_{0}\in (0,\ \delta],\ deg(Z_{1}-T_{0\epsilon_{0}}(Z_{1}),\ \Omega_{M},\ 0)\neq0,\] then\ $\forall \epsilon\in (0,\ \delta]$, there exists\ $Z_{1\epsilon}\in \Omega_{M}$, such that\ $Z_{1\epsilon}=T_{0\epsilon}(Z_{1\epsilon})$.\end{corollary}
{\it Proof of corollary 4.1}. We can get the proof by Theorem 4.1.\qed\\
We admit that the priori estimation of\ $Z_{1}$\ will helps to obtain the nonzero Leray-Schauder degree in (4.41). However,\ $Z_{1\epsilon}=T_{0\epsilon}(Z_{1\epsilon})$\ is only an ordinary integral equation. It is possible to calculate that degree directly. We will discuss more in Section 5.\\
Next we see the solution of\ $Z_{1}^{\ast}=T_{0}(Z_{1}^{\ast})$\ as the following.
\begin{theorem} \label{Theorem4-2} (1)(Strong solution) $\forall \epsilon_{k}>0,\ \epsilon_{k}\rightarrow0,\ k\rightarrow+\infty,\ Z_{1\epsilon_{k}}\in \Omega_{M}$,\ $Z_{1\epsilon_{k}}=T_{0\epsilon_{k}}(Z_{1\epsilon_{k}}),\ k\geq1$, there exist sub-series\ $n_{k},\ k\geq1$, and\ $Z_{1}^{\ast}\in H^{-m_{1}}(\Omega)$, such that\[ \lim_{k\rightarrow+\infty}\|Z_{1\epsilon_{n_{k}}}-Z_{1}^{\ast}\|_{-m_{1}}=0,\ \mbox{moreover}\ \lim_{k\rightarrow+\infty}\|Z_{1\epsilon_{n_{k}}}-T_{0}(Z_{1\epsilon_{n_{k}}})\|_{-m_{1}}=0.\]If\ $Z_{1}^{\ast}$\ is locally integrable, then\ $\|Z_{1}^{\ast}\|_{L^{\infty}}\leq M$, where\[\|Z_{1}^{\ast}\|_{L^{\infty}}=\max_{1\leq i\leq9m}\|Z_{1,\ i}^{\ast}\|_{L^{\infty}}.\](2)($L^{\infty}$\ solution)If there exist\ $\epsilon_{k}>0,\ Z_{1\epsilon_{k}}\in \Omega_{M}$,\ $Z_{1\epsilon_{k}}=T_{0\epsilon_{k}}(Z_{1\epsilon_{k}}),\ k\geq1$, moreover\[\lim_{k\rightarrow+\infty}\epsilon_{k}=0,\ \lim_{k\rightarrow+\infty}Z_{1\epsilon_{k}}\mbox{exists almost everywhere on}\ \overline{\Omega},\] then there exists\ $Z_{1}^{\ast}\in L^{\infty}(\overline{\Omega}),\ \|Z_{1}^{\ast}\|_{L^{\infty}}\leq M$, such that\ $Z_{1}^{\ast}=T_{0}(Z_{1}^{\ast})$.\end{theorem}
{\it Proof of theorem 4.2}. (1)From (3) in lemma 4.1, we can get that there exist sub-series\ $n_{k},\ k\geq1$, and\ $Z_{1}^{\ast}\in H^{-m_{1}}(\Omega)$, such that\[\lim_{k\rightarrow +\infty}\epsilon_{n_{k}}=0,\ \lim_{k\rightarrow +\infty}\|T_{0}(Z_{1\epsilon_{n_{k}}})-Z_{1}^{\ast}\|_{-m_{1}}=0.\] From\ $Z_{1\epsilon_{n_{k}}}=T_{0\epsilon_{n_{k}}}(Z_{1\epsilon_{n_{k}}}),\ k\geq1$, and (2) in lemma 4.1, we can obtain (4.42) holds.\\
If\ $\|Z_{1}^{\ast}\|_{L^{\infty}}> M$, then there exists\ $Z_{1,\ i}^{\ast}$, such that\ $\|Z_{1,\ i}^{\ast}\|_{L^{\infty}}>M$.\\
Hence there exists\ $\epsilon_{0}>0$, such that\ $m(\Omega_{i}(\epsilon_{0}))>0$, where\[\Omega_{i}(\epsilon_{0})=\{ X\in \overline{\Omega}:\ |Z_{1,\ i}^{\ast}(X)|\geq M+\epsilon_{0}\}.\] Because\ $m(\Omega_{i}(\epsilon_{0}))=m(\Omega_{i}^{+}(\epsilon_{0}))+m(\Omega_{i}^{-}(\epsilon_{0}))$, where
\begin{eqnarray*}\Omega_{i}^{+}(\epsilon_{0})&=&\{ X\in \overline{\Omega}:\ Z_{1,\ i}^{\ast}(X)\geq M+\epsilon_{0}\},\\ \Omega_{i}^{-}(\epsilon_{0})&=&\{ X\in \overline{\Omega}:\ Z_{1,\ i}^{\ast}(X)\leq -(M+\epsilon_{0})\},\end{eqnarray*} we know at least one of\ $m(\Omega_{i}^{+}(\epsilon_{0})),\ m(\Omega_{i}^{-}(\epsilon_{0}))$\ is bigger than\ $0$.\\ We assume as well\ $m(\Omega_{i}^{+}(\epsilon_{0}))>0$. If we choose\ $\varphi_{0}\in C_{0}^{\infty}(\Omega_{i}^{+}(\epsilon_{0}))$, and\ $\varphi_{0}\geq0,\ supp(\varphi_{0})\neq\emptyset$, then we can get that\[\|Z_{1\epsilon_{n_{k}}}-Z_{1}^{\ast}\|_{-m_{1}}\geq \cfrac{\int_{\Omega_{i}^{+}(\epsilon_{0})}\epsilon_{0}\varphi_{0}dX}{\|\varphi_{0}\|_{m_{1}}}>0.\] That's contradict with (4.42).\\
(2)From (4.44), we know there exist\ $Z_{1}^{\ast}$, such that\[ \lim_{k\rightarrow +\infty}Z_{1\epsilon_{k}}=Z_{1}^{\ast},\ \mbox{a.e.}\]
Because\ $Z_{1\epsilon_{k}}\in \Omega_{M}$, we get\ $Z_{1}^{\ast}\in L^{\infty}(\overline{\Omega})$,
moreover\ $\|Z_{1}^{\ast}\|_{L^{\infty}}\leq M.$\\ From\ $\|Z_{1\epsilon_{k}}-Z_{1}^{\ast}\|_{-m_{1}}\leq \|Z_{1\epsilon_{k}}-Z_{1}^{\ast}\|_{L^{2}}$,
where\[\|Z_{1\epsilon_{k}}-Z_{1}^{\ast}\|_{L^{2}}=\max_{1\leq i\leq 9m}\|Z_{1,\ i,\ \epsilon_{k}}-Z_{1,\ i}^{\ast}\|_{L^{2}},\]and Lebesgue dominated convergence theorem, we obtain\[\lim_{k\rightarrow +\infty}\|Z_{1\epsilon_{k}}-Z_{1}^{\ast}\|_{-m_{1}}=0.\]
At last we see\ $\|T_{0}(Z_{1\epsilon_{k}})-T_{0}(Z_{1}^{\ast})\|_{-m_{1}}$. From
\[<T_{0}(Z_{1\epsilon_{k}})-T_{0}(Z_{1}^{\ast}),\ \varphi>=<F[T_{0}(Z_{1\epsilon_{k}})-T_{0}(Z_{1}^{\ast})],\ F^{-1}(\varphi)>,\]
we only need to discuss\ $\psi(Z_{1\epsilon_{k}})-\psi(Z_{1}^{\ast})$\ in\ $F[T_{0}(Z_{1\epsilon_{n_k}})-T_{0}(Z_{1}^{\ast})]$. Because\ $\psi$\ is continuous, we can get there exists a constant\ $C_{M}>0$, such that
\[ \|\psi(Z_{1\epsilon_{k}})-\psi(Z_{1}^{\ast})\|_{L^{\infty}}\leq C_{M}.\]
Again from Lebesgue dominated convergence theorem, we obtain
\[\lim_{k\rightarrow +\infty}\|T_{0}(Z_{1\epsilon_{k}})-T_{0}(Z_{1}^{\ast})\|_{-m_{1}}=0.\]Together with\ $Z_{1\epsilon_{k}}=T_{0\epsilon_{k}}(Z_{1\epsilon_{k}})$, we obtain\ $\|Z_{1}^{\ast}-T_{0}(Z_{1}^{\ast})\|_{-m_{1}}=0$. This means\ $Z_{1}^{\ast}=T_{0}(Z_{1}^{\ast}),\ \mbox{a.e.}$\ which completed the statement.\qed\\
From the previous theorem and corollary, we know that the strong solution will exist locally under the condition (4.36) and exist globally under the condition (4.41).\\
Finally, we discuss a little more for\ $L^{\infty}$\ solution as follows.
\begin{theorem} \label{Theorem4-3} A necessary and sufficient condition for (4.44) holds is that there exist\ $\eta_{k}>0,\ Z_{1\eta_{k}}\in \Omega_{M}$,\ $Z_{1\eta_{k}}=T_{0\eta_{k}}(Z_{1\eta_{k}}),\ k\geq1$, moreover\[\lim_{k\rightarrow+\infty}\eta_{k}=0,\ \lim_{k,\ l\rightarrow+\infty}\|Z_{1\eta_{k}}-Z_{1\eta_{l}}\|_{L^{2}}=0.\]\end{theorem}
{\it Proof of theorem 4.3}. Necessity. If (4.44) holds, then there exists\ $Z_{1}^{\ast}\in L^{\infty}(\overline{\Omega})$, such that
\[\lim_{k\rightarrow +\infty}\|Z_{1\epsilon_{k}}-Z_{1}^{\ast}\|_{L^{2}}=0.\] If we let\ $\eta_{k}=\epsilon_{k},\ k\geq1$, then (4.54) stands.\\
Sufficiency. If (4.54) holds, then\ $Z_{1\eta_{k}},\ k\geq 1$, are convergent by Lebesgue measure as follows,
\[ \forall \epsilon>0,\ \lim_{k,\ l\rightarrow +\infty}m(\{X\in \overline{\Omega}: \|Z_{1\eta_{k}}-Z_{1\eta_{l}}\|_{\infty}\geq \epsilon\})=0.\]
From Riesz theorem on page 142 in [14], we will see that there exist sub-series\ $n_{k},\ k\geq1$, such that
\[ \lim_{k\rightarrow+\infty}\eta_{n_{k}}=0,\ \lim_{k\rightarrow+\infty}Z_{1\eta_{n_{k}}}\mbox{exists almost everywhere on}\ \overline{\Omega}.\]
If we let\ $\epsilon_{k}=\eta_{n_{k}}$, then (4.44) stands.\qed\\
We know\ $L^{2}(\Omega)$\ is not completed under norm\ $\|\cdot\|_{-m_{1}}$. After being completed, it will be\ $H^{-m_{1}}(\Omega)$. So there exist\ $f_{k}\in L^{2}(\Omega),\ k\geq1$, such that\[ \lim_{k,\ l\rightarrow +\infty}\|f_{k}-f_{l}\|_{-m_{1}}=0,\ \mbox{and}\ \|f_{k}-f_{l}\|_{L^{2}}\geq c>0,\ \forall\ k\neq l.\]
From (4.54), we know that (4.44) will not always stand. But we can see that it will hold in many cases as the following.
\begin{theorem} \label{Theorem4-4}A sufficient condition for (4.44) holds is that there exist\ $\epsilon_{k}>0,\ \epsilon_{k}\rightarrow0,\ k\rightarrow+\infty,\ Z_{1\epsilon_{k}}\in \Omega_{M}$,\ $Z_{1\epsilon_{k}}=T_{0\epsilon_{k}}(Z_{1\epsilon_{k}}),\ k\geq1$, such that\ $\forall c>0,\ \exists d>0,\ \forall k,\ l,\ i,\ 1\leq i\leq 9m$, we have\[m(\Omega_{k,\ l,\ i}^{+}(d)\cup \Omega_{k,\ l,\ i}^{-}(d))\leq c,\]
where\begin{eqnarray*}\Omega_{k,\ l,\ i}^{+}&=&\{X\in \overline{\Omega}:Z_{1,\ i,\ \epsilon_{k}}-Z_{1,\ i,\ \epsilon_{l}}\geq 0\},\\
\Omega_{k,\ l,\ i}^{-}&=&\{X\in \overline{\Omega}:Z_{1,\ i,\ \epsilon_{k}}-Z_{1,\ i,\ \epsilon_{l}}< 0\},\\ \Omega_{k,\ l,\ i}^{+}(d)&=&\{X\in \overline{\Omega}:\ dist(X,\ \partial \Omega_{k,\ l,\ i}^{+})\leq d\},\\ \Omega_{k,\ l,\ i}^{-}(d)&=&\{X\in \overline{\Omega}:\ dist(X,\ \partial \Omega_{k,\ l,\ i}^{-})\leq d\}.\end{eqnarray*}\end{theorem}
{\it Proof of theorem 4.4}. From Theorem 4.2, we know there exist sub-series\ $n_{k},\ k\geq1$, such that
\[\lim_{k,\ l\rightarrow +\infty}\|Z_{1\epsilon_{n_{k}}}-Z_{1\epsilon_{n_{l}}}\|_{-m_{1}}=0.\]
We will prove that\ $\forall i,\ 1\leq i\leq 9m,\ Z_{1,\ i,\ \epsilon_{n_{k}}},\ k\geq1$, is convergent by Lebesgue measure.\\
If that is not true, then there exist\ $i_{0},\ 1\leq i_{0}\leq 9m,\ \exists b_{1}>0,\ \exists b_{2}>0$, and sub-series\ $k_{j},\ l_{j},\ j\geq1$, such that
\[ m(\{X\in \overline{\Omega}:\ |Z_{1,\ i_{0},\ \epsilon_{n_{k_{j}}}}-Z_{1,\ i_{0},\ \epsilon_{n_{l_{j}}}}|\geq b_{1}\})\geq b_{2},\ j\geq1.\]
We denote it in an easy way as follows,\ $\forall j,\ j\geq1$,\begin{eqnarray*}f_{j}&=&Z_{1,\ i_{0},\ \epsilon_{n_{k_{j}}}}-Z_{1,\ i_{0},\ \epsilon_{n_{l_{j}}}},\\
\Omega_{j}^{+}&=&\{X\in \overline{\Omega}:f_{j}\geq 0\},\\
\Omega_{j}^{-}&=&\{X\in \overline{\Omega}:f_{j}< 0\},\\ \Omega_{j}^{+}(d)&=&\{X\in \overline{\Omega}:\ dist(X,\ \partial \Omega_{j}^{+})\leq d\},\\ \Omega_{j}^{-}(d)&=&\{X\in \overline{\Omega}:\ dist(X,\ \partial \Omega_{j}^{-})\leq d\}.\end{eqnarray*}
From (4.59), we can obtain that there exists\ $d_{0}>0$, such that
\[ m(\Omega_{j}^{+}(d_{0})\cup\Omega_{j}^{-}(d_{0}))\leq \cfrac{b_{1}b_{2}}{4M}.\]
If we let\ $\varphi_{j}\in C_{0}^{\infty}(\Omega),\ j\geq1$, as the following,
\[ \varphi_{j}(X)=\int_{\Omega_{j}^{+}\setminus \Omega_{j}^{+}(d_{0})} \alpha_{d_{0}}(X-Y)dY-\int_{\Omega_{j}^{-}\setminus \Omega_{j}^{-}(d_{0})} \alpha_{d_{0}}(X-Y)dY,\ j\geq1,\]
where$$\alpha_{d_{0}}(X)=\cfrac{1}{d_{0}^{3}}\alpha(\cfrac{X}{d_{0}}),\ \alpha(X)= \begin{cases}
                                                                                                                              Ce^{1/(|X|^{2}-1)},\ |X|<1, \\
                                                                                                                              0,\ |X|\geq 1.
                                                                                                                            \end{cases} ,\ C=(\int_{|X|<1}e^{1/(|X|^{2}-1)}dX)^{-1},$$
then\ $\varphi_{j}=1$\ on\ $\Omega_{j}^{+}\setminus \Omega_{j}^{+}(d_{0})$,\ $\varphi_{j}\in [0,\ 1]$\ on\ $\Omega_{j}^{+}(d_{0})$, and\ $\varphi_{j}=-1$\ on\ $\Omega_{j}^{-}\setminus \Omega_{j}^{-}(d_{0})$,\ $\varphi_{j}\in [-1,\ 0]$\ on\ $\Omega_{j}^{-}(d_{0})$.
\\Moreover we can get as follows, \[ |\partial^{\gamma}\varphi_{j}|\leq 2\int_{\Omega}|\partial^{\gamma}\alpha_{d_{0}}(X-Y)|dY,\ 0\leq |\gamma|\leq m_{1}.\]
Hence there exists\ $d_{1}>0$, such that\[ \|\varphi_{j}\|_{m_{1}}\leq d_{1},\ \forall j,\ j\geq1.\]
From $$ |<f_{j},\ \varphi_{j}>|\geq|<f_{j},\ sign f_{j}>|-|<f_{j},\ sign f_{j}-\varphi_{j}>|\geq b_{1}b_{2}-2M\cfrac{b_{1}b_{2}}{4M}=\cfrac{b_{1}b_{2}}{2},$$where
$$ sign f_{j}=\begin{cases}
1,\ f_{j}\geq0,\\
-1,\ f_{j}<0. \end{cases},$$ we can obtain that\[ \|f_{j}\|_{-m_{1}}\geq \cfrac{|<f_{j},\ \varphi_{j}>|}{\|\varphi_{j}\|_{m_{1}}}\geq \cfrac{b_{1}b_{2}}{2d_{1}}.\]
That's contradict with (4.60).\\
So\ $\forall i,\ 1\leq i\leq 9m,\ Z_{1,\ i,\ \epsilon_{n_{k}}},\ k\geq1$, is convergent by Lebesgue measure. From Riesz theorem on page 142 in [7], we know that (4.44) holds. \qed\\ If (4.59) is not true, then\ $\forall \epsilon_{k}>0,\ \epsilon_{k}\rightarrow0,\ k\rightarrow+\infty,\ Z_{1\epsilon_{k}}\in \Omega_{M}$,\ $Z_{1\epsilon_{k}}=T_{0\epsilon_{k}}(Z_{1\epsilon_{k}}),\ k\geq1$,\ $\exists c>0,\ \forall d_{j}>0,\ \exists k_{j},\ l_{j},\ i_{j},\ 1\leq i_{j}\leq 9m$, such that\[m(\Omega_{k_{j},\ l_{j},\ i_{j}}^{+}(d_{j})\cup \Omega_{k_{j},\ l_{j},\ i_{j}}^{-}(d_{j}))\geq c.\]
If we let\ $ d_{j}\rightarrow 0,\ j\rightarrow +\infty$, then we can get the following,\[ \lim_{j\rightarrow +\infty} S(\partial\Omega_{k_{j},\ l_{j},\ i_{j}}^{+})=+\infty,\] where $$S(\partial\Omega_{k_{j},\ l_{j},\ i_{j}}^{+})=\int_{\{X\in \overline{\Omega}:\ Z_{1,\ i_{j},\ \epsilon_{k_{j}}}-Z_{1,\ i_{j},\ \epsilon_{l_{j}}}= 0\}}dS. $$ And the following is another sufficient condition for (4.44) holds,
\[ \exists\ \epsilon_{k}>0,\ Z_{1\epsilon_{k}}\in \Omega_{M},\ Z_{1\epsilon_{k}}=T_{0\epsilon_{k}}(Z_{1\epsilon_{k}}),\ k\geq1,\ \lim_{k\rightarrow +\infty}\epsilon_{k}=0,\ \mbox{and}\ \sup_{k,\ l,\ i}S(\partial\Omega_{k,\ l,\ i}^{+})< +\infty,\]where $$S(\partial\Omega_{k,\ l,\ i}^{+})=\int_{\{X\in \overline{\Omega}:\ Z_{1,\ i,\ \epsilon_{k}}-Z_{1,\ i,\ \epsilon_{l}}= 0\}}dS,\ k\neq l,\ 1\leq i\leq 9m. $$
We obtain\ $u^{\ast}\in W^{2,\ +\infty}(\overline{\Omega})$, if\ $Z_{1}^{\ast}=T_{0}(Z_{1}^{\ast})$,\ $Z_{1}^{\ast}\in L^{\infty}(\overline{\Omega})$, where\ $Z_{1}^{\ast}=(v_{m+1}^{\ast},\ \cdots,\ v_{10m}^{\ast})^{T}$. Here\ $\ W^{2,\ +\infty}(\overline{\Omega})$\ is a Sobolev space defined on page 153 in [2]. From the condition that domain\ $\Omega$\ satisfies a uniform exterior and interior cone, if\ $\Omega$\ is bounded,\ $\partial\Omega\in C^{1,\ \beta},\ 0<\beta\leq 1$, we can get that\ $u^{\ast}\in C^{1,\ 1}(\overline{\Omega})$\ by imbedding. By using Morrey's inequality defined on page 163 in [2], we get\ $u^{\ast}$\ is twice classically differentiable almost everywhere in\ $\overline{\Omega}$.\\
\section{Leray-Schauder degree}\setcounter{equation}{0}
In this section, we will discuss the
the Leray-Schauder degree of nonlinear integral equation of Hammerstein type as follows,
\[ f(X)=g(X)+T(f(X)),\ \forall X\in \overline{\Omega},\]
where\[T(f(X))=\int_{\Omega}k(X,\ Y)\psi(Y,\ f(Y))dY,\]
$X=(x,\ y,\ z)^{T},\ Y=(x_{1},\ y_{1},\ z_{1})^{T}$,\ $f(X)$\ is an unknown continuous function on\ $\overline{\Omega}$,\ $g(X)$\ is a known continuous function on\ $\overline{\Omega}$,\ $k(X,\ Y)$\ is a known continuous function on\ $\overline{\Omega}\times\overline{\Omega}$,\ $\psi$\ is a known continuous function on\ $\overline{\Omega}\times[-M,\ M]$.\\
From the definition of the Leray-Schauder degree on page 138 to page 139 in [8], we can work out the Leray-Schauder degree of Eq(5.1) directly as follows.\\
By the Weirstrass theorem, we know there exist two polynomials\[k_{N}(X,\ Y)=\sum_{|\alpha|=0}^{N}C_{\alpha}(Y)X^{\alpha},\ g_{N}(X)=\sum_{|\alpha|=0}^{N}g_{\alpha}X^{\alpha},\] where\ $\alpha=(\alpha_{1},\ \alpha_{2},\ \alpha_{3})^{T},\ |\alpha|=\alpha_{1}+\alpha_{2}+\alpha_{3},\ X^{\alpha}=x^{\alpha_{1}}y^{\alpha_{2}}z^{\alpha_{3}},\ C_{\alpha}(Y),\ 0\leq|\alpha|\leq N$, are all polynomials of\ $Y$, moreover\ $\alpha_{1},\ \alpha_{2},\ \alpha_{3}$\ are all nonnegative whole numbers,\ $g_{\alpha},\ 0\leq|\alpha|\leq N$,
are all real numbers, such that\ $\forall\ f(X)\in\Omega_{M}=\{f(X):\|f(X)\|_{\infty}< M\}$, we have
\[ \|\int_{\Omega}(k(X,\ Y)-k_{N}(X,\ Y))\psi(Y,\ f(Y))dY\|_{\infty}\leq \cfrac{\tau}{3},\ \|g(X)-g_{N}(X)\|_{\infty}\leq \cfrac{\tau}{3},\]
$\tau$\ is defined as follows,\[ \tau=\inf_{\|f(X)\|_{\infty}= M}\|f(X)-T(f(X))-g(X)\|_{\infty}>0.\]
From the Combination theory, the number of the solutions of\ $|\alpha|=\alpha_{1}+\alpha_{2}+\alpha_{3}$\ is that\ $C_{2+|\alpha|}^{2}$. The number of all the items\ $X_{\alpha},\ 0\leq|\alpha|\leq N$, is that\ $L_{N}=C_{2}^{2}+C_{3}^{2}+\cdots+C_{2+N}^{2}=C_{3+N}^{3}$.\\By the homotopic, we have the following,
\[ deg(f(X)-T(f(X)),\ \Omega_{M},\ g(X))=deg(f(X)-T(f(X)),\ \Omega_{M},\ g_{N}(X)).\]
Assuming that\[T_{N}(f(X))=\int_{\Omega}k_{N}(X,\ Y)\psi(Y,\ f(Y))dY,\] then from (5.4) we have\[ \|T(f(X))-T_{N}(f(X))\|_{\infty}\leq \cfrac{\tau}{3},\ \forall\ f(X)\in\Omega_{M}.\]
Moreover,\[T_{N}(f(X))=\sum_{|\alpha|=0}^{N}(\int_{\Omega}C_{\alpha}(Y)\psi(Y,\ f(Y))dY)X^{\alpha}\in E_{N},\]
where\ $E_{N}$\ is the sub-space with finite dimensions generated by\ $X^{\alpha},\ 0\leq|\alpha|\leq N$.\\From the definition of the Leray-Schauder degree, we can obtain the following,
\[deg(f(X)-T(f(X)),\ \Omega_{M},\ g_{N}(X))=deg(f(X)-T_{N}(f(X)),\ \Omega_{M,\ 1},\ g_{N}(X)),\]where\ $\Omega_{M,\ 1}=\Omega_{M}\cap E_{N}$.\\
If we denote\begin{eqnarray} f(X)&=&\tilde{X}^{T}D_{N},\ \tilde{X}=(X^{\alpha},\ 0\leq|\alpha|\leq N)^{T},\\
D_{N}&=&(D_{\alpha},\ 0\leq|\alpha|\leq N)^{T}\in R^{L_{N}},\end{eqnarray} \begin{eqnarray}
T_{N}(f(X))&=&\int_{\Omega}k_{N}(X,\ Y)\psi(Y,\ f(Y))dY\\&=&\sum_{|\alpha|=0}^{N}(\int_{\Omega}C_{\alpha}(Y)\psi(Y,\ \tilde{Y}^{T}D_{N})dY)X^{\alpha},\\&=&\sum_{|\alpha|=0}^{N}\phi_{\alpha}(D_{N})X^{\alpha}=\tilde{X}^{T}\phi(D_{N}),\end{eqnarray} \begin{eqnarray}
\phi_{\alpha}(D_{N})&=&\int_{\Omega}C_{\alpha}(Y)\psi(Y,\ \tilde{Y}^{T}D_{N})dY,\\
\phi(D_{N})&=&(\phi_{\alpha}(D_{N}),\ 0\leq|\alpha|\leq N)^{T},\\ g_{N}&=&(g_{\alpha},\ 0\leq|\alpha|\leq N)^{T},\\
\Omega_{M,\ 2}&=&\{D_{N}\in R^{L_{N}}: \|\tilde{X}^{T}D_{N}\|_{\infty}< M\},
\end{eqnarray} then we obtain\ $ f(X)=T_{N}(f(X))+g_{N}(X),\ f(X)\in \Omega_{M,\ 1}$, is equivalent to
\[ D_{N}=\phi(D_{N})+g_{N},\ D_{N}\in \Omega_{M,\ 2}.\]Hence, we obtain
\[ deg(f(X)-T(f(X)),\ \Omega_{M},\ g_{N}(X))=deg(D_{N}-\phi(D_{N}),\ \Omega_{M,\ 2},\ g_{N}).\]
From (5.16), we can see that\ $\phi(D_{N})$\ is explicit. At this time, we can work out the finite dimensions Brouwer degree\ $deg(D_{N}-\phi(D_{N}),\ \Omega_{M,\ 2},\ g_{N})$\ directly. From the definition of the Brouwer degree on page 89 in [8], we only need to calculate the double integral on\ $\Omega_{M,\ 2}\times\Omega$. This is just what we want, the Leray-Schauder degree of Eq(5.1).\\As\ $\tau$\ is fixed, we don't need that\ $N$\ is sufficient big. This is the reason why we haven't discuss the priori estimation yet.\\By the same way, we can also work out the Leray-Schauder degree of Eq(5.1) even if
\[ T(f(X))=\int_{\Omega}G(X,\ Y,\ f(Y))dY,\] where\ $G$\ is a known continuous function on\ $\overline{\Omega}\times\overline{\Omega}\times[-M,\ M]$.\\
We only need to change\ $k_{N}(X,\ Y)$\ into the following,
\[ k_{N}(X,\ Y,\ f(Y))=\sum_{|\alpha|=0}^{N}C_{\alpha,\ 1}(Y,\ f(Y))X^{\alpha},\]where\ $C_{\alpha,\ 1}(Y,\ f(Y)),\ 0\leq|\alpha|\leq N$, are all polynomials of\ $Y,\ f(Y)$.\\
 \section{Declarations}
Availability of supporting data\\
This paper is available to all supporting data.\\
Competing interests\\
We declare that we have no competing interests.\\
Funding\\
The author was supported by mathematics research fund of Hohai university, No. 1014-414126.\\
Authors' contributions\\
We have transformed nonlinear second order partial differential equations resolved with any derivatives into the equivalent generalized integral equations. Moreover, we discuss the existence for the classical solution by Leray-Schauder degree and Sobolev space.\\
Acknowledgements\\
We would like to give our best thanks to Prof. Mark Edelman in Yeshiva University, Prof. Caisheng Chen in Hohai University, Prof. Junxiang Xu in Southeast University, Prof. Zuodong Yang in Nanjing Normal University, Prof. Jishan Fan in Nanjing Forestry University for their guidance and other helps.\\
Authors' information\\ Jianfeng Wang, Prof. in the Mathematics department of Hohai University, email: wjf19702014@163.com.\\Acadmic email: 20020001@hhu.edu.cn\\
Orcid: 0000-0002-3129-756x\\


\begin{thebibliography}{2}
\bibitem{1}L. H$\ddot{o}$rmander, The analysis of linear partial differential operators I, Springer-Verlag, 1983, page 181 to page 183.
\bibitem{2}Gilbarg D and Trudinger N S. Elliptic partial differential equations of second order [M]. 2nd ED. Berlin: Springer-Verlag, 1983, page 144 to page 177.
\bibitem{3}Yazhe, Chen, Elliptic partial differential equation and Elliptic partial differential equations of second order, Science education press, 1991.
\bibitem{4}Lawrence C. Evans, Partial Differential Equations, Graduate Studies in Mathematics, Volume 19, American Mathematical Society, 1997.
\bibitem{5}Gongqing Zhang, Yuanqu Lin, Functional analysis, Peking
university press, 2010, page 195.
\bibitem{6}Chaohao Gu, Zhengfan Xu, Daqian Li, Zongyi Hou, Likang Li, The mathematical-physics equations, (second edition), Science and Technology press of Shanghai, 1978, page 173 to page 201.
\bibitem{7}Mingqiang Zhou, Theory on real function, Peking
university press, 2000, page 136.
\bibitem{8}Dajun Guo, Nonlinear Functional Analysis, (second edition), Science and Technology press of Shandong, 2001, page 135 to page 165.
\bibitem{9}Daqian Li, Tiehu
Qing, Physics and partial differential equations, (second edition),
Higher education press, 2005, page 112 to page 115.
\bibitem{10}Shuxing Chen, The general theory of the partial
differential equation, Higher education press, 1981, page 118 to page 119, page 130.
\bibitem{11}Caisheng Chen, Gang Li, Jidong Chou, Wenchu Wang,
The mathematical-physics equations, Science education press, 2008, page 219.
\bibitem{12}Tian Ma, The theory and method of partial differential equations, Science education press, 2011.
\end{thebibliography}
\end{document}